\documentclass[]{article}

\usepackage[utf8]{inputenc}
\usepackage[T1]{fontenc}
\usepackage{geometry}
\usepackage[hang, small,labelfont=bf,up,textfont=it,up]{caption} 
\usepackage{booktabs} 
\usepackage{abstract} 

\usepackage{titlesec}
\titleformat{\section}[block]{\large\scshape\centering}{\thesection.}{1em}{} 
\titleformat{\subsection}[block]{\large}{\thesubsection.}{1em}{} 

\usepackage{amssymb, epsfig,amssymb, latexsym}
\usepackage{bm}
\usepackage{amsfonts,psfrag,amsmath,bbm,color,url}
\usepackage{graphicx}
\graphicspath{ {./img/} }
\usepackage{subfig}

\usepackage{amsfonts}
 \usepackage{amsmath}
 \usepackage{xcolor}
\usepackage{tikz}
\usepackage{pgfplots}
\usepackage{comment}
\usepackage{bm}
\usepackage{mathabx}
\DeclareMathAlphabet{\mathbbx}{U}{bbold}{m}{n}

\usepackage{cases}
\usepackage[output-decimal-marker={,}]{siunitx}
\usepackage{lineno,hyperref}
\definecolor{vargreen}{rgb}{0.0, 0.5, 0.0}
\definecolor{navyblue}{rgb}{0.0, 0.0, 0.5}

\newcommand{\blue}[1]{{\color{black} #1}}
\definecolor{mediumorchid}{rgb}{0.73, 0.33, 0.83}
\definecolor{crimson}{rgb}{0.86, 0.08, 0.24}
\definecolor{lightseagreen}{rgb}{0.13, 0.7, 0.67}
\definecolor{royalblue}{rgb}{0.25, 0.41, 0.88}
\definecolor{hotpink}{rgb}{1.0, 0.41, 0.71}
\definecolor{magenta}{rgb}{1.0, 0.0, 1.0}
\definecolor{goldenrod}{rgb}{0.85, 0.65, 0.13}

\definecolor{plum(traditional)}{rgb}{0.56, 0.27, 0.52}
\usepackage{appendix}

\usepackage[export]{adjustbox}
\usepackage{mdframed}
\usepackage{array}
\usepackage{tikz}
\usetikzlibrary{shapes}
\usetikzlibrary{plotmarks}
\pgfplotsset{compat=1.17}
\usepackage[mode=buildnew]{standalone}
\usepackage[draft,inline,marginclue]{fixme}
\usepackage{xpatch}
\usepackage{nomencl} 
\makenomenclature
\RequirePackage{ifthen}
\let\oldref\ref
\renewcommand{\ref}[1]{(\oldref{#1})}

\DeclareMathAlphabet\mathbfcal{OMS}{cmsy}{b}{n}

\xpatchcmd{\thenomenclature}{%
  \section{\nomname}
}{}{\typeout{Success}}{\typeout{Failure}}

\renewcommand{\nomname}{B.$\quad$List of abbreviations and symbols}

\usepackage{ifthen}
\renewcommand{\nomgroup}[1]{%
  \ifthenelse{\equal{#1}{A}}{\item[\textbf{Abbreviations}]}{%
    \ifthenelse{\equal{#1}{G}}{\item[\textbf{Symbols}]}{%
      \ifthenelse{\equal{#1}{C}}{\item[\textbf{Abbreviations}]}{%
        \ifthenelse{\equal{#1}{S}}{\item[\textbf{Subscripts}]}{%
          \ifthenelse{\equal{#1}{Z}}{\item[\textbf{Mathematical Symbols}]}{}
        }
      }
    }
  }
}

\FXRegisterAuthor{gs}{ags}{\color{red}GS}

\newcommand{\Toffline}{\ensuremath{\text{T}_{\text{offline}}}}
\newcommand{\Tonline}{\ensuremath{\text{T}_{\text{online}}}}

\newcommand{\nutrain}{\ensuremath{\mathcal{N}_{\text{train}}}}
\newcommand{\nutest}{\ensuremath{\mathcal{N}_{\text{test}}}}

\usepackage{siunitx}
\usepackage{multirow}
\usepackage{comment}
\usepackage{pifont}
\newcommand{\xmark}{\ding{55}}%

\pgfplotsset{compat=1.17}

\title{\textbf{Parametric Intrusive Reduced Order Models enhanced with Machine Learning Correction Terms}}

\date{ }

\author{Anna Ivagnes  \\ \small SISSA, International School for Advanced Studies, \\ \small Mathematics Area, mathLab, Trieste, Italy. \\ \small  \href{mailto:aivagnes@sissa.it}{aivagnes@sissa.it} \normalsize \and Giovanni Stabile \\ \small Sant'Anna School of Advanced Studies
\\ \small  The Biorobotics Institute, 
Pontedera, Pisa, Italy. \\ \small  \href{mailto:giovanni.stabile@santannapisa.it}{giovanni.stabile@santannapisa.it} \normalsize \and Gianluigi Rozza  \\ \small  SISSA, International School for Advanced Studies, \\\small  Mathematics Area, mathLab, Trieste, Italy. \\ \small \href{mailto:grozza@sissa.it}{grozza@sissa.it}}

\begin{document}

\maketitle

\nomenclature{$\text{FOM}$}{Full Order Model}
\nomenclature{$\text{ROM}$}{Reduced Order Model}
\nomenclature{$\text{NSE}$}{Navier--Stokes Equations}
\nomenclature{$\text{POD}$}{Proper Orthogonal Decomposition}
\nomenclature{$\text{PCA}$}{Principal Component Analysis}
\nomenclature{$\text{SVD}$}{Singular Value Decomposition}
\nomenclature{$\text{PPE-ROM}$}{ROM with pressure Poisson equation stabilisation}
\nomenclature{$\text{DD-ROM}$}{Data-driven Reduced Order Model}
\nomenclature{$\text{URANS}$}{Unsteady Reynolds--Averaged Navier--Stokes}
\nomenclature{$\text{MSE}$}{Mean Squared Error}
\nomenclature{$\text{NN}$}{Neural Network}
\nomenclature{$\text{LSTM}$}{Long--Short Term Memory neural network}
\nomenclature{$\text{SinNN}$}{Neural network with ad-hoc architecture (specified in Figure \ref{fig:sinNN})}
\nomenclature{$\text{LES}$}{Large Eddy Simulation}
\nomenclature{$\text{EFR}$}{Evolve--Filter--Relax}
\nomenclature{$\text{SST}$}{Shear Stress Transport}
\nomenclature{$\text{ODE}$}{Ordinary Differential Equation}

\nomenclature[G]{$\bm{u}$}{velocity field}
\nomenclature[G]{${p}$}{pressure field}
\nomenclature[G]{$\bm{\mu}$}{parameters considered for the ROM}
\nomenclature[G]{$\nu_t$}{eddy viscosity field}
\nomenclature[G]{$\bm{u}_r$}{reduced velocity field}
\nomenclature[G]{$p_r$}{reduced pressure field}
\nomenclature[G]{$\nu_{t_r}$}{reduced eddy viscosity field}

\nomenclature[G]{$\bm{a}$}{reduced vector of unknowns for velocity}
\nomenclature[G]{$\bm{b}$}{reduced vector of unknowns for pressure}
\nomenclature[G]{$\bm{g}$}{reduced vector of unknowns for eddy viscosity}
\nomenclature[G]{$\bm{a}_r^{proj}$}{projected reduced velocity coefficients}
\nomenclature[G]{$\bm{b}_q^{proj}$}{projected reduced pressure coefficients}

\nomenclature[G]{$Re$}{Reynolds Number}
\nomenclature[G]{$\mathcal{R}_{ij}$}{$(ij)_{th}$ component of the Reynolds stress tensor}
\nomenclature[G]{$\bm{E}_{ij}$}{$(ij)_{th}$ component of the averaged strain rate tensor}
\nomenclature[G]{${\nu}$}{kinematic viscosity}
\nomenclature[G]{${N_u^h}$}{number of unknowns for velocity at full-order level}
\nomenclature[G]{${N_p^h}$}{number of unknowns for pressure at full-order level}
\nomenclature[G]{$N_u$ or $r$}{number of unknowns for velocity at reduced order level}
\nomenclature[G]{$N_p$ or $q$}{number of unknowns for pressure at reduced order level}
\nomenclature[G]{$d$}{number of velocity modes re-introduced in the ROM with correction terms}
\nomenclature[G]{$h$}{number of pressure modes re-introduced in the ROM with correction terms}
\nomenclature[G]{$\bm{\tau}_u$}{correction term for the momentum equation}
\nomenclature[G]{$\bm{\tau}_p$}{correction term for the PPE}
\nomenclature[G]{$\Toffline$}{training/offline set of time instances}
\nomenclature[G]{$\Tonline$}{online set of time instances}
\nomenclature[G]{$\nutrain$}{training/offline set of viscosity values}
\nomenclature[G]{$\nutest$}{online set of viscosity values}
\nomenclature[G]{$N_{\text{networks}}$}{number of neural networks used to compute the confidence interval}
\nomenclature[G]{$\varepsilon_{\bm{u}}^{\nu^{\star}}(t)$}{relative error of velocity field with respect to the full order counterpart, for the parameter $\nu^{\star}$}
\nomenclature[G]{$\varepsilon_{p}^{\nu^{\star}}(t)$}{relative error of pressure field with respect to the full order counterpart, for the parameter $\nu^{\star}$}

\nomenclature[G]{$N_P$}{number of snapshots}
\nomenclature[G]{$N_T$}{number of time instances}
\nomenclature[G]{$N_M$}{number of viscosity values}
\nomenclature[G]{$\mathbb{V}^u_{\text{POD}}$}{reduced basis space for velocity}
\nomenclature[G]{$\mathbb{V}^p_{\text{POD}}$}{reduced basis space for pressure}
\nomenclature[G]{${\Omega}$}{bounded domain}
\nomenclature[G]{${\Gamma}$}{boundary of $\Omega$}
\nomenclature[G]{${\Gamma_{D_i}}$}{i-th boundary of $\Omega$, with non-homogeneous Dirichlet boundary condition}
\nomenclature[G]{$N_{\text{BC}}$}{number of boundaries of $\Omega$ with Dirichlet non-homogeneous boundary conditions}
\nomenclature[G]{$\bm{n}$}{outward normal vector}
\nomenclature[G]{$\bm{\varphi_i}$}{i-th POD basis function for velocity}
\nomenclature[G]{${\chi_i}$}{i-th POD basis function for pressure}
\nomenclature[G]{$\eta_{i}$}{i-th POD basis function for eddy viscosity}

\nomenclature[G]{$\mathcal{G}(\cdot)$}{mapping used to compute the reduced eddy viscosity field}
\nomenclature[G]{$\mathcal{M}(\cdot)$}{mapping used to compute the correction terms}

\nomenclature[G]{$\bm{{\mathcal{S}_u}}$}{snapshots matrix for the velocity field}
\nomenclature[G]{$\bm{{\mathcal{S}_p}}$}{snapshots matrix for the pressure field}
\nomenclature[G]{$\bm{f}(\cdot;\bm{\mu})$}{residual of the reduced momentum conservation equation}
\nomenclature[G]{$\bm{c}(\cdot; \bm{\mu})$}{residual of the reduced continuity equation}
\nomenclature[G]{$\bm{h}(\cdot; \bm{\mu})$}{residual of the reduced PPE}
\nomenclature[G]{$\bm{\tau}^{exact}$}{exact correction term}
\nomenclature[G]{$\bm{\tau}^{approx}$}{approximated correction term}
\nomenclature[G]{$\tau$}{weight used for boundary penalty term}
\nomenclature[G]{$N_{\text{seq}}$}{sequence length in the LSTM neural network}

\nomenclature[G]{$\otimes$}{tensor product}
\nomenclature[G]{$\bm{\nabla}\cdot$}{divergence operator}
\nomenclature[G]{$\bm{\nabla}\times$}{curl operator}
\nomenclature[G]{$\bm{\nabla}$}{gradient operator}
\nomenclature[G]{$\Delta$}{laplacian operator}
\nomenclature[G]{$\left\lVert \cdot\right\rVert_{L^2(\Omega)}$}{norm in $L^2(\Omega)$}
\nomenclature[G]{$( \cdot , \cdot )_{L^2(\Omega)}$}{inner product in $L^2(\Omega)$}

\begin{abstract}
\noindent In this paper, we propose an equation-based parametric Reduced Order Model (ROM), whose accuracy is improved with data-driven terms added into the reduced equations.
These additions have the aim of reintroducing contributions that in standard ROMs are not taken into account. In particular, in this work we consider two types of contributions: the \textbf{turbulence modeling}, added through a reduced-order approximation of the eddy viscosity field, and the \textbf{correction model}, aimed to re-introduce the contribution of the discarded modes.
Both approaches have been investigated in previous works such as~\cite{mohebujjaman2019physically, hijazi2020data, ivagnes2023pressure, ivagnes2023hybrid} and the goal of this paper is to extend the model to a \textbf{parametric} setting making use of ad-hoc machine learning procedures.
More in detail, we investigate different neural networks' architectures, from simple dense feed-forward to Long-Short Term Memory neural networks, in order to find the most suitable model for the re-introduced contributions.
We tested the methods on two test cases with different behaviors: the periodic turbulent flow past a circular cylinder and the unsteady turbulent flow in a channel-driven cavity. In both cases, the parameter considered is the Reynolds number and the machine learning-enhanced ROM considerably improved the pressure and velocity accuracy with respect to the standard ROM.
\emph{A visual version of the abstract can be found in Figure \ref{fig:graphical-abstract}.}
\begin{figure}
    \centering
    \includegraphics[width=\textwidth]{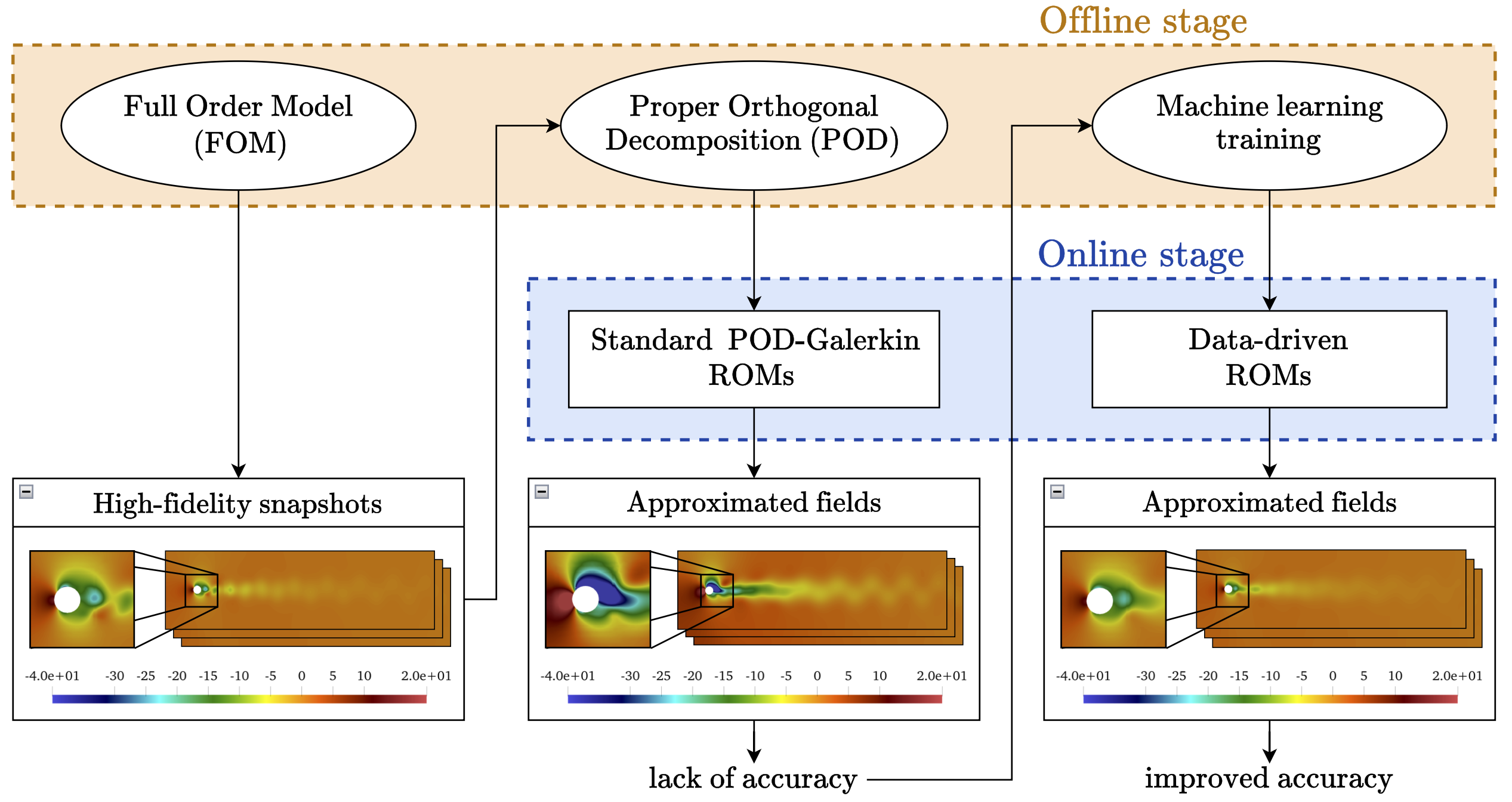}
    \caption{Visual abstract.}
    \label{fig:graphical-abstract}
\end{figure}
\end{abstract}

\listoffixmes



\section{Introduction}
\label{sec:intro}
Reduced Order Models (ROMs) \cite{degruyter1, degruyter2, degruyter3,rozza2008reduced, rozza2013reduced,aromabook} are a powerful tool used to reduce the computational effort of 
time-demanding simulations.
One of the fields where this class of techniques is widespread is Computational Fluid Dynamics,
where high-fidelity simulations may take days or weeks, even in the case of parallel computations
on many cores. For this reason, a simplified model is necessary to efficiently compute the solutions
for unseen configurations.

Most of the ROMs are built upon an \emph{offline-online} paradigm \cite{aromabook}.
The \emph{offline} stage consists of the computation of a large number of expensive high-fidelity simulations,
performed after setting the Full Order Model (FOM), which is usually a discretized version of complex PDEs, like the Navier--Stokes Equations (NSE).
The goal of this stage is the collection of the so-called \emph{snapshots}, namely the solutions of the simulations.
On the other hand, in the \emph{online} stage the full-order manifold is projected into a space with reduced dimensionality,
resulting in a reduced representation of the snapshots.
This reduction step may be assessed both with linear or nonlinear approaches. In particular, we employ the Proper Orthogonal Decomposition (POD) \cite{kosambi2016statistics, chatterjee2000introduction}, a linear technique which may be considered equivalent to the Principal Component Analysis (PCA) \cite{pearson1901liii, wold1987principal} and Singular Value Decomposition (SVD) \cite{lawson1995solving, golub1971singular, golub1965calculating}.

While the offline stage usually employs the numerical resolution of complex PDEs, in \emph{intrusive} ROMs the online stage allows for the resolution of ODEs, which consists of a reduced and simplified version of the FOM.

In particular, in this contribution, we focus on POD-Galerkin ROMs, that are based on a Galerkin projection \cite{noack1994low, bergmann2009enablers, kunisch2002galerkin, Azaiez2017}.
The main assumption of these models is that the solution may be approximated as a convex combination of a reduced number of global basis (the \emph{modes}),
whose coefficients are the reduced representations, namely the unknown variables in the ROM.
The above-mentioned linearity assumption may make the model inaccurate in the case of advection-dominated flows and where transport phenomena are dominant.

In the POD framework, the choice of the latent dimension is influenced by the singular values' decay and it has to be done \emph{a priori}. 

The goal of the ROM community is usually to retain a small number of modes, and to have efficient online simulations, namely a \emph{marginally-resolved} regime. In this regime, the number of modes is enough to capture the dynamics of the system, but the POD-Galerkin ROM may lead to inaccurate results \cite{rowley2004model, iollo2000stability}.

Indeed, in advection-dominated cases, where the decay is slow, POD-Galerkin ROMs may require hundreds of modes to provide accurate results.

The inaccuracy may be caused by the formation of spurious oscillations in the approximated solution or by the ill-conditioning of the reduced system.
This motivates the use of data-driven techniques to stabilize and/or enhance the system solution. Some examples of works integrating ROMs with machine learning strategies are \cite{zancanaro2021hybrid, papapicco2022neural, pichi2023artificial, ivagnes2024enhancing}.

More in detail, a widespread class of stabilization approaches is the regularized ROM, which aims to smooth out the noisy and inaccurate solution acting as a \emph{filter}. Some examples are the Evolve--Filter--Relax (EFR) ROMs \cite{gunzburger2019evolve, strazzullo2022consistency, strazzullo2023new}, Leray-ROMs, or more in general ROMs using filtering approaches inspired by Large Eddy Simulation (LES) \cite{xie2018data, mohebujjaman2018physically, girfoglio2023hybrid}.

Moreover, recently several researchers in the ROM community are putting their efforts into performing a system \emph{closure}, by adding to the system data-driven terms aimed to close the gap between the system solution and a reference optimal solution \cite{xie2020closure, ahmed2021closures, akhtar2012new}.

The work presented in this manuscript has the goal of improving the performance of the standard POD-Galerkin ROMs, briefly recalled in \ref{sec:pod-g-roms}, by the addition of some extra-terms into the reduced system. 

The data-driven terms are of two fundamental types:
\begin{itemize}
    \item[(i)] \emph{physics-based} data-driven terms, that re-introduce the turbulence modeling inside the ROM, inspired on \cite{hijazi2020data}, presented in \ref{subsec:met-eddy-viscosity};
    \item[(ii)] \emph{purely} data-driven terms, reintegrating the contributions of the neglected modes into the system, as done in \cite{xie2018data, mohebujjaman2018physically, mohebujjaman2019physically}, presented in \ref{subsec:met-corrections}.
\end{itemize}
The combination of the two strategies was already successfully applied to the test case of the periodic flow past a cylinder in \cite{ivagnes2023hybrid}.

The novelty of the present contribution consists of a machine-learning-based extension of the previous work to a parametrized setting, and on more challenging test cases.
In particular, we show the numerical results for the periodic turbulent flow past a circular cylinder in Subsection \ref{subsec:test-case-a}, and for the unsteady channel-driven cavity flow in Subsection \ref{subsec:test-case-b}.
In both test cases, different neural networks have been investigated to find the extra-terms (i) and (ii).

Finally, the discussion section \ref{subsec:discussion} and the conclusive part \ref{sec:conclusions} summarize the key results obtained in the contribution.

\section{Numerical methods}
\label{sec:methods}

This Section is dedicated to the presentation of the numerical methods used in
this project. In particular, it is divided in the following parts:
\begin{itemize}
    \item the presentation of the Full Order Model (FOM) in Section
    \ref{sec:fom}, which is used to compute the \emph{offline} solutions, named
    \emph{snapshots};
    \item a brief overview of the POD-Galerkin ROM approach in Section
    \ref{sec:pod-g-roms}, which is the starting model of our analysis, and the
    Pressure Poisson Equation (PPE) stabilization technique;
    \item the machine learning techniques used to enhance the standard
    POD-Galerkin approach, in Section \ref{sec:ml-roms}. We will briefly
    describe the \emph{physics-based data-driven} approach in Section
    \ref{subsec:met-eddy-viscosity} and the \emph{purely data-driven} approach
    in Section \ref{subsec:met-corrections}.
\end{itemize}

\subsection{Full Order Model (FOM)}
\label{sec:fom}

We call the fluid domain $\Omega \in \mathbb{R}^d$ with $d=2$ (since we are
considering only 2-dimensional domains), $\Gamma$ its boundary. We indicate by
$\bm{\mu}$ the flow parameters. In our specific test cases, $\bm{\mu} \in \mathbb{R}^2$, and it includes time $t \in [0, T]$ and the
kinematic viscosity $\nu$. $\bm{u}=\bm{u}(\bm{x}, \bm{\mu})$ is the flow velocity
vector field, $p=p(\bm{x}, \bm{\mu})$ is the pressure scalar field normalized by
the fluid density. In both the test cases we consider in the numerical results,
the FOM is based on the Unsteady Reynolds--Averaged Navier--Stokes (URANS)
formulation for incompressible flows, which consists of a time-averaged version
of the Navier--Stokes equations.

The main hypothesis that characterizes the RANS approach is the \emph{Reynolds
decomposition}~\cite{reynolds1895iv}. This theory is based on the assumption
that each flow field $\bm{s}$ can be expressed as the sum of its time
averaged value, indicated as $\overline{\bm{s}}$, and fluctuating parts,
indicated with $\bm{s}^{\prime}$. Such mean has different definitions depending
on the case of application. Classical choices are, for example, time averaging,
averaging along homogeneous directions or ensemble averaging.

We briefly recall here the standard URANS formulation for the incompressible
Navier-Stokes equations:
\begin{equation}
\begin{cases}
\dfrac{\partial \overline{u}_{i}}{\partial x_i} = 0 ,\smallskip\\
    \overline{u}_{j} \dfrac{\partial \overline{u}_{i} }{\partial x_j}=-\dfrac{\partial \overline{p}}{\partial x_i}+ \dfrac{\partial( 2 \nu\overline{\bm{E}}_{ij}- \mathcal{R}_{ij})}{\partial x_j},
    \label{RANS}
\end{cases}
\end{equation}

where the Einstein notation has been adopted, $\mathcal{R}_{ij}=\overline{u'_i
u'_j}$ is the Reynolds stress tensor, and
$\overline{\bm{E}}_{ij}=\dfrac{1}{2}\left(\dfrac{\partial
\overline{u}_{i}}{\partial x_j} + \dfrac{\partial \overline{u}_{j}}{\partial
x_i}\right)$ is the averaged strain rate tensor.

The URANS formulation in \eqref{RANS} needs to be coupled with a turbulence
model to close system \eqref{RANS}. In particular, we adopt the $\kappa-\omega$
\emph{Shear Stress Transport} (SST) model~\cite{menter1994two}.

This model belongs to the class of \emph{eddy viscosity models}, which are based
on the Boussinesq hypothesis, i.e. the turbulent stresses are related to the
mean velocity gradients as follows:
\[
 -\mathcal{R}_{ij}=2 \nu_t \overline{\bm{E}}_{ij} - \dfrac{2}{3} \kappa \delta_{ij},
    \label{bouss}
\]
where $\kappa=\frac{1}{2}\overline{u'_i u'_i}$ is the turbulent kinetic energy
and $\nu_t$ is the eddy viscosity. In general, in this case the URANS model is enriched with two additional transport equations, for $\kappa$ and $\omega$, respectively. For the complete model we refer the reader to
the original paper \cite{menter1994two}, and the extended RANS model including
the SST $\kappa-\omega$ equations can be found in \cite{hijazi2020data}.

The full order solutions of \ref{RANS} are computed by means of the open-source
software \emph{OpenFOAM}, which employs a finite volume discretization of the
URANS equations \cite{moukalled2016finite, jasak1996error}.

Following the finite volume method, developed and implemented in
\cite{moukalled2016finite, jasak1996error}, the computational domain is
discretized in polygonal control volumes and the partial differential equations
\ref{RANS} are integrated in each control volume and converted to algebraic
equations. In particular, the volume integrals are converted into surface
integrals by the divergence theorem, and discretized as sums of the fluxes at
the boundary faces of the control volumes~\cite{moukalled2016finite}. 

\subsection{POD-Galerkin Reduced Order Models}
\label{sec:pod-g-roms}

In this section, we provide a brief overview of the standard POD-Galerkin ROM
approach.
\medskip

Once all the high-fidelity simulations are run, considering completed the
offline stage, all the FOM snapshots, i.e., the FOM solutions for different parameter values  $\{\bm{\mu}_j\}_{j=1}^{N_{\mu}}$, in our case for different time
instances $\{t_i\}_{i=1}^{N_T}$ and for different viscosity values
$\{\nu_k\}_{k=1}^{N_M}$, where $N_{\mu}=N_T\times N_M$. In particular, each parameter is $[t_i, \nu_k], \, i=1, \dots, N_T, \, k=1, \dots, N_M$.
The POD is then applied on the full order snapshots' matrices:
\begin{equation*}
\mathcal{S}_u=\{\bm{u}(\bm{x},\bm{\mu}_1),...,\bm{u}(\bm{x},\bm{\mu}_{N_{\mu}})\}  \in \mathbb{R}^{N_u^h \times N_{\mu}}, \quad \mathcal{S}_p=\{p(\bm{x},\bm{\mu}_1),...,p(\bm{x},\bm{\mu}_{N_{\mu}})\}  \in \mathbb{R}^{N_p^h \times N_{\mu}}.
\end{equation*}
If we call $N_c$ the number of cells of the mesh considered, $N_u^h=d\, N_c$ and $N_p^h=N_c$ are the numbers of spatial degrees of freedom for the
velocity and pressure fields, respectively.

After the POD is applied to the snapshots' matrices, the following reduced POD
spaces are found:
\begin{equation}
\mathbb{V}^u_{\text{POD}}=\mbox{span}\{[\boldsymbol{\phi}_i]_{i=1}^{N_u}
\},\quad 
\mathbb{V}^p_{\text{POD}}=\mbox{span}\{[\chi_i]_{i=1}^{N_p}\},
\end{equation}
where $N_u \ll N_u^h$ and $N_p \ll N_p^h$, and
$[\boldsymbol{\phi}_i]_{i=1}^{N_u}$ and $[\chi_i]_{i=1}^{N_p}$ indicate the
velocity and pressure POD modes, respectively.

The main hypothesis of the POD is that each field can be then approximated as a
convex combination of its modes, namely:
\begin{equation}
\bm{u}(\bm{x},\bm{\mu}) \approx \bm{u}_r(\bm{x}, \bm{\mu})=\sum_{i=1}^{N_u} a_i(\bm{\mu})\boldsymbol{\phi}_i(\bm{x}), \quad
p(\bm{x},t)\approx p_r(\bm{x},\bm{\mu})=\sum_{i=1}^{N_p} b_i(\bm{\mu})\chi_i(\bm{x}).
\label{eq:appfields}
\end{equation}

Supposing that the reduced fields provide an accurate approximation of the
original snapshots, the reduced order system can be written in a compact form
as:
\begin{equation}
    \begin{cases}
        \bm{f}(\bm{a}, \bm{b}; \bm{\mu})=\bm{0},\\
        \bm{c}(\bm{a}; \bm{\mu}) = \bm{0}.
    \end{cases}
    \label{eq:compact-stand-rom}
\end{equation}
In \eqref{eq:compact-stand-rom} and in the following parts of the manuscript we
will omit the dependencies $\bm{a}(\bm{\mu})$ and $\bm{b}(\bm{\mu})$ for the
sake of brevity.

The above-mentioned model represents a system of DAEs whose unknowns are the vectors of coefficients
$\bm{a}=(\bm{a})_{i=1}^{N_u}$ (for the velocity), and
$\bm{b}=(\bm{b})_{i=1}^{N_p}$ (for the pressure).

However, it may happen that the model \eqref{eq:compact-stand-rom} is not stable
in terms of velocity-pressure coupling and this may lead to spurious
oscillations, causing the residuals to increase as time evolves. For this
reason, there exist in literature different stabilization approaches for the
velocity-pressure coupling. Two of those approaches are the supremizer method
(SUP-ROM) and the Pressure Poisson Equation approach (PPE-ROM). The first
approach consists in enriching the velocity reduced space with additional modes
named \emph{supremizer modes}, which are directly computed either from the
pressure modes (\emph{exact} supremizer approach) or from the pressure snapshots
(\emph{approximated} supremizer approach) \cite{rozza2007stability,
ballarin2015supremizer, stabile2018finite, ali2020stabilized}. The second method
simply replaces the reduced continuity equation with a reduced version of the
Pressure Poisson Equation, which is computed by taking the divergence of the
conservation of the momentum equation
\cite{akhtar2009stability,Stabile2017CAIM,stabile2018finite, noack2005need}.

In this project, we choose to adopt the PPE-ROM stabilization for two specific
reasons:
\begin{itemize}
    \item[(i)] it relies on a dedicated pressure equation and it allows for the
    introduction of specific \emph{pressure correction} terms only acting on the
    pressure accuracy, as will be specified in Section \ref{sec:ml-roms};
    \item[(ii)] it is consistent with velocity-pressure coupling at the FOM
    level.
\end{itemize}

In the case of PPE-ROMs we report here the standard POD-Galerkin ROM
formulations, which we will refer to as \emph{standard ROM}.

\begin{equation}
\begin{cases}
    \bm{M} \dot{\bm{a}}=\nu(\bm{B}+\bm{B_T})\bm{a}-\bm{a}^T \bm{C} \bm{a}-\bm{H}\bm{b}+\tau \left( \sum_{k=1}^{N_{\text{BC}}}(U_{BC,k}\bm{D}^k-\bm{E}^k \bm{a})\right)\, ,\\
    \bm{D}\bm{b}+ \bm{a}^T \bm{G} \bm{a} - \nu \bm{N} \bm{a}- \bm{L}=\bm{0}\, .
    \end{cases}
    \label{eq:ppe-rom}
\end{equation}

In the above formulation, the POD operators read as follows:
\begin{equation}
\begin{split}
&(\bm{M})_{ij}=(\boldsymbol{\phi}_i,\boldsymbol{\phi}_j)_{L^2(\Omega)}, \quad (\bm{P})_{ij}=(\chi_i,\nabla \cdot \boldsymbol{\phi}_j)_{L^2(\Omega)}\, ,\quad (\bm{B})_{ij}=(\boldsymbol{\phi}_i,\nabla \cdot \nabla \boldsymbol{\phi}_j)_{L^2(\Omega)}, \\
&(\bm{B_T})_{ij}=(\boldsymbol{\phi}_i,\nabla \cdot (\nabla \boldsymbol{\phi}_j)^T)_{L^2(\Omega)},\quad (\bm{C})_{ijk}=(\boldsymbol{\phi}_i,\nabla \cdot (\boldsymbol{\phi}_j \otimes \boldsymbol{\phi}_k))_{L^2(\Omega)}, \quad (\bm{H})_{ij}=(\boldsymbol{\phi}_i,\nabla \chi_j)_{L^2(\Omega)}\, \\
&(\bm{D})_{ij}=(\nabla \chi_i,\nabla \chi_j)_{L^2(\Omega)}, \quad 
(\bm{G})_{ijk}=(\nabla \chi_i,\nabla \cdot (\boldsymbol{\phi}_j \otimes \boldsymbol{\phi}_k))_{L^2(\Omega)}, \\ &(\bm{N})_{ij}=(\bm{n} \times \nabla \chi_i,\nabla \boldsymbol{\phi}_j)_\blue{{L^2(\Gamma)}}, \quad (\bm{L})_{ij}=(\chi_i,\bm{n} \cdot \boldsymbol{R}_t)_\blue{{L^2(\Gamma)}}\,,
\end{split}
\label{eq:operators-rom}
\end{equation}
vector $\bm{n}$ is the normal unitary vector to the domain boundary.

Moreover, the additional term in the momentum equation is a \emph{weak}
enforcement of the non-homogeneous Dirichlet boundary conditions, according to
the \emph{penalty method}. The penalization factor, namely $\tau$  tuned through
a sensitivity analysis on the specific problem considered \cite{hijazi2020data,
star2019extension}. In general, bigger values of the penalization factor lead to
a stronger enforcement of the boundary conditions. At the same time, large values of $\tau$ also lead to larger condition number for system \eqref{eq:ppe-rom} and, hence, to less stability.

The matrices $\bm{E}^k$ and
vectors $\bm{D}^k$ are defined as follows:
\[(\bm{E}^k)_{ij}=(\boldsymbol{\phi}_i,
\boldsymbol{\phi}_j)_{L^2(\Gamma_{D_k})}, \quad
(\bm{D}^k)_{i}=(\boldsymbol{\phi}_i)_{L^2(\Gamma_{D_k})}, \text{ for all
}k=1,...,N_{\text{BC}},\]

where $\Gamma_{D_i}$ is the i-th boundary with non-homogeneous Dirichlet conditions, $N_{\text{BC}}$ is the total number of boundaries where non-homogeneous Dirichlet boundary conditions are imposed.

We can finally rewrite the formulation in System \eqref{eq:ppe-rom} in the
following compact form, which is inherited from the form in
\eqref{eq:compact-stand-rom}:
\begin{equation}
\begin{cases}
    \bm{f}(\bm{a}, \bm{b}; \bm{\mu})=\bm{0},\\
    \bm{h}(\bm{a}, \bm{b}; \bm{\mu})=\bm{0},
\end{cases}
    \label{eq:compact-ppe-rom}
\end{equation}


where the expression $\bm{h}(\bm{a}, \bm{b}; \bm{\mu})$ includes all the already
mentioned contributions of the PPE and replaces the continuity equation $\bm{c}(\bm{a})=\bm{0}$ in \eqref{eq:compact-stand-rom}.

Although the PPE-ROM approach is well-known in the POD-Galerkin framework to
provide an efficient stabilization, it may fail in capturing the flow dynamics
especially in unsteady and turbulent test cases, like the cases we will consider
in Section \ref{sec:results}. This motivates the following Section, describing
the data-driven techniques used to improve the accuracy.

\subsection{Machine-learning-enhanced Reduced Order Models}
\label{sec:ml-roms}
This part of the manuscript is dedicated to the description of the machine
learning techniques adopted to enhance the ROM results. In particular, this
Section is divided into two parts:
\begin{itemize}
    \item the description of the models used to model the reduced approximation
    of the \emph{eddy viscosity} field, which is here introduced to take into
    account the turbulent behaviour of the flow (section
    \ref{subsec:met-eddy-viscosity});
    \item the description of the models used to approximate the
    \emph{correction/closure} terms added into the System \eqref{eq:ppe-rom}
    (section \ref{subsec:met-corrections}).
\end{itemize}

\subsubsection{Physics-based data-driven methods}
\label{subsec:met-eddy-viscosity}
The main goal of this Section is to add the turbulence modeling inside the
reduced formulation, building the so-called \emph{physics-based data-driven
model}. As specified in Section \ref{sec:fom}, the URANS approach used to model
the FOM usually includes turbulence by using an \emph{eddy viscosity} model. For
instance, in our case we considered the SST $\kappa-\omega$ model, which adds to
the URANS equations the transport equations for $\kappa$ and $\omega$.

It is important to highlight that only the projections of the momentum and the PPE equations are considered at the reduced order level. Since the projections of the $\kappa$ and $\omega$ transport equations are not taken into account, it is important to include the turbulent viscosity modeling in the momentum and PPE equations.

Indeed, an approximation version of the eddy viscosity field
can be modeled and included in the reduced system \eqref{eq:ppe-rom}, as done in
\cite{hijazi2020data}, as follows:
\begin{equation*}
\end{equation*}
where $\eta_i(\bm{x})$ is the $i$-th eddy viscosity mode evaluated through a POD
procedure, just as done for the velocity and pressure fields, $g_i(\bm{\mu})$ is
the corresponding coefficient, $N_{\nu_t}$ is the \emph{a-priori} selected
number of modes retained for the eddy-viscosity field.

When using a PPE approach, the URANS system \eqref{RANS} at the full order level
can be written as follows:
\begin{equation}
    \begin{cases}
    \dfrac{\partial \overline{\bm{u}}}{\partial t}+ \nabla \cdot (\overline{\bm{u}} \otimes \overline{\bm{u}})=\nabla \cdot \left[-\overline{p} \bm{I} +(\nu+\nu_t) \left(\nabla \overline{\bm{u}} + (\nabla \overline{\bm{u}})^T \right) \right] & \text{ in } \Omega \times [0,T]\, ,\\
    \Delta \overline{p}=-\nabla \cdot (\nabla \cdot (\overline{\bm{u}} \otimes \overline{\bm{u}})) +\nabla \cdot \left[ \nabla \cdot \left( \nu_t \left(\nabla \overline{\bm{u}} +(\nabla \overline{\bm{u}})^T \right) \right) \right] & \text{ in }\Omega \, ,\\
    + \text{ Boundary Conditions }& \text{ on } \Gamma \times [0,T]\, ,\\
    + \text{ Initial Conditions } & \text{ in } (\Omega,0)\, .
    \end{cases}
    \label{RANS-PPE}
\end{equation}
Consequently, the dynamical system \eqref{eq:ppe-rom} takes the following form:
\begin{equation}
    \begin{cases}
    \bm{M} \dot{\bm{a}}=\nu(\bm{B}+\bm{B_T})\bm{a}-\bm{a}^T \bm{C} \bm{a}+ \bm{g}^T (\bm{C}_{\text{T1}} +\bm{C}_{\text{T2}}) \bm{a}-\bm{H}\bm{b}+\tau \left( \sum_{k=1}^{N_{\text{BC}}}(U_{\text{BC},k}\bm{D}^k-\bm{E}^k \bm{a})\right) \, ,\\
    \bm{D}\bm{b}+ \bm{a}^T \bm{G} \bm{a} -\bm{g}^T(\bm{C}_{\text{T3
    }} +\bm{C}_{\text{T4
    }})\bm{a} - \nu \bm{N} \bm{a}- \bm{L}=\bm{0}\, ,
    \end{cases}
    \label{eq:ppe-rom-turb}
\end{equation}
where:
\[
\begin{split}
    &(\bm{C}_{\text{T1}})_{ijk}=(\boldsymbol{\phi}_i, \eta_j \nabla \cdot \nabla \boldsymbol{\phi}_k)_{L^2(\Omega)} \, ,\quad
    (\bm{C}_{\text{T2}})_{ijk}=(\boldsymbol{\phi}_i, \nabla \cdot \eta_j (\nabla \boldsymbol{\phi}_k)^T)_{L^2(\Omega)}\, ,\\
    &(\bm{C}_{\text{T3}})_{ijk}=(\nabla \chi_i, \eta_j \nabla \cdot \nabla \boldsymbol{\phi}_k)_{L^2(\Omega)}\, , \quad (\bm{C}_{\text{T4}})_{ijk}=(\nabla \chi_i, \nabla \cdot \eta_j(\nabla \boldsymbol{\phi}_k)^T)_{L^2(\Omega)}\,.
\end{split}
\]
In the dynamical system defined in \eqref{eq:ppe-rom-turb} the number of
unknowns is $N_u$ for the velocity, $N_p$ for the pressure, and $N_{\nu_t}$ for
the eddy viscosity. However, the number of equations is $N_u+N_p$. Thus, there
are more unknowns than equations and the system is not closed. In order to close
the systems, the eddy viscosity coefficients $[g_i(\bm{\mu})]_{i=1}^{N_{\nu_t}}$
can be computed considering the mapping $\bm{g}=\mathcal{G}(\bm{a}, \bm{\mu})$
through either \emph{interpolation}~\cite{lazzaro2002radial,
micchelli1986interpolation} or \emph{regression} techniques. An interpolation
technique was exploited in \cite{hijazi2020data}, following the POD-I approach
\cite{wang2012comparative,walton2013reduced,salmoiraghi2018free}. In this paper,
the mapping is a feed-forward fully-connected neural network whose loss function
is the mean squared error (MSE) between the output of the neural network
$\mathcal{G}(\bm{a}, \bm{\mu})$ and the known eddy viscosity coefficients
$\bm{g}_r^{proj}(\bm{\mu})$ found from the POD procedure. The compact form of
system \eqref{eq:ppe-rom-turb} is:
\begin{equation}
\begin{cases}
    \bm{f}(\bm{a}, \bm{b}, \bm{g};\bm{\mu}) = \bm{0}\, ,\\
    \bm{h}(\bm{a}, \bm{b}, \bm{g};\bm{\mu}) = \bm{0}\, .
\end{cases}
    \label{eq:compact-ppe-rom-turb}
\end{equation}

\subsubsection{Purely data-driven methods}
\label{subsec:met-corrections}
The main goal here is to identify a strategy to model the
\emph{correction/closure} terms already exploited in previous works like
\cite{mohebujjaman2019physically, mou2021data, san2013proper, ahmed2021closures,
ivagnes2023pressure, dar2023artificial, ivagnes2023hybrid}, but in a parametric
setting.

The procedure used is as follows.
\begin{enumerate}
    
    \item[\textbf{1}.] Select a reduced dimension for the velocity, $r \equiv
    N_u$, and for the pressure, $q \equiv N_p$. The sum $r+q$ is the dimension
    of the reduced system \eqref{eq:ppe-rom}.
    \item[\textbf{2}.] Select two bigger dimensions $d>r$ and $h>q$, for
    instance in literature we have $d=\mathrm{k}r$, $h=\mathrm{k}q$, with
    $\mathrm{k}$ an integer.
    \item[\textbf{3}.] Select one or more operators, that we name $\mathcal{C}$.
        For instance, if we only consider the non-linear operators in System
        \eqref{eq:ppe-rom} we obtain: \begin{equation} \mathcal{C}(\bm{a})=
        \begin{bmatrix}
            \bm{a}^T \bm{C} \bm{a}\\
            \bm{a}^T \bm{G} \bm{a} \end{bmatrix}.
    \end{equation}
    \item[\textbf{4}.] Compute the term $\mathcal{C}(\bm{a}_d^{proj}) \in \mathbb{R}^{d+h}$:
    \begin{equation}
    \mathcal{C}(\bm{a}_d^{proj})=
    \begin{bmatrix}
            (\bm{a}_d^{proj})^T \bm{C}_d \bm{a}_d^{proj}\\
            (\bm{a}_d^{proj})^T \bm{G}_h \bm{a}_d^{proj}
        \end{bmatrix},
        \label{eq:operator-big}
    \end{equation}
    where $\bm{a}_d^{proj}$ is the velocity coefficients' vector found by
    directly projecting the field on the POD subspace of dimension $d$. The
    operators $\bm{C}_d$ and $\bm{G}_h$ are expressed as:
    \[
    \begin{split}
    &(\bm{C}_d)_{ijk}=(\boldsymbol{\phi}_i,\nabla \cdot (\boldsymbol{\phi}_j \otimes \boldsymbol{\phi}_k))_{L^2(\Omega)}, \quad i, j, k=1, \dots, d, \\
    &(\bm{G}_h)_{ijk}=(\nabla \chi_i,\nabla \cdot (\boldsymbol{\phi}_j \otimes \boldsymbol{\phi}_k))_{L^2(\Omega)}, \quad i=1, \dots, h; j, k=1, \dots, d.
    \end{split}
    \]
    \item[\textbf{5}.] Compute the term $\mathcal{C}(\bm{a}_r^{proj}) \in \mathbb{R}^{r+q}$:
    \begin{equation}
    \mathcal{C}(\bm{a}_r^{proj})=
    \begin{bmatrix}
            (\bm{a}_r^{proj})^T \bm{C} \bm{a}_r^{proj}\\
            (\bm{a}_r^{proj})^T \bm{G} \bm{a}_r^{proj}
        \end{bmatrix},
        \label{eq:operator-reduced}
    \end{equation}
    where $\bm{a}_r^{proj}$ is the velocity coefficients' vector found by
    directly projecting the field on the POD subspace of dimension $r$. The
    operators $\bm{C}$ and $\bm{G}$ are the reduced operators appearing in
    \eqref{eq:ppe-rom} and already defined in \eqref{eq:operators-rom}.
    \item[\textbf{6}.] Compute the \emph{exact} correction term as follows:
    \begin{equation}
        \boldsymbol{\tau}^{exact} = \overline{\mathcal{C}(\bm{a}_d^{proj})}^r - \mathcal{C}(\bm{a}_r^{proj}),
        \label{eq:exact-correction}
    \end{equation}
    where $\overline{(\cdot)}^r$ acts like a filter and indicates that only the
    first $r$ components should be retained. In our example, we will have $\boldsymbol{\tau}^{exact} \in \mathbb{R}^{r+q}$:
     \begin{equation}
    \boldsymbol{\tau}^{exact}=
    \begin{bmatrix}
            \overline{(\bm{a}_d^{proj})^T \bm{C}_d \bm{a}_d^{proj}}^r - (\bm{a}_r^{proj})^T \bm{C} \bm{a}_r^{proj}\\
            \overline{(\bm{a}_d^{proj})^T \bm{G}_p \bm{a}_d^{proj}}^q - (\bm{a}_r^{proj})^T \bm{G} \bm{a}_r^{proj}
        \end{bmatrix},
        \label{eq:exact-nonlinear-correction}
    \end{equation}
    where we retain the first $r$ and $q$ components for the first and second
    term, respectively.
    \item[\textbf{7}.] Finally, we have to choose a \emph{model} $\mathcal{M}$
    able to approximate with good accuracy the exact correction term, also for
    unseen configurations, namely time and/or parameter extrapolation. In
    particular: $\boldsymbol{\tau}^{approx}=\boldsymbol{\tau}^{approx}(\bm{a},
    \bm{b}, \bm{\mu})=\mathcal{M}(\bm{a}, \bm{b}, \bm{\mu})$, where $\bm{\mu}$
    includes the parameters of the test case taken into account. Some of the
    possible choices to model the mapping $\mathcal{M}$ are the following.
    \begin{itemize}
        \item Choose an \emph{ansatz}, for example a quadratic ansatz with
        respect to the coefficients' vectors:
        \begin{equation}
            \mathcal{M}(\bm{a}, \bm{b})= 
                \tilde{\bm{A}} \begin{bmatrix}\bm{a}\\ \bm{b} \end{bmatrix}+ \begin{bmatrix}\bm{a}^T& \bm{b}^T \end{bmatrix} \tilde{\bm{B}}  \begin{bmatrix}\bm{a}\\ \bm{b} \end{bmatrix}.
                \label{eq:quadratic-ansatz}
            \end{equation}
            where the matrix $\tilde{\bm{A}} \in \mathbb{R}^{(r+q)\times(r+q)}$
            and the operator $\tilde{\bm{B}} \in \mathbb{R}^{(r+q)\times(r+q)\times(r+q)}$. This approach is easily addressed in
            time-dependent problems by solving a minimization problem:
            \begin{equation}
                \tilde{\bm{A}}, \tilde{\bm{B}}=\text{arg} \min_{\bm{A}, \bm{B}}\left(\sum_{j=1}^{M} \left\lVert {\bm{A}} \begin{bmatrix}(\bm{a}_r^{proj}(t_j))^T\\ (\bm{b}_q^{proj}(t_j))^T\end{bmatrix}+ \begin{bmatrix}\bm{a}_r^{proj}(t_j)& \bm{b}_q^{proj}(t_j) \end{bmatrix} {\bm{B}}  \begin{bmatrix}\bm{a}_r^{proj}(t_j)\\ \bm{b}_q^{proj}(t_j) \end{bmatrix} - \bm{\tau}^{exact}(t_j) \right\rVert^2 \right).
                \label{eq:min-least-squares}
            \end{equation}
            In previous works, like \cite{mou2021data,
            mohebujjaman2019physically, ivagnes2023hybrid, ivagnes2023pressure},
            the minimization problem \eqref{eq:min-least-squares} is typically
            re-written as a least squares problem. However, this can be done if
            time is the only parameter of the problem. Indeed, matrices
            $\tilde{\bm{A}}$ and $\tilde{\bm{B}}$ does not depend on time.
            Parametric test cases require more advanced mappings depending also
            on the remaining parameter(s).
            
        \item Train a \emph{neural network} (NN), which takes as input the
        coefficients and the parameters of the problem $(\bm{a}, \bm{b},
        \bm{\mu})$ and gives as output the approximated correction coefficients
        $\bm{\tau}^{approx}$. In this case, the mapping $\mathcal{M}$ can be
        modeled considering different architectures, like a feed-forward
        fully-connected NN or a recurrent neural network, for instance a
        Long-Short Term Memory (LSTM). Moreover, also other \emph{ad-hoc}
        architectures can be created to model the initial mapping, depending on
        the test case we are considering, as we will see in the results' section
        \ref{sec:results}. The comparisons among the results obtained with the
        different neural networks are reported in the numerical results' Section
        \ref{sec:results}, while the details on the architectures are reported
        in the supplementary parts \ref{appendix-correction-case-a} and
        \ref{appendix-correction-case-b}, for test case \textbf{a} and
        \textbf{b}, respectively. This second approach is suitable for test case
        that are only time-dependent, but also for problems depending on more
        parameters, whose dependency is easily addressed by adding the
        parameters as input to the neural network.
    \end{itemize}
\end{enumerate}

When we employ a reduced eddy-viscosity model, as in
\eqref{eq:compact-ppe-rom-turb}, the operator $\mathcal{C}$ can also include the
turbulence contributions, namely terms $\bm{g}^T (\bm{C}_{\text{T1}} +
\bm{C}_{\text{T2}})\bm{a}$ and $\bm{g}^T (\bm{C}_{\text{T3}} +
\bm{C}_{\text{T4}})\bm{a}$.

Moreover, as in the previous works \cite{ivagnes2023pressure,
ivagnes2023hybrid}, we identify the first $r$ elements of the correction term
$\bm{\tau}=\mathcal{M}(\bm{a}, \bm{b}, \bm{\mu})$ as the \emph{velocity}
correction, whereas the last $q$ components are identifies as \emph{pressure}
correction. Hence, the $\bm{\tau}$ vector can be decomposed as $(\bm{\tau}_u,
\bm{\tau}_p)$, where $\bm{\tau}_u \in \mathbb{R}^r$ and $\bm{\tau}_p \in \mathbb{R}^q$.

We are finally able to write the so-called \emph{hybrid data-driven ROM}, which
combines the physics-based and the purely data-driven based approaches. Its
final compact formulation reads as follows:
\begin{equation}
    \begin{cases}
    \bm{f}(\bm{a}, \bm{b}, \bm{g};\bm{\mu})+\bm{\tau}_u(\bm{a}, \bm{b},\bm{\mu}) = \bm{0}\, ,\\
    \bm{h}(\bm{a}, \bm{b}, \bm{g};\bm{\mu})+\bm{\tau}_p (\bm{a}, \bm{b}, \bm{\mu})= \bm{0}\, .
\end{cases}
\label{eq:hyb-dd-rom}
\end{equation}

In particular, when we consider a \emph{purely data-driven} approach, the exact
correction $\bm{\tau}^{exact}$ considered in Step \textbf{6}. is computed
considered all the operators in the reduced system \eqref{eq:ppe-rom}. When
considering a \emph{hybrid data-driven} approach, we include in the exact
correction also the contribution depending on the reduced eddy viscosity
$\bm{g}$.

It is important to remark that in this case the reduced eddy viscosity is not an
input to the mapping $\mathcal{M}$ but it is only embedded in the exact
correction.


\section{Numerical results}
\label{sec:results}
\begin{figure}[htpb!]
    \centering
    \subfloat[Parameters for test case $(\bm{a})$]{\includegraphics[width=0.53\textwidth]{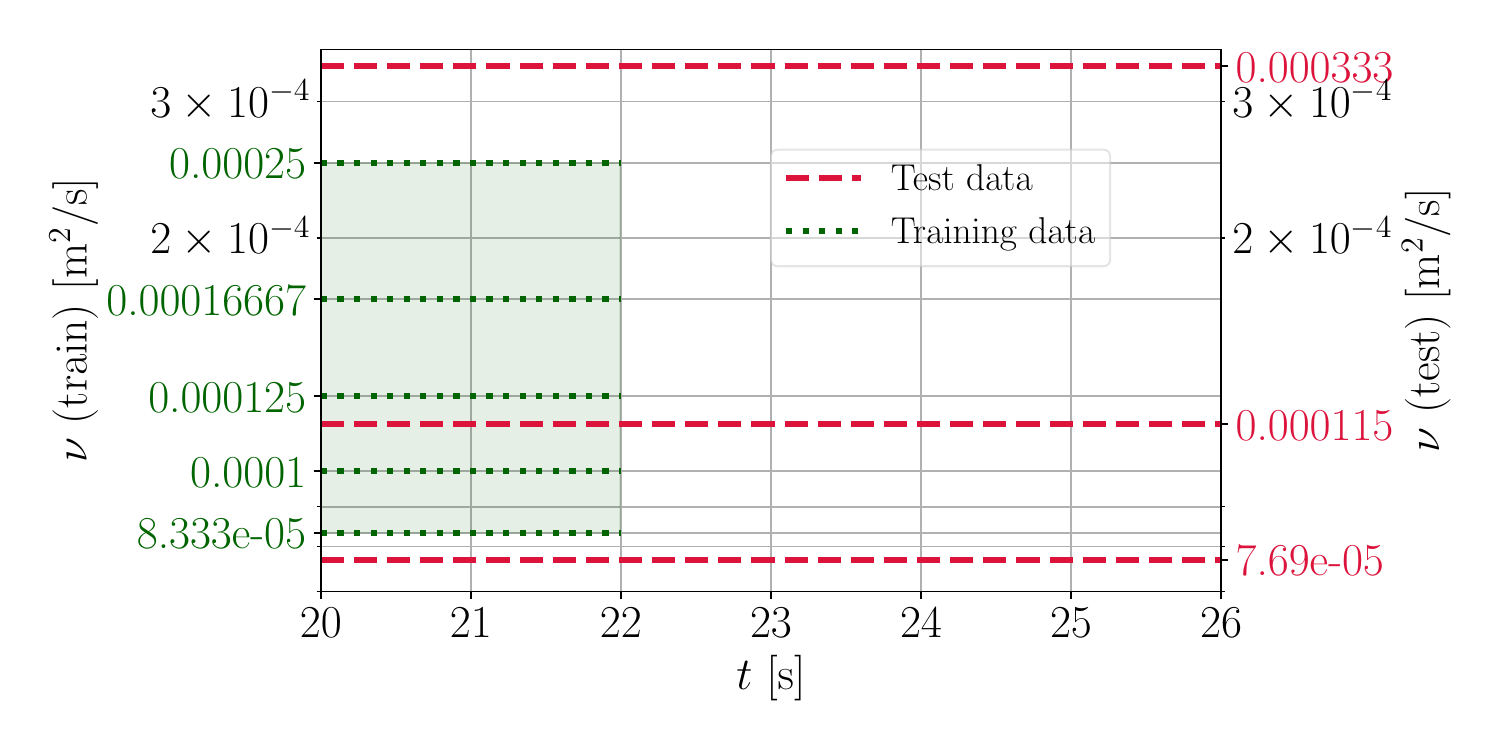}}
    \subfloat[Parameters for test case $(\bm{b})$]{\includegraphics[width=0.53\textwidth]{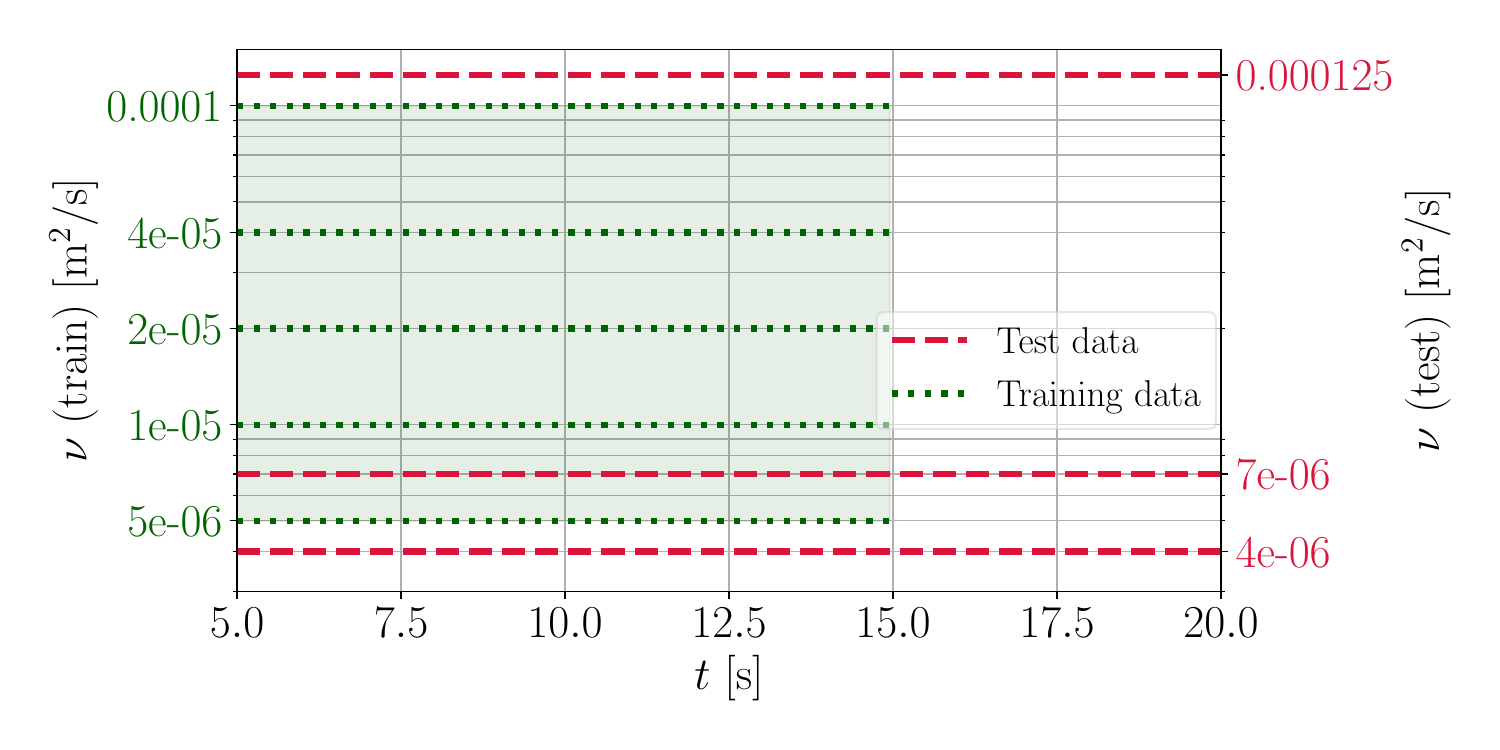}}
    \caption{Graphical representations of the sets of parameters used in the offline stage ($\Toffline$ and $\nutrain$), in green, and in the online stage ($\Tonline$ and $\nutest$), in red. The green shadow describes the accessible offline/training area. The parameters are also specified in Table \ref{tab:table-params}.}
    \label{fig:params-plot}
\end{figure}

In this paper, we test all the numerical methods presented in \ref{sec:methods} on two different test cases:
\begin{enumerate}
    \item[$(\bm{a})$] the periodic flow past a circular cylinder;
    \item[$(\bm{b})$] the channel-driven cavity flow.
\end{enumerate}
In both test cases, we consider as parameters both time and the viscosity value, i.e. the Reynolds number.
In particular, we consider three different sets of snapshots for training and testing the reduced order models, that are represented in Figure \ref{fig:params-plot} and specified in Table \ref{tab:table-params}, for both the test cases. In particular, we use the following notation: 
\begin{itemize}
    \item $\Toffline$ and $\nutrain$: the time window and set of viscosities, respectively, used to collect the snapshots for the parametric POD. In the test cases we consider, the collected data are also used to train the neural networks exploited to compute the data-driven terms;
    \item $\Tonline$ and $\nutest$: the time window and set of viscosities, respectively, used to run the online simulations. These sets are used to test the parametric data-driven approaches in different extrapolation settings for both parameters.
\end{itemize}

\begin{table}[htpb!]
    \centering
    \begin{tabular}{ccc}
       \toprule
        \multirow{4}{*}{Test case ($\bm{a}$)} &$\Toffline [\unit{\second}]$  & [20, 22]\\
        \cline{2-3}
        &$\Tonline [\unit{\second}]$ & [20, 28] \\
        \cline{2-3}
        &$\nutrain [\unit{\metre^2 \per \second}]$ & $\{\num{8.33e-5}; \num{1e-4}; \num{1.25e-4}; \num{1.67e-4}; \num{2.5e-4}\} $ \\
        \cline{2-3}
        &$\nutest [\unit{\metre^2 \per \second}]$ & $\{\num{7.69e-5}; \num{1.15e-4}; \num{3.33e-4}\}$\\
        \hline
        \multirow{4}{*}{Test case ($\bm{b}$)}&$\Toffline [\unit{\second}]$  & [5, 15]\\
        \cline{2-3}
        &$\Tonline [\unit{\second}]$ & [5, 20] \\
        \cline{2-3}
        &$\nutrain [\unit{\metre^2 \per \second}]$ &  $\{\num{5e-6}; \num{1e-5}; \num{2e-5}; \num{4e-5}; \num{1e-4}\} $\\
        \cline{2-3}
        &$\nutest [\unit{\metre^2 \per \second}]$ & $\{\num{4e-6};\num{7e-6};\num{1.25e-4}\}$\\
        \bottomrule
    \end{tabular}
    \caption{Different sets used to train and test the data-driven ROMs for both test cases.}
    \label{tab:table-params}
\end{table}

Moreover, we specify that we collect the snapshots every $\Delta t_{\text{offline}}=\Delta t_{\text{online}}=\num{4e-3} \unit{\second}$ for test case $(\bm{a})$, and every $\Delta t_{\text{offline}}=\Delta t_{\text{online}}=\num{0.05} \unit{\second}$ for test case $(\bm{b})$.

The present Section is divided into two parts, one for each test case, and in both Subsections we focus on the following aspects:
\begin{enumerate}
    \item[$(\bm{i})$] Description of the FOM, POD eigenvalues decay, and selection of a reduced number of modes for the online stage;
    \item[$(\bm{ii})$] Neural networks' predictions of the data-driven correction terms;
    \item[$(\bm{iii})$] Comparison among standard ROM and DD-ROMs in terms of relative $L^2$ errors of the velocity and pressure fields with respect to the high-fidelity solutions;
    \item[$(\bm{iv})$] Graphical representation of the results at the final online time instance.
\end{enumerate}

In the analysis of step $(\bm{ii})$, we consider multiple neural networks with the same architecture but trained with different (random) weights initialization. Hence, we will show the performance of the regression models by displaying the average prediction (on all the pre-trained networks), and the true coefficients, namely the coefficients obtained by projection. We will also include the confidence interval of our model.
We introduce the following notation: $\bm{s}^{\nu^{\star}} (t)=\bm{s}(t;\nu^{\star} )$ represents the full-order field of interest at the viscosity value $\nu^{\star}$ and evolving in time, $\tilde{\bm{s}}$ indicates the ROM approximation of the field. In particular, in our case the fields are velocity or pressure: $\bm{s}=\bm{u}$ or $\bm{s}=p$.
Following this notation, the $99.7 \%$ confidence interval $\mathcal{I}_{99.7 \%}$ is computed as:
\begin{equation}
    \mathcal{I}_{99.7 \%}(t)=[i_{-}(t), i_{+}(t)], \text{ where }i_{\pm}(t)=\dfrac{\mu(\tilde{\bm{s}}^{\nu^{\star}} (t)) \pm 3 \sigma(\tilde{\bm{s}}^{\nu^{\star}} (t))}{N_{\text{networks}}},
    \label{eq:confidence}
\end{equation}
where $N_{\text{networks}}$ is the number of networks considered for the statistics in our case $10$. The confidence interval reported in \eqref{eq:confidence} is well-defined under the assumption of a Gaussian distribution. In particular, the difference among the $N_{\text{networks}}$ models considered is the (random) initialization of the weights.


\medskip

For what concerns step $(\bm{iii})$, the expression of the time-dependent relative errors used to compare the ROMs performances is computed as follows.
\begin{equation}
    \mathcal{E}_s^{\nu^{\star}} (t) = \dfrac{\| \tilde{\bm{s}}^{\nu^{\star}} (t) - \bm{s}^{\nu^{\star}} (t) \|_{L^2(\Omega)}}{\| \bm{s}^{\nu^{\star}} (t) \|_{L^2(\Omega)}}\, .
    \label{eq:l2-errs-exp}
\end{equation}
Also in this case, when dealing with DD-ROMs we will represent the statistics, namely the accuracy obtained considering in System \eqref{eq:hyb-dd-rom} the averaged prediction of the neural networks, and the corresponding confidence interval.


\subsection{Test case \textbf{a}: periodic flow past a cylinder}
\label{subsec:test-case-a}

The periodic flow past a circular cylinder is a wide-known benchmark in fluid dynamics and in the field of reduced order modeling. In this Section we will focus on a parametric version of this test case, where the parameter considered is the viscosity that indirectly parametrizes the Reynolds number.


\subsubsection{Offline stage}
\label{subsubsec:fom-a}
\begin{figure}[htpb!]
    \centering
    \subfloat[Domain with notation]{\includegraphics[width=0.5\textwidth]{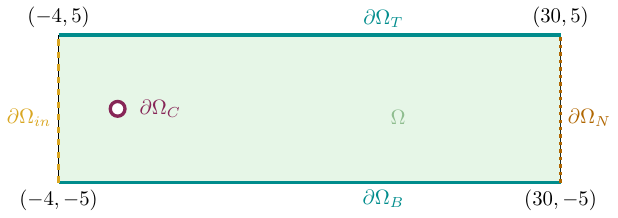}}
    \subfloat[Full order mesh]{\includegraphics[width=0.5\textwidth, trim={0cm 10cm 0cm 10cm}, clip]{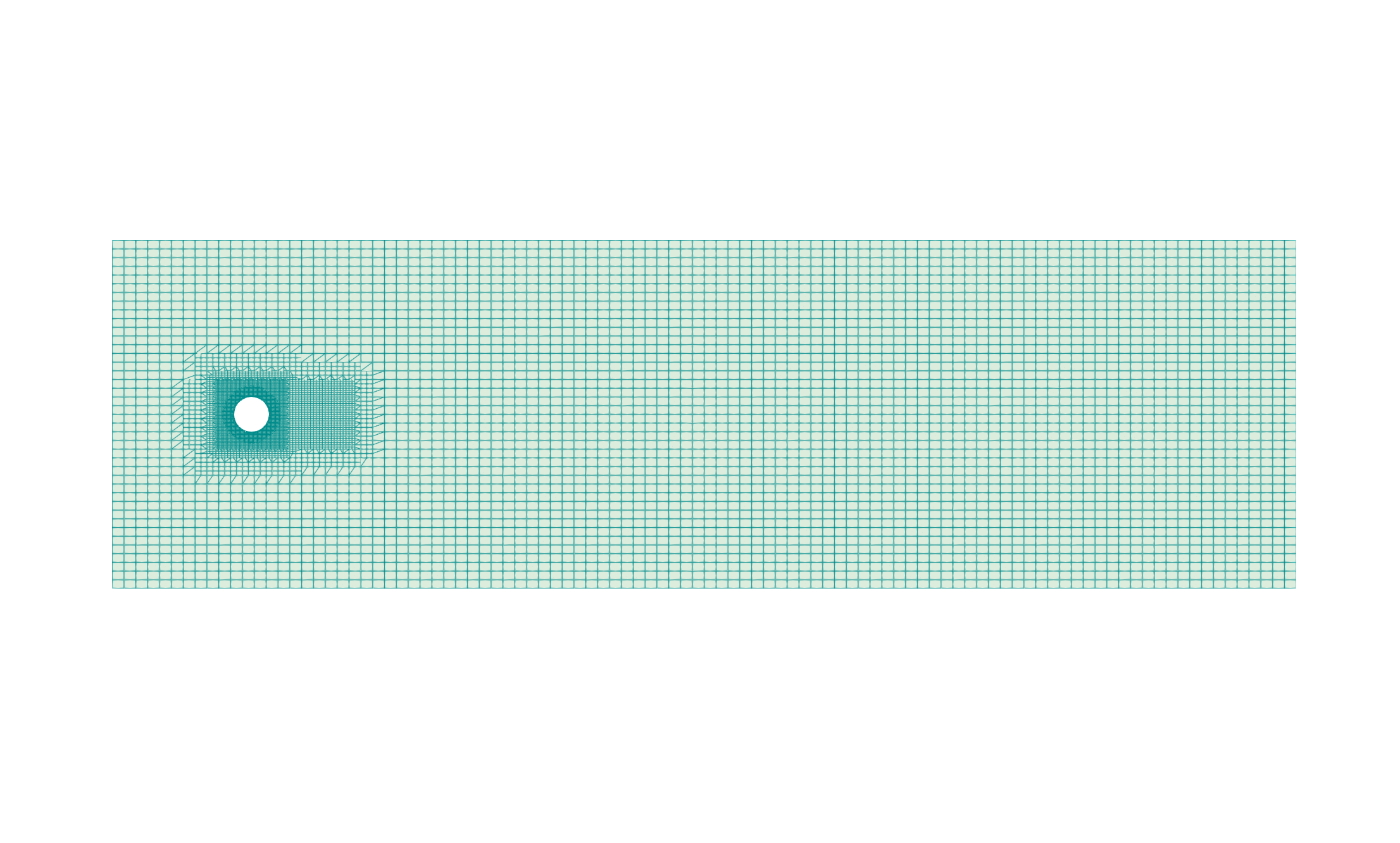}}
    \caption{The domain and full order mesh considered for the periodic flow around a circular cylinder.}
    \label{fig:cyl-domain}
\end{figure}

The domain and mesh for this test case are represented in Figure \ref{fig:cyl-domain}. Following the notation in the above-mentioned Figure, the boundary conditions read as follows.

\begin{equation*}
\text{On }\partial \Omega_{in}:
    \begin{cases}
        \bm{u} = (U_{in}, 0),\\
        \nabla p \cdot \bm{n} = 0;
    \end{cases}
    \quad
    \text{On }\partial \Omega_T \cup \partial \Omega_B:
    \begin{cases}
        \bm{u} \cdot \bm{n} = 0,\\
        \nabla p \cdot \bm{n} = 0;
        
    \end{cases}
\end{equation*}

\begin{equation*}
        \text{On }\partial \Omega_N:
    \begin{cases}
        \nabla \bm{u} \cdot \bm{n} = 0,\\
        p = 0;
        
    \end{cases}
    \quad
            \text{On }\partial \Omega_C:
    \begin{cases}
       \bm{u} = \bm{0},\\
        \nabla p \cdot \bm{n} = 0.
    \end{cases}
\end{equation*}

\begin{figure}[htpb!]
    \centering
    \subfloat[POD cumulative eigenvalues]{\includegraphics[height=4.7cm]{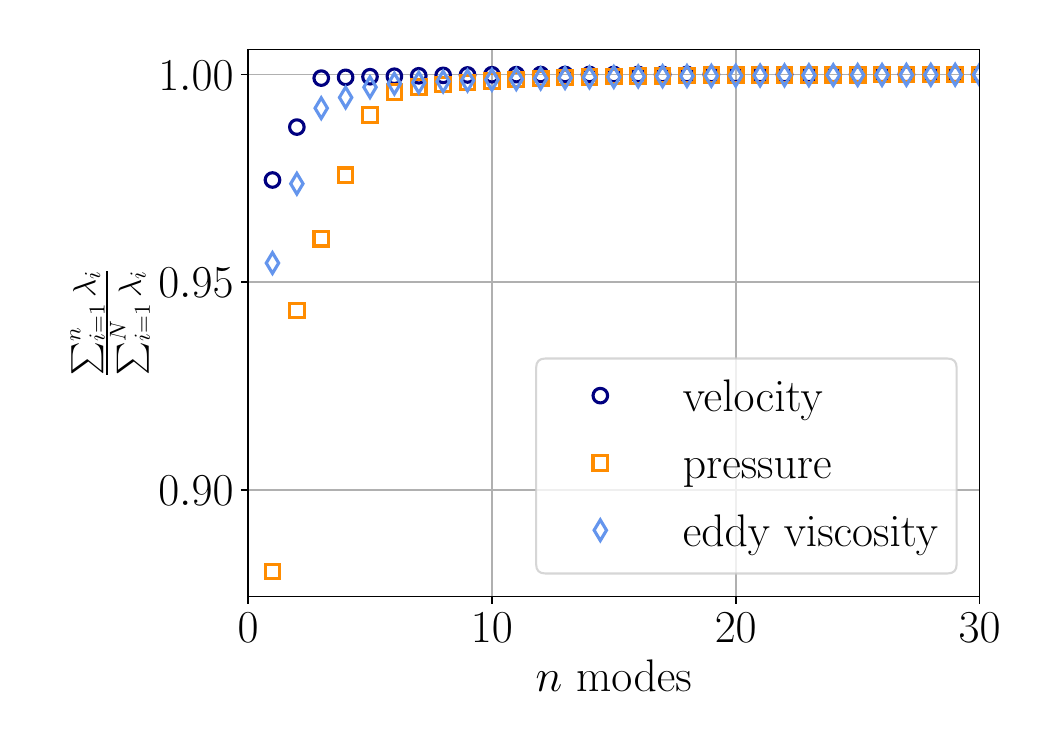}}
    \subfloat[POD eigenvalues decay]{\includegraphics[height=4.7cm]{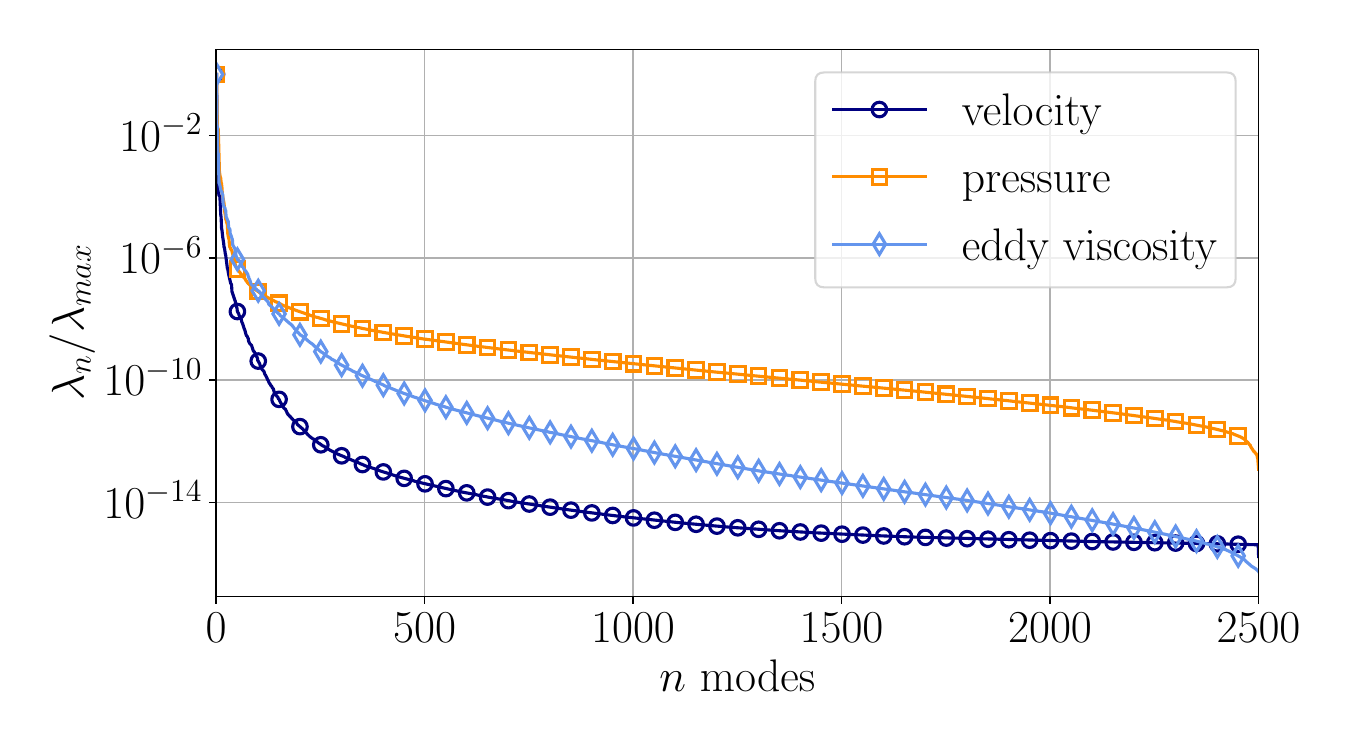}}
    \caption{Cumulative eigenvalues and eigenvalues decay for the test case \textbf{a}.}
    \label{fig:eig-a}
\end{figure}

The time and parameter ranges considered to run the FOM and collect the POD snapshots have been introduced in Figure \ref{fig:params-plot}. Moreover, we collect a total amount of $500$ snapshots for each viscosity value, for a total of $2500$ snapshots. The snapshots are taken every $\SI{0.004}{\second}$, the POD cumulative eigenvalues and the corresponding decay are represented in Figure \ref{fig:eig-a}. We can see a fast decay, especially for the velocity field. However, as we will see in the following part, the standard POD-Galerkin approach is characterized by lack of accuracy in the \emph{marginally-resolved} regime, even if the number of modes considered is enough to capture the dynamics of the problem.

We decide to focus our analysis on a number of modes $N_u=N_{\nu_t}=r=3$ and $N_p=q=3$, namely in the \emph{marginally-resolved} regime.

\subsubsection{Neural networks' performance}
\label{subsec:networks-a}

We focus here on the analysis of the accuracy of the neural networks used to compute both the \emph{eddy viscosity} and the \emph{correction} coefficients.

For what concerns the eddy viscosity coefficients, we consider as mapping $\mathcal{G}(\bm{a}, \nu, t)$ a multi-layer perceptron neural network, whose hyper-parameters' details can be found in the supplementary results in Section \ref{appendix-correction-case-a}.

\begin{figure}[htpb!]
    \centering
    \includegraphics[width=\textwidth]{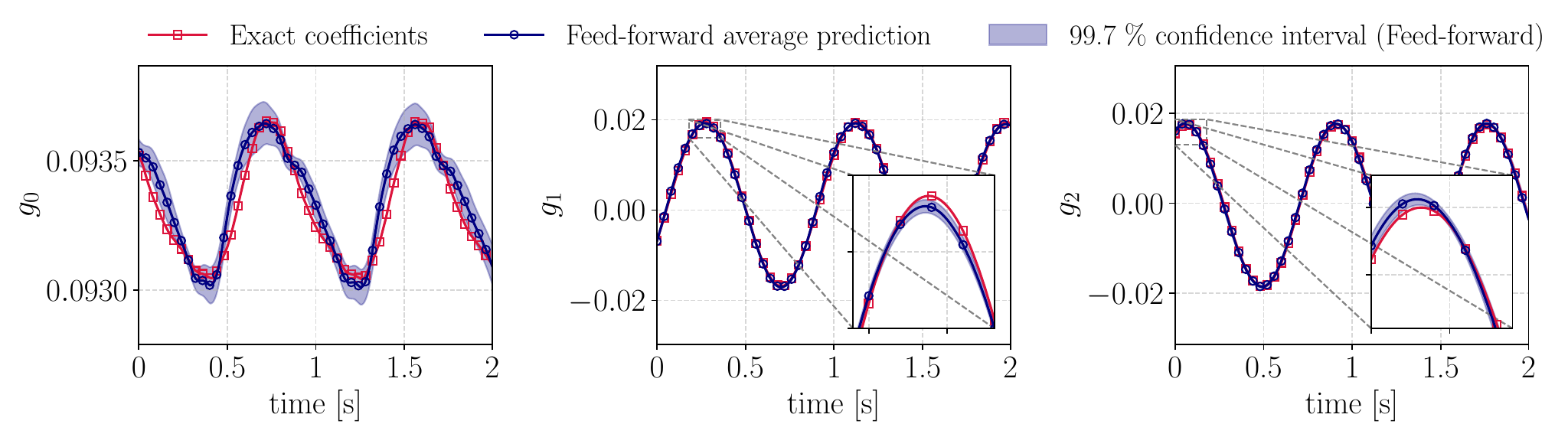}
    \caption{Eddy viscosity coefficients for $\nu=\SI{8.33e-5}{\metre \per \second}$ (in $\nutrain$): neural networks predictions' average and corresponding confidence interval, and POD projected exact coefficients.}
    \label{fig:cyl-coeffs-eddy-0}
\end{figure}

Figure \ref{fig:cyl-coeffs-eddy-0} shows the averaged prediction of the neural networks considered for the $\mathcal{G}$ mapping, the corresponding confidence interval, and the exact coefficients. The prediction is here computed for a viscosity value in the train set $\nutrain$ and for the same time window considered offline for collecting the POD snapshots. We can notice that the neural network has a good performance, since the prediction is close to the truth and the confidence interval is small. In particular, the performance is better for the eddy viscosity coefficients $g_1$ and $g_2$ than for the first component $g_0$. However, this would not affect the results in terms of accuracy, as we will show in the analysis in \ref{subsubsec:dd-roms-a}.

\begin{figure}[htpb!]
    \centering
    \includegraphics[width=\textwidth]{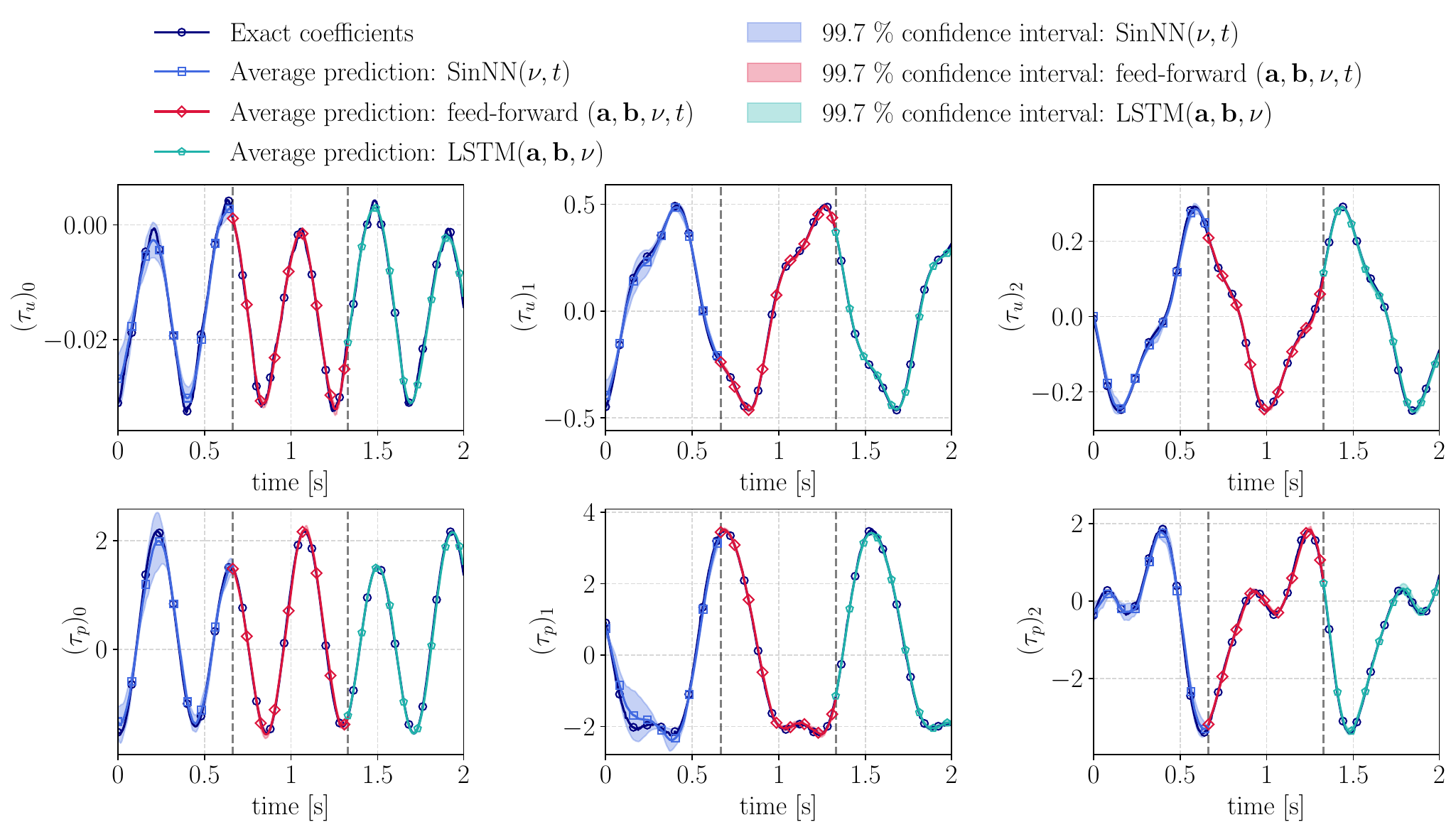}
    \caption{Correction terms' coefficients for $\nu=\SI{8.33e-5}{\metre \per \second}$ (in $\nutrain$): neural networks predictions' average and corresponding confidence interval, and POD projected exact coefficients.}
    \label{fig:cyl-coeffs-turb-0}
\end{figure}

\begin{figure}[htpb!]
    \centering
    \includegraphics[width=\textwidth]{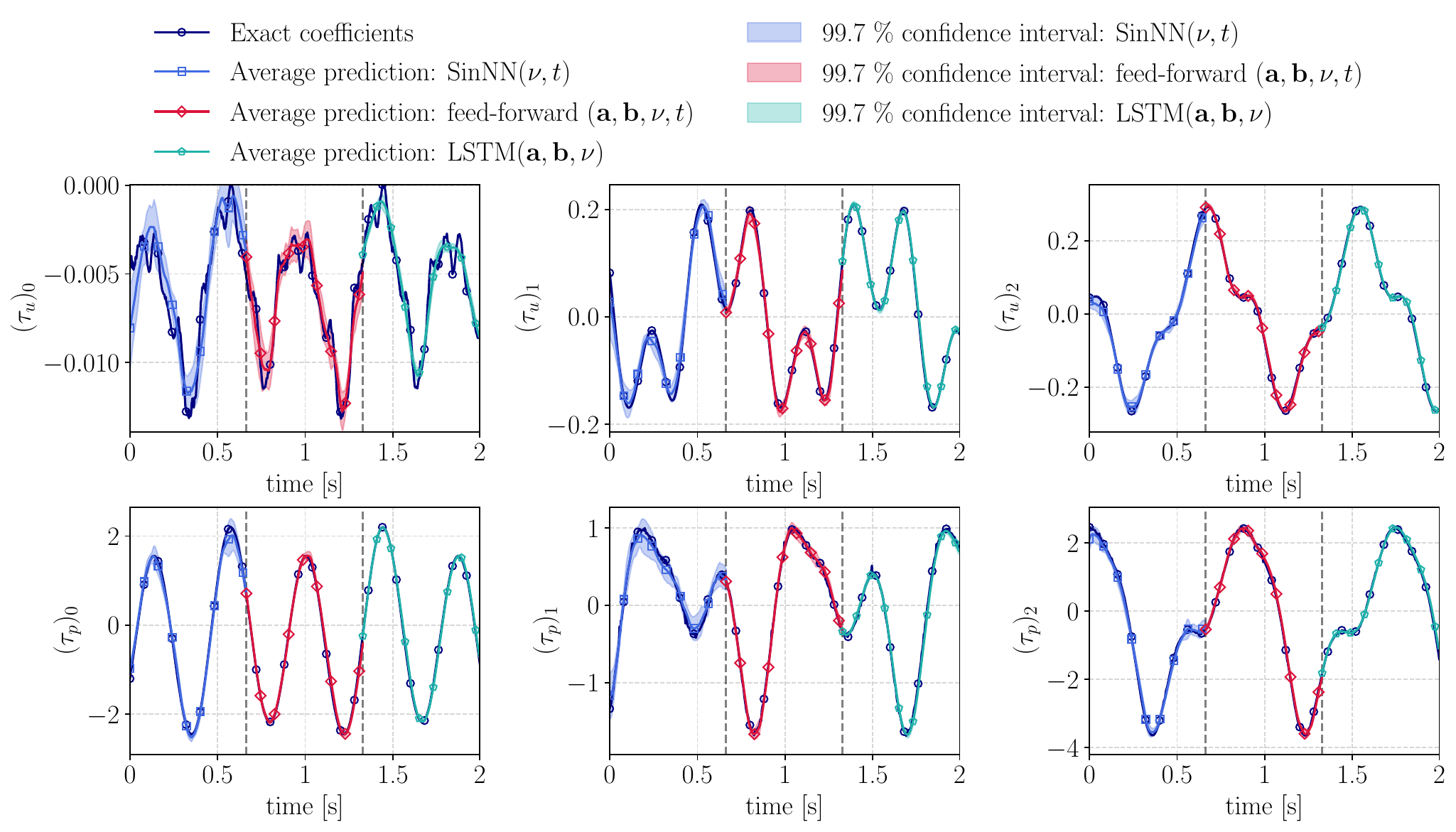}
    \caption{Correction terms' coefficients for $\nu=\SI{1.25e-4}{\metre \per \second}$ (in $\nutrain$): neural networks predictions' average and corresponding confidence interval, and POD projected exact coefficients.}
    \label{fig:cyl-coeffs-turb-2}
\end{figure}

Figures \ref{fig:cyl-coeffs-turb-0} and \ref{fig:cyl-coeffs-turb-2} represent the performance of three different types of neural networks considered for the mapping $\mathcal{M}$, on two different training viscosities. In particular, we named the networks as:
\begin{itemize}
    \item Feed-forward$(\bm{a}, \bm{b}, \nu, t)$, where $(\bm{a}, \bm{b}, \nu, t)$ are the inputs, with a multi-layer perceptron standard architecture;
    \item LSTM$(\bm{a}, \bm{b}, \nu)$, where the time parameter is not considered since it is already embedded in the Long--Short Term Memory framework;
    \item SinNN$(\nu, t)$, which is an \emph{ad-hoc} architecture, presented in the supplementary part, in Figure \ref{fig:sinNN}. It is inspired by the fact that the coefficients are in our case periodic, as can be seen from the eddy viscosity coefficients in Figure \ref{fig:cyl-coeffs-eddy-0}. Hence, the periodicity is integrated in the network by adding sinusoidal expressions in it.
\end{itemize}
The hyperparameters of the above-mentioned neural networks are reported in Section \ref{appendix-correction-case-a}.

From the comparison between Figures \ref{fig:cyl-coeffs-turb-0} and \ref{fig:cyl-coeffs-turb-2}, it can be seen that the shape of the first coefficient of $\bm{\tau}_u$ is more complex in \ref{fig:cyl-coeffs-turb-2}. Therefore, also the prediction of the neural network is less accurate for this specific parameter. Moreover, the correction coefficients represented in these Figures also include the eddy viscosity contributions, but similar considerations can be withdrawn in the case not including the turbulent terms.

\subsubsection{DD-ROMs performance}
\label{subsubsec:dd-roms-a}

\begin{figure}[htpb!]
    \centering
    \includegraphics[width=\textwidth]{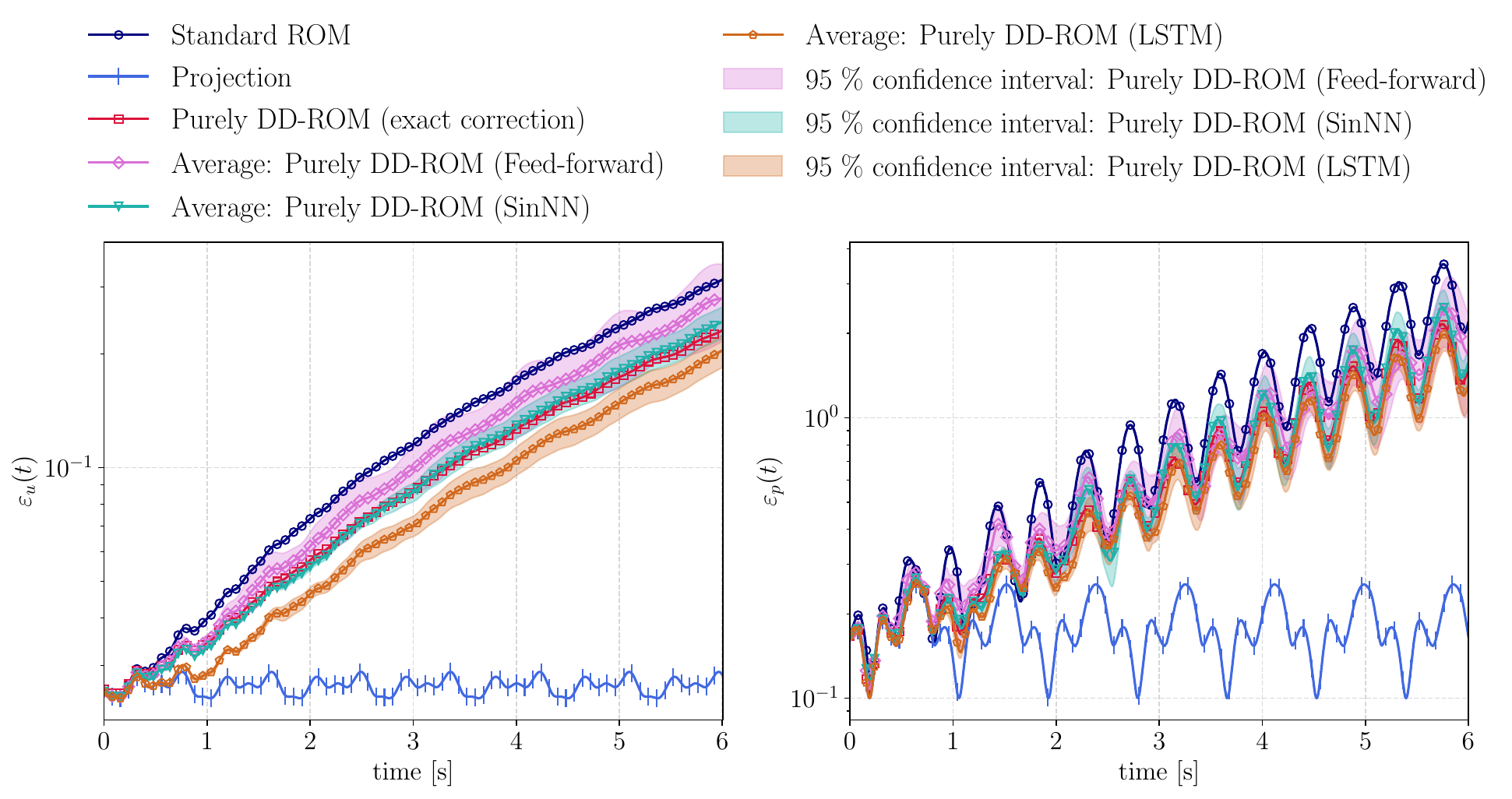}
    \caption{Relative errors for $\nu \in \nutest$ (inside the train dataset range) for the pressure and the velocity magnitude fields, for the purely DD-ROMs, standard POD-Galerkin ROMs. The reconstruction error is here also represented.}
    \label{fig:err-noturb-test-int}
\end{figure}

\begin{figure}[htpb!]
    \centering
    \includegraphics[width=\textwidth]{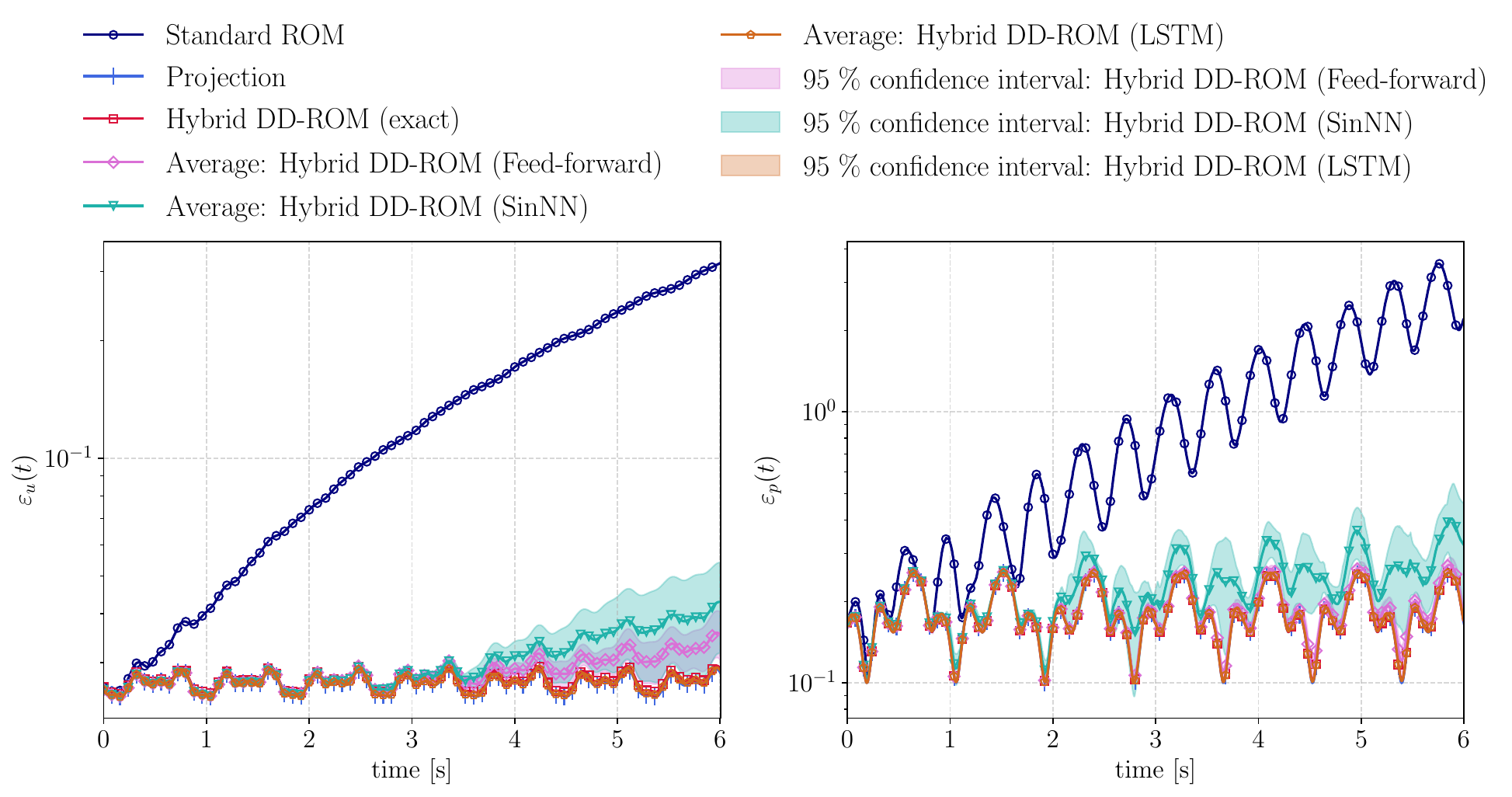}
    \caption{Relative errors for $\nu \in \nutest$, $\nu=\SI{1.15e-4}{\metre \per \second}$ for the pressure and the velocity magnitude fields, for the physics-based DD-ROMs, hybrid DD-ROMs, standard POD-Galerkin ROMs. The reconstruction error is here also represented.}
    \label{fig:err-turb-test-int}
\end{figure}

Figures \ref{fig:err-noturb-test-int} and \ref{fig:err-turb-test-int} show the relative errors obtained with all the DD-ROMs here considered, making a comparison with the standard ROM and with the reconstruction error. 

We report also the result of the DD-ROM obtained inserting in the reduced system \eqref{eq:hyb-dd-rom} the \emph{exact} correction term. From now on we will call this model the \emph{exact purely DD-ROM} or \emph{exact hybrid DD-ROM}, depending if the turbulence is included or not.

The viscosity parameter is a test parameter internal to the range considered for collecting the snapshots. Moreover, the online simulation is run until $6$ seconds, namely in a time extrapolation regime, since the snapshots are collected in the window $[0, 2]$ s.

If we consider the purely DD-ROM (Figure \ref{fig:err-noturb-test-int}), all the models slightly improve the results of the standard ROM. In particular, the LSTM outperforms both the exact purely DD-ROM and the standard ROM. 

Indeed, the mapping $\mathcal{M}$ is built such that its output is as close as possible to the exact correction, but the final DD-ROM accuracy may also be better than the exact model.

Figure \ref{fig:err-turb-test-int} shows that the hybrid DD-ROM significantly improves the standard ROM results in an unseen setting for both time and viscosity. The LSTM is the best-performing machine learning algorithm for this test parameter. Indeed, the average prediction of the LSTM leads to an accuracy close to the reconstruction error and to the exact hybrid DD-ROM. Moreover, the model provides a more confident prediction than the others, since the confidence interval is smaller than the others. 
In addition, the SinNN has the worst performance among the machine learning model considered. The reason for this fact may be that the SinNN does not take as inputs the velocity coefficients. Hence, it does not take into account the evolution of the dynamics leading to a worse performance in a test setting.

\begin{figure}[htpb!]
    \centering
    \includegraphics[width=0.6\textwidth]{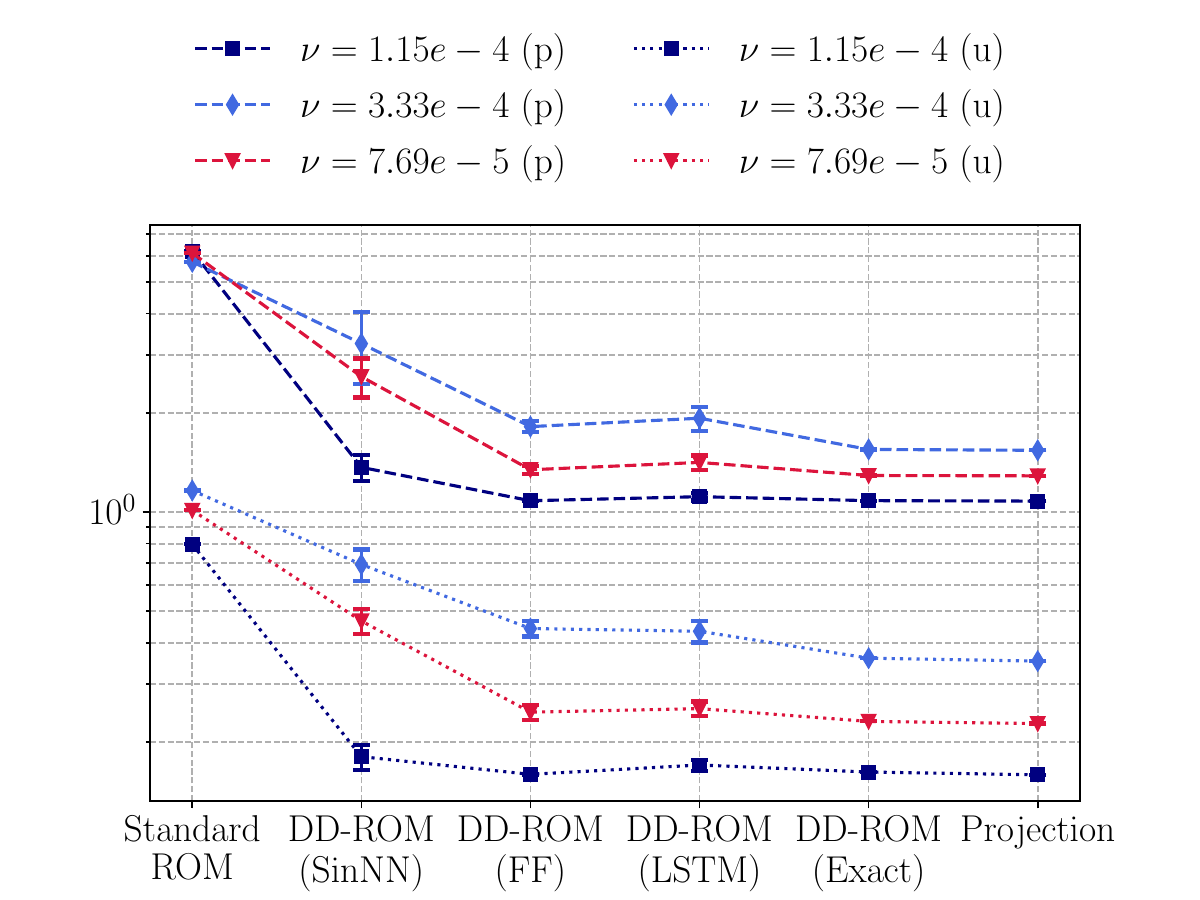}
    \caption{Integrals of relative errors over time for $\nu \in \nutest$ for the pressure and the velocity magnitude fields, in the standard ROM, hybrid DD-ROMs, and for the projection.}
    \label{fig:integrals-cylinder}
\end{figure}

For a global comparison of the different ROM methods, we refer to Figure \ref{fig:integrals-cylinder}, showing the integrals of the relative errors over the online time window $[20, 26]$ s, for test viscosity values.
The performance of all the methods is better for $\nu=\SI{1.15e-4}{\metre^2 \per \second}$ since it is within the parameters' range considered for the POD.
All the hybrid DD-ROM models outperform the standard approach, but the most accurate methods are the standard multi-layer perceptron and the LSTM, as previously pointed out. 

\subsubsection{Graphical results}
\label{subsubsec:graph-a}
\begin{figure}[htpb!]
    \centering
    \subfloat[]
    {\includegraphics[width=0.5\textwidth]{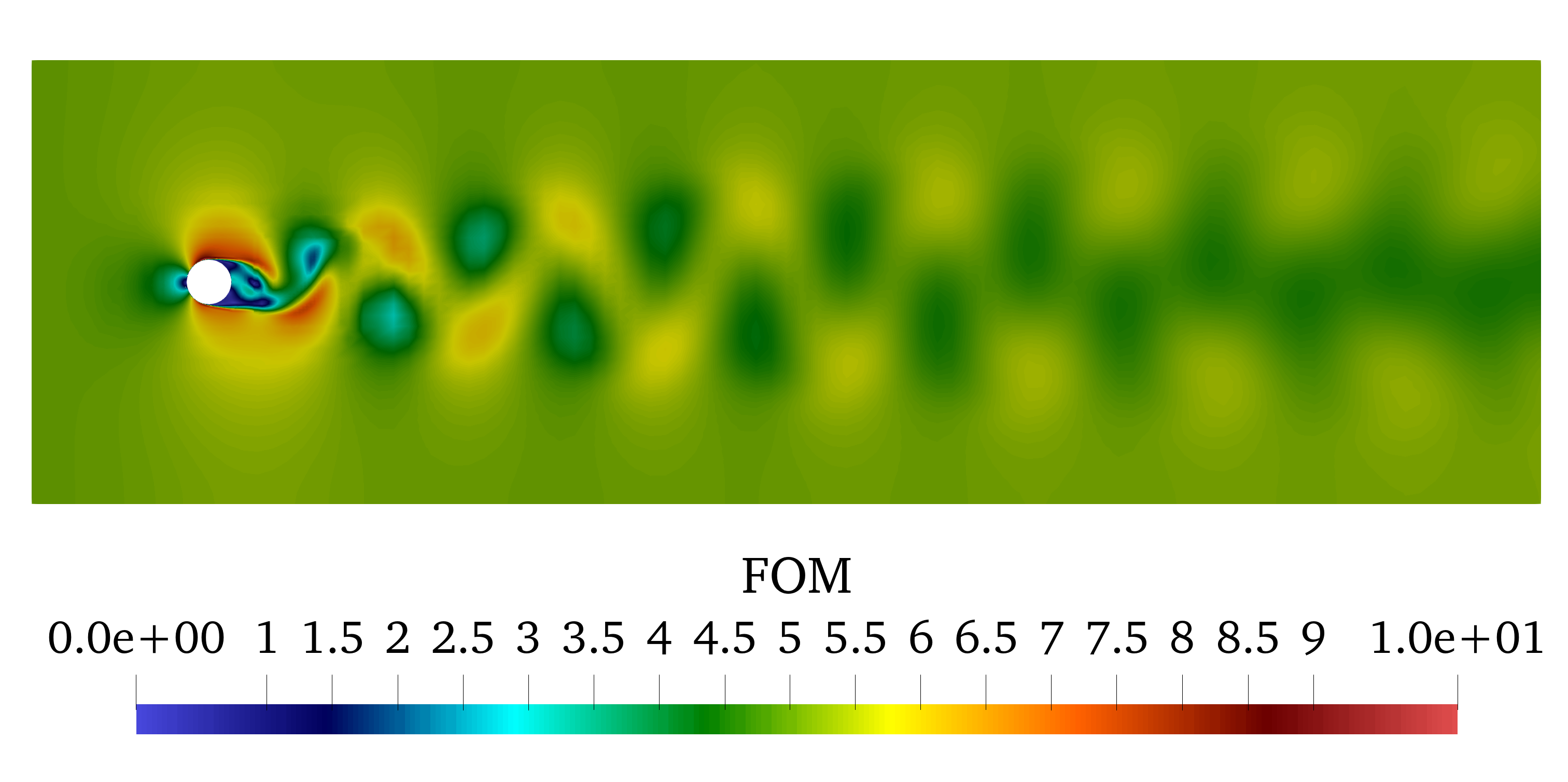}}\\
    \subfloat[]
    {\includegraphics[width=0.5\textwidth]{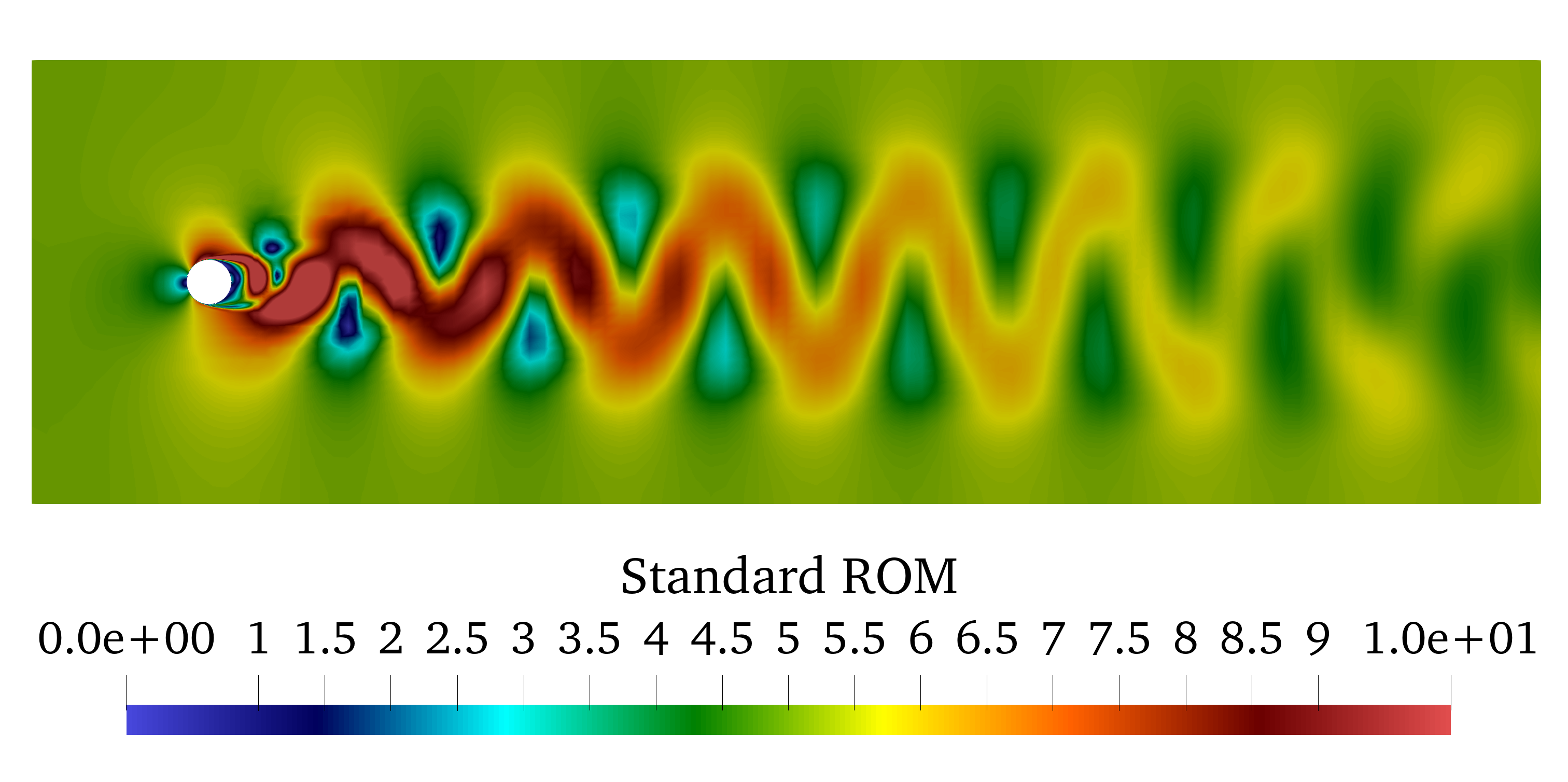}}
    \subfloat[]
    {\includegraphics[width=0.5\textwidth]{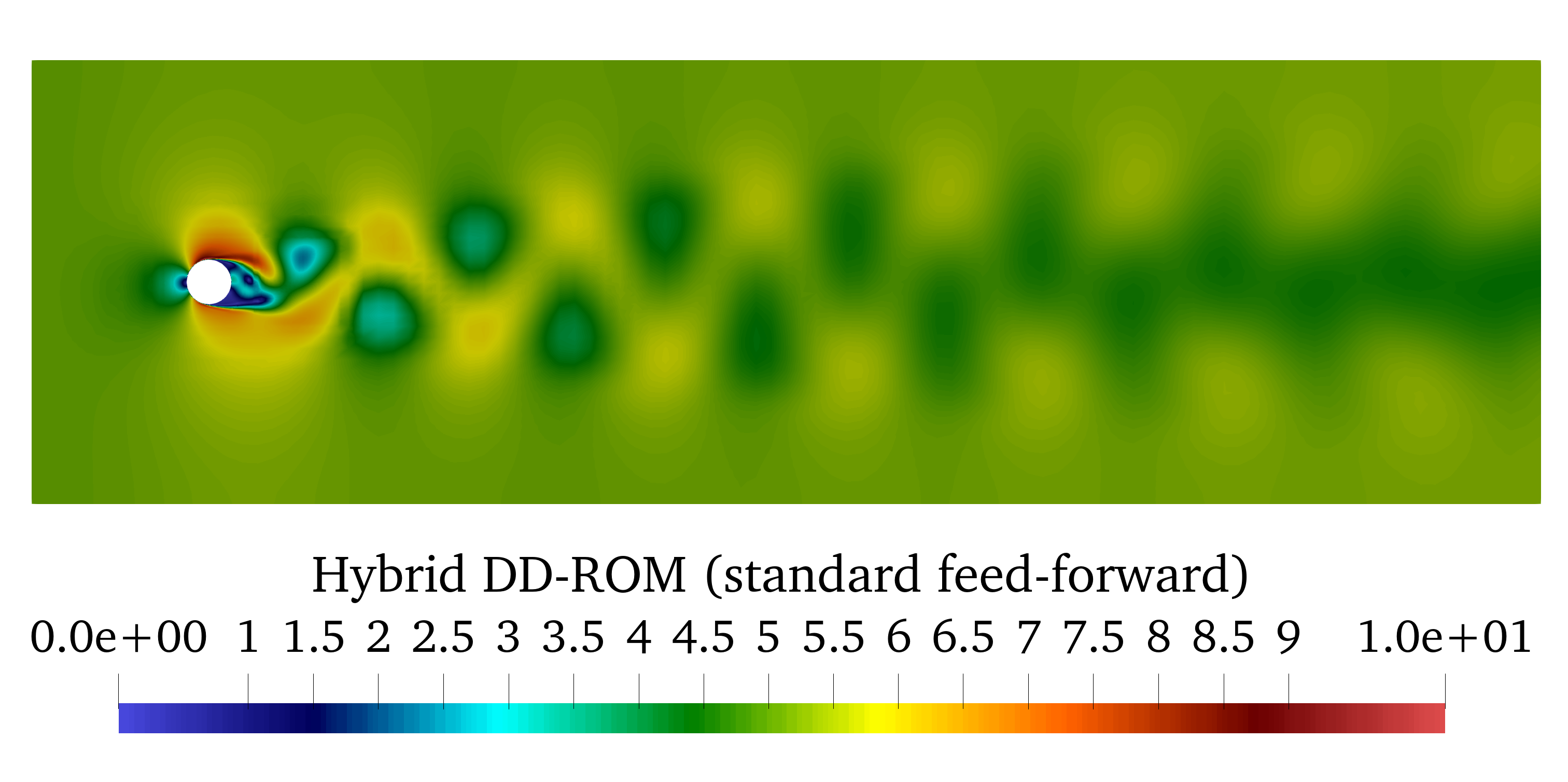}}\\
    \subfloat[]
    {\includegraphics[width=0.5\textwidth]{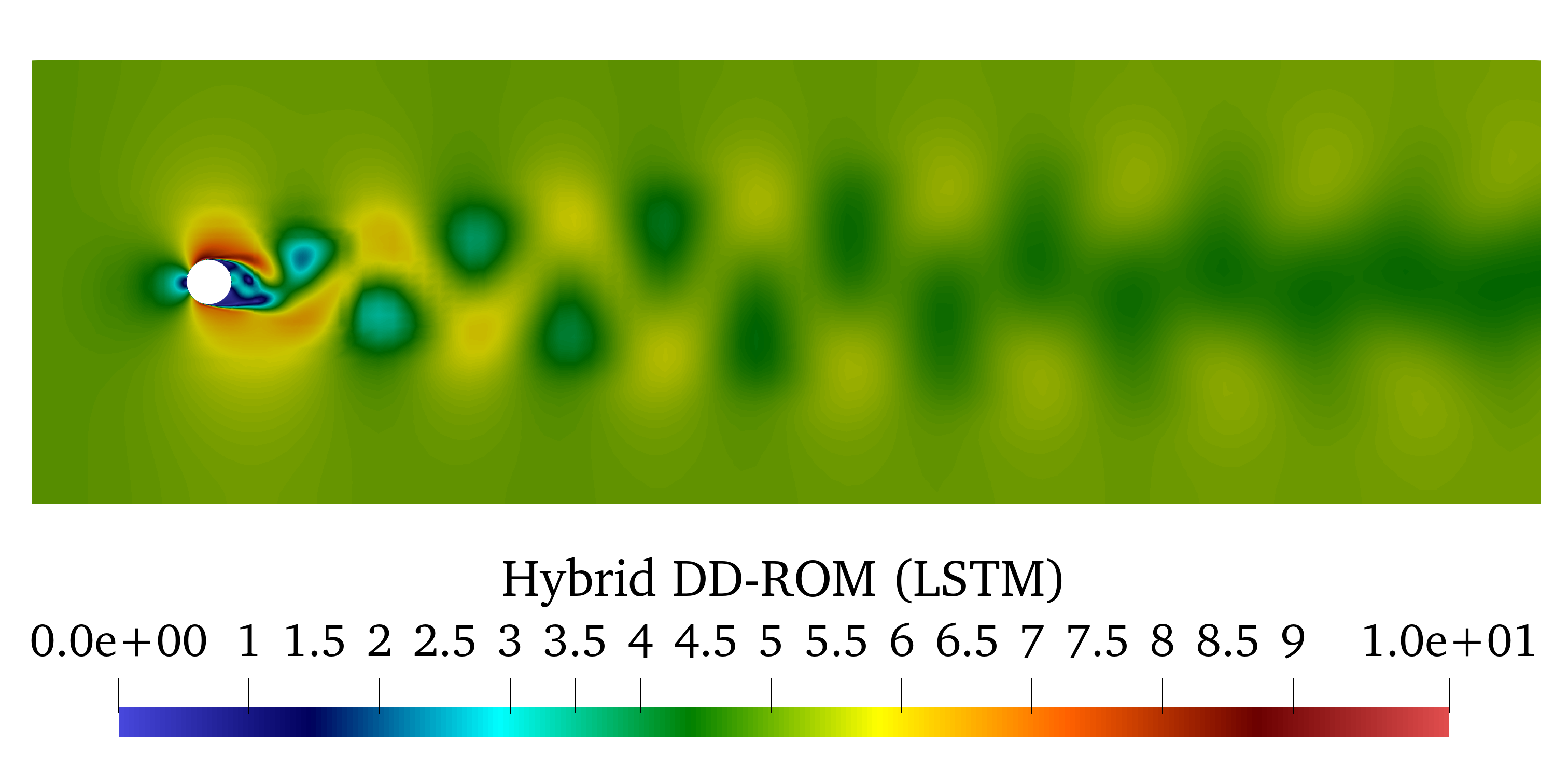}}
    \subfloat[]
    {\includegraphics[width=0.5\textwidth]{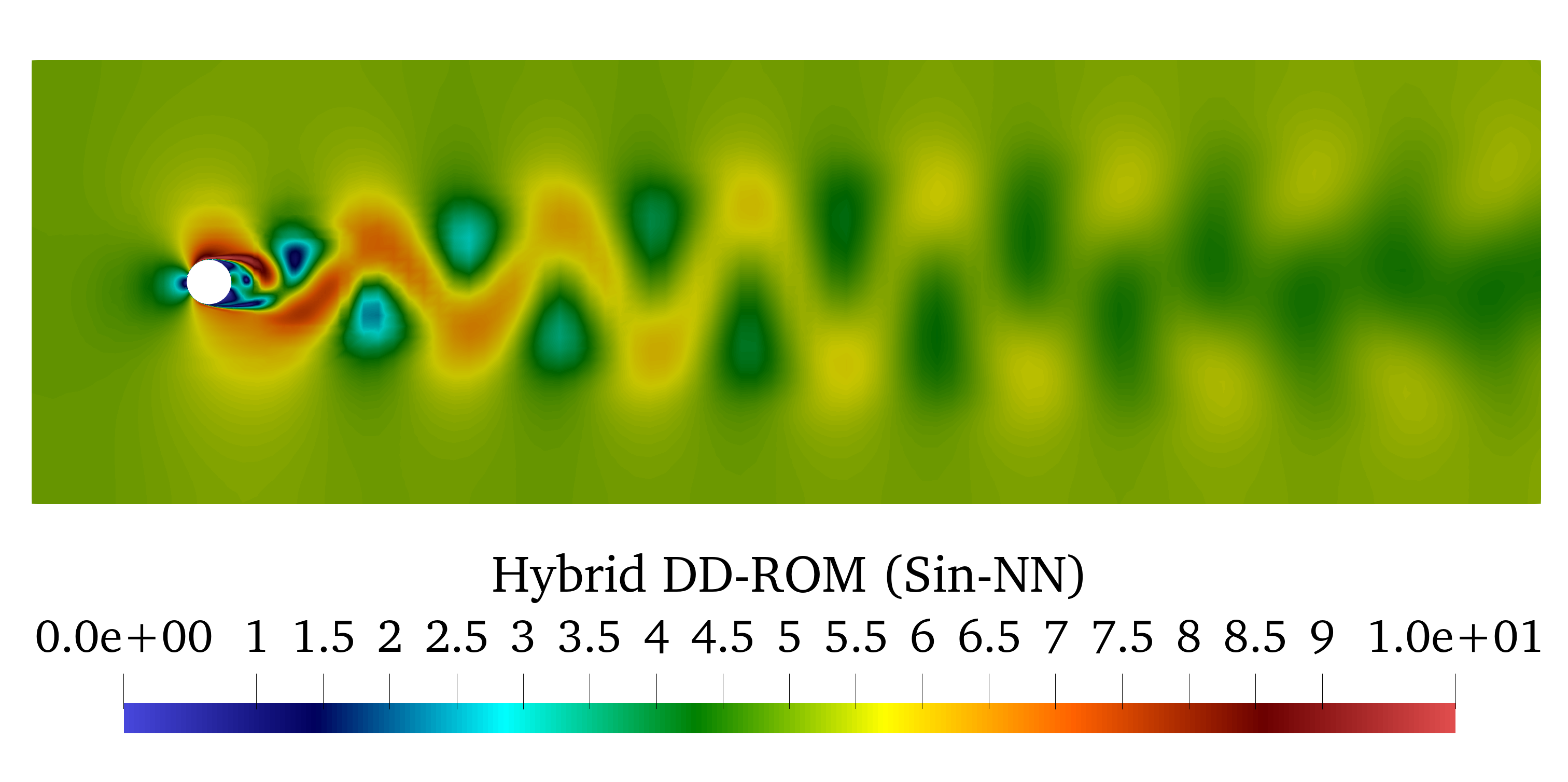}}
    \caption{Velocity magnitude fields at final instant of ROM.}
    \label{fig:vel-cyl-graph}
\end{figure}

\begin{figure}[htpb!]
    \centering
    \subfloat[]
    {\includegraphics[width=0.5\textwidth]{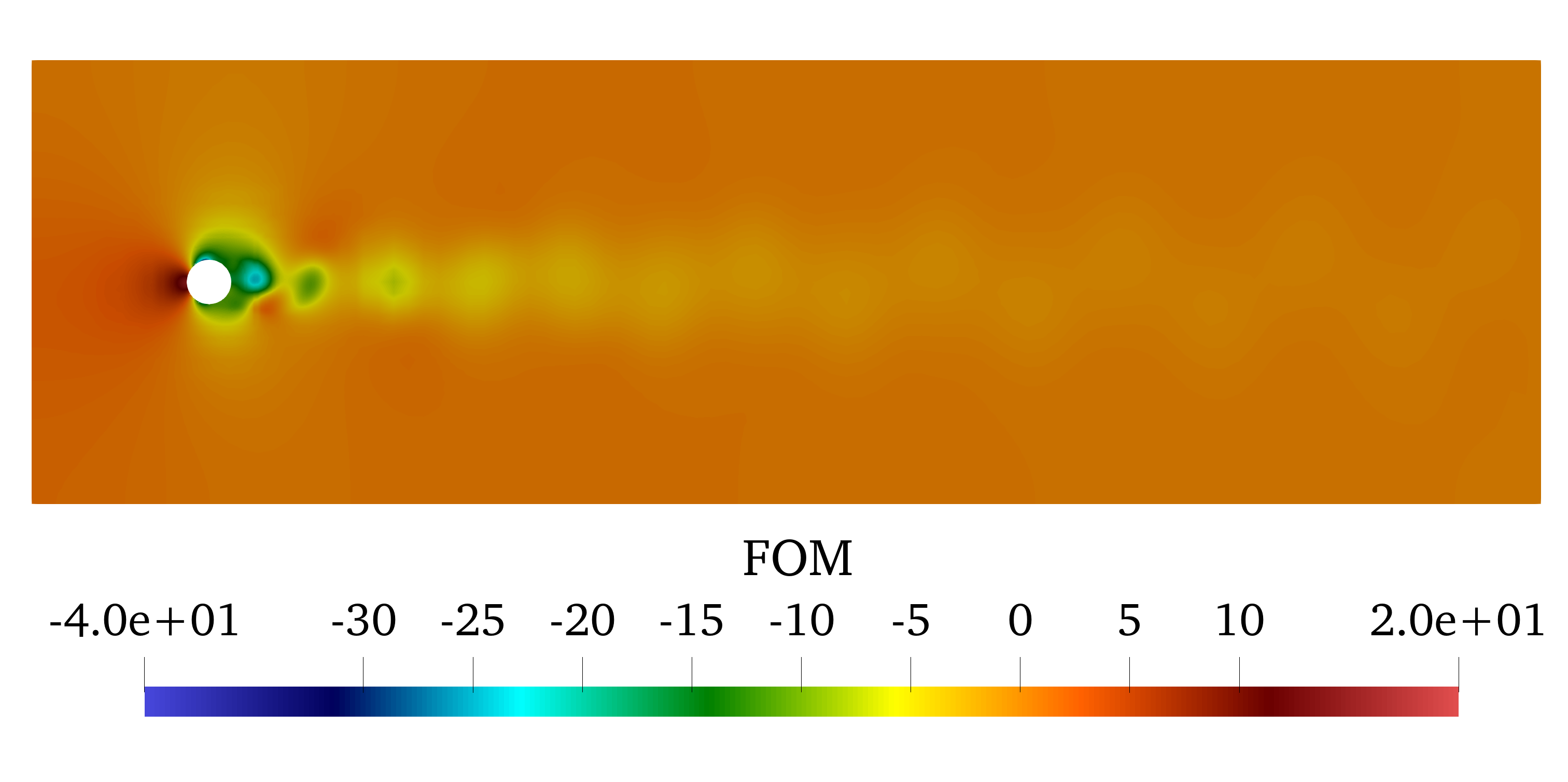}}\\
    \subfloat[]
    {\includegraphics[width=0.5\textwidth]{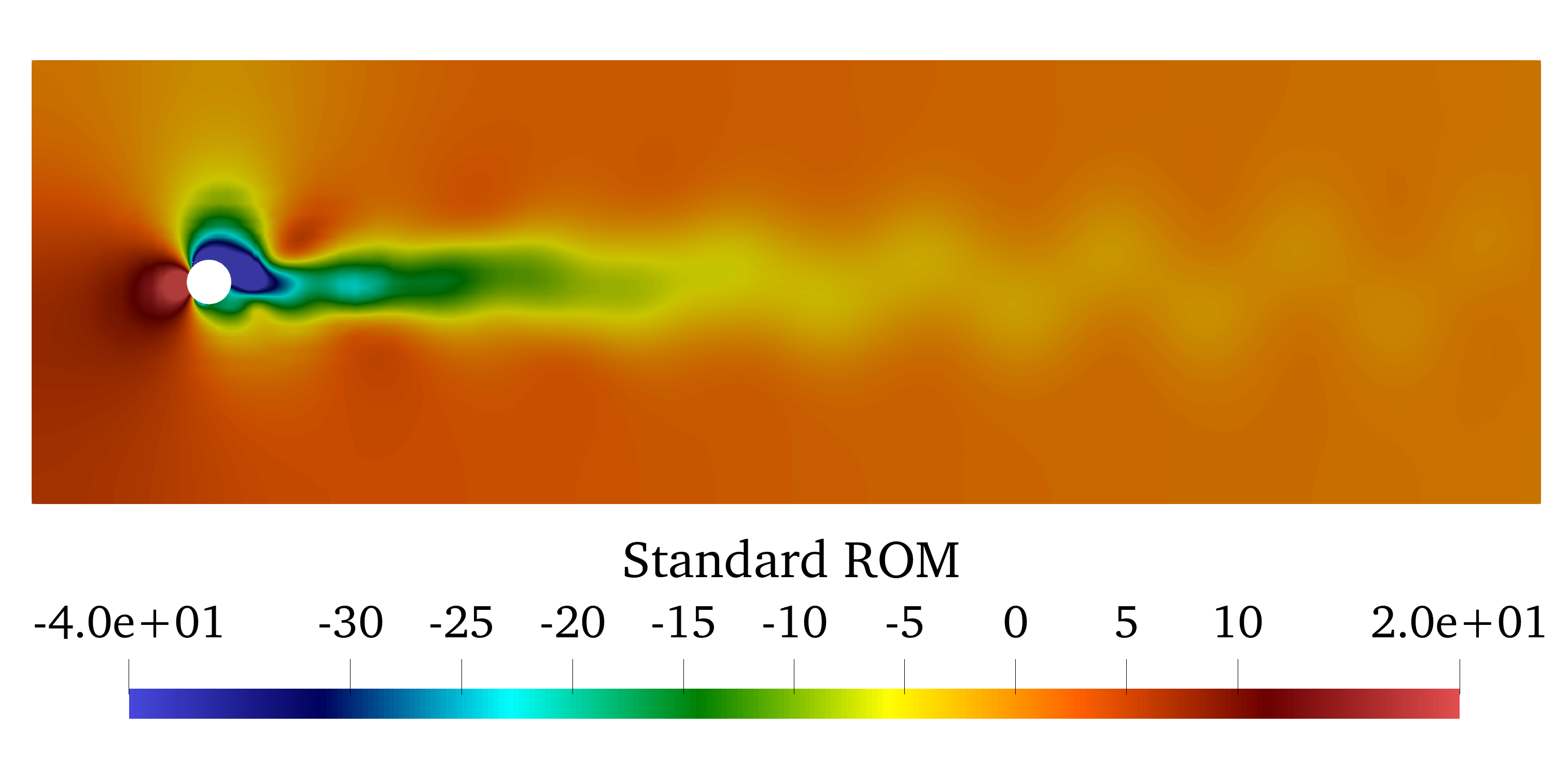}}
    \subfloat[]
    {\includegraphics[width=0.5\textwidth]{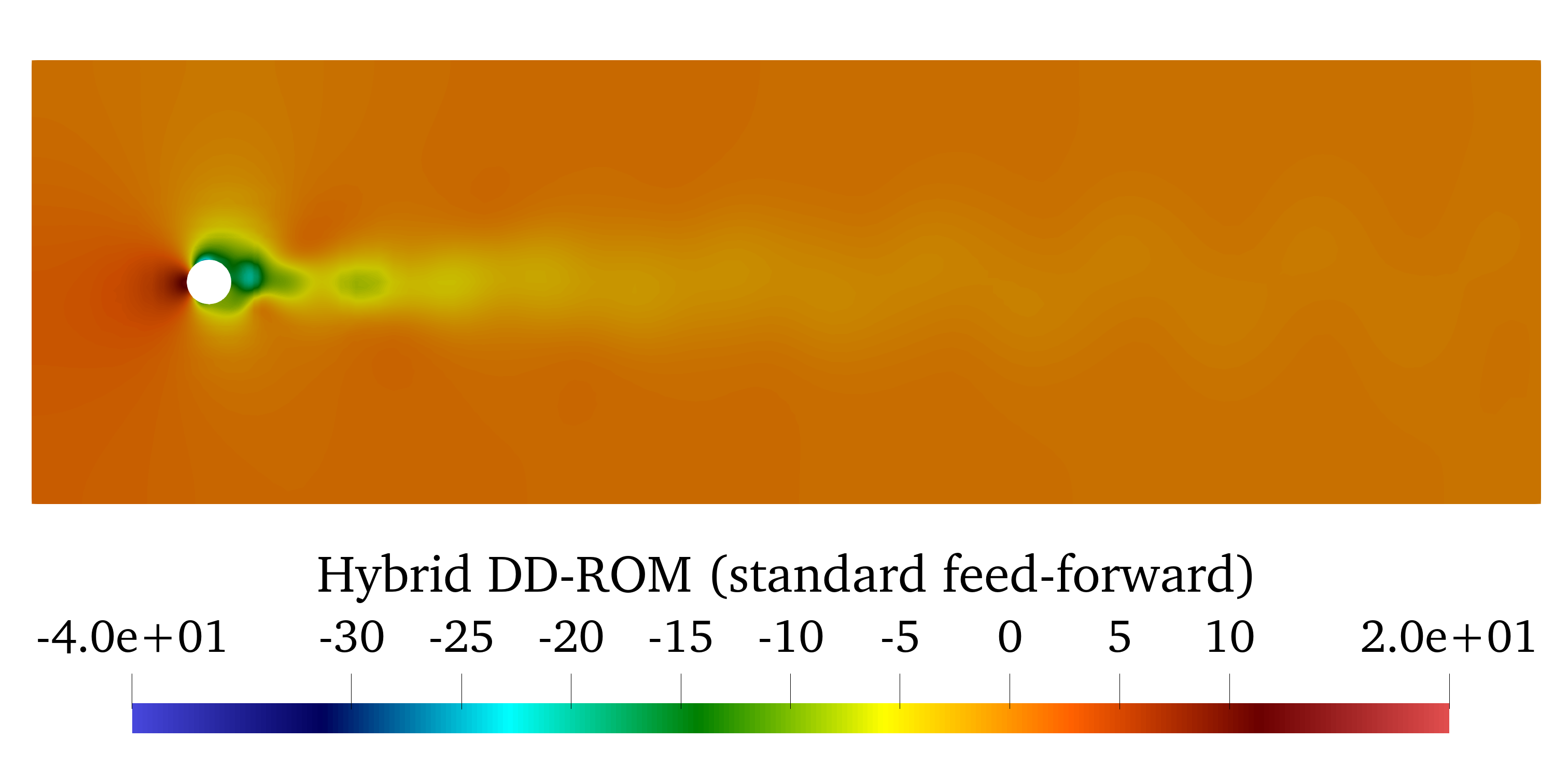}}\\
    \subfloat[]
    {\includegraphics[width=0.5\textwidth]{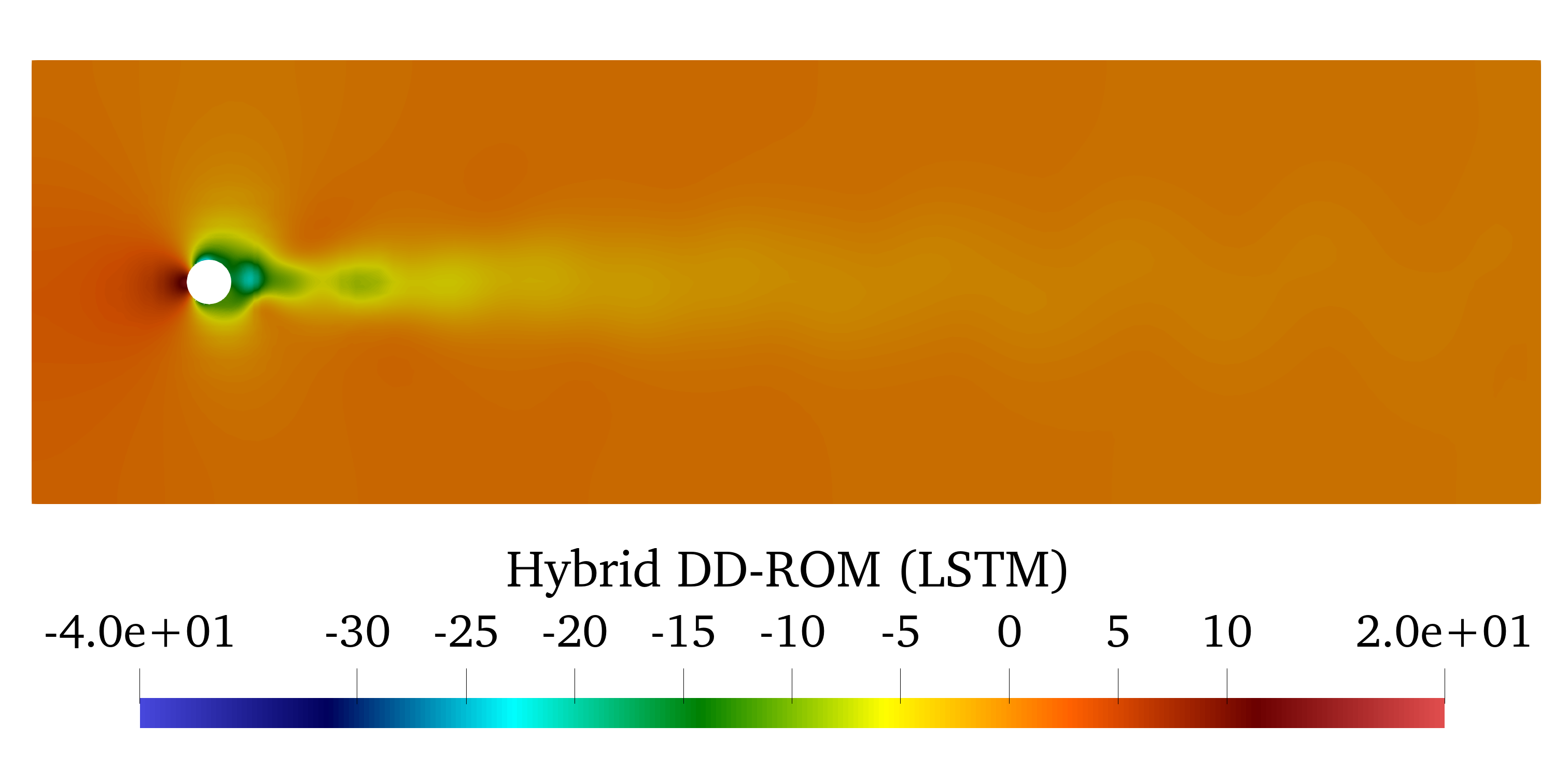}}
    \subfloat[]
    {\includegraphics[width=0.5\textwidth]{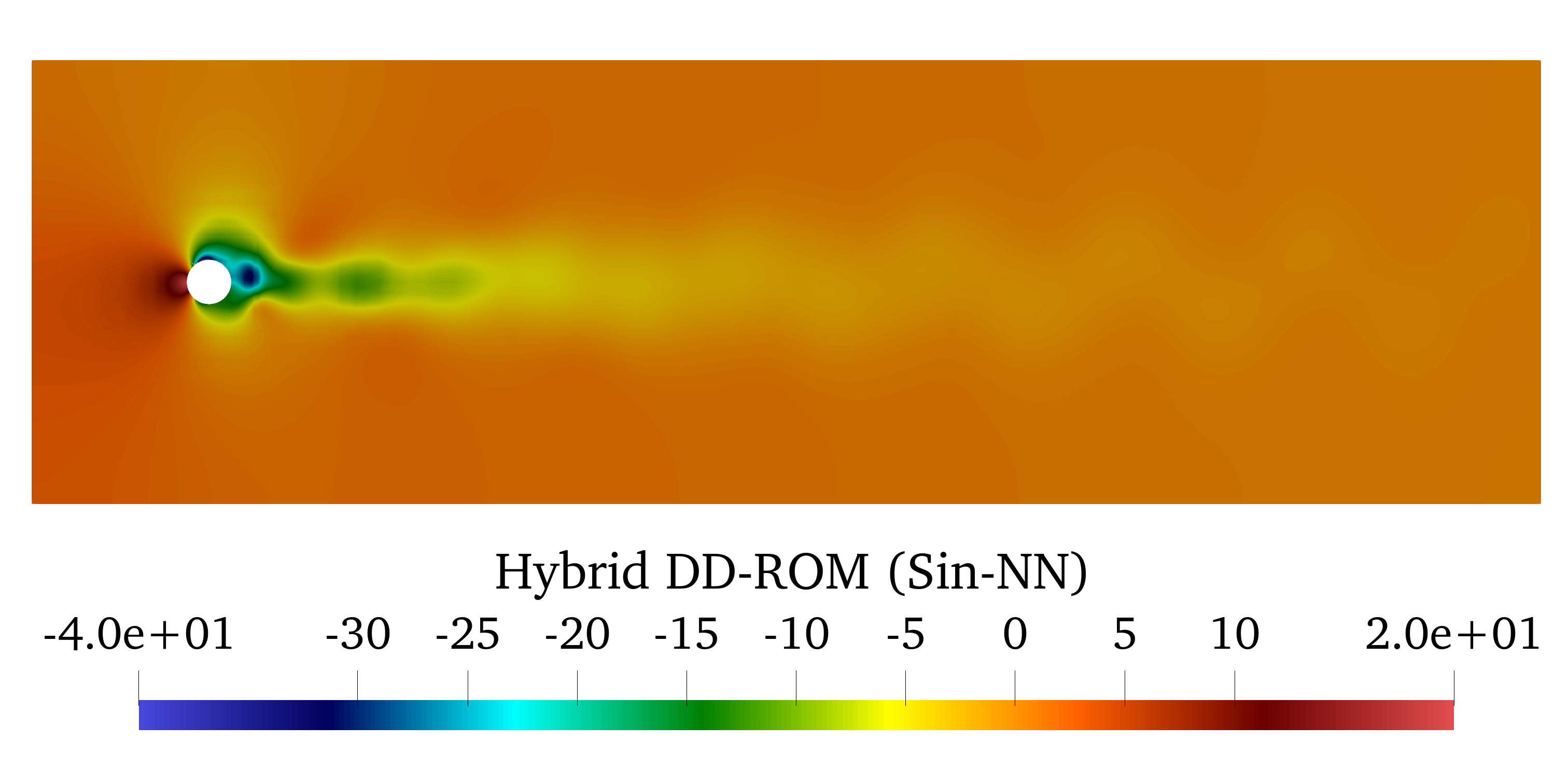}}
    \caption{Pressure fields at final instant of ROM.}
    \label{fig:pres-cyl-graph}
\end{figure}

Figures \ref{fig:vel-cyl-graph} and \ref{fig:pres-cyl-graph} show a graphical comparison of the velocity and pressure fields, respectively, at the final instance of the online simulation.

The POD-Galerkin ROM reconstruction leads to poor approximations in the region near the cylinder, for both the pressure and the velocity field. This may result in inaccurate predictions of the force fields, like the drag and the lift acting on the cylinder.
On the other hand, as pointed out in Section \ref{subsubsec:dd-roms-a}, the hybrid models are significantly closer to the full-order solution. However, the SinNN is less accurate around the cylinder, as previously shown in the errors' plots.

\newpage

\subsection{Test case \textbf{b}: flow in a channel-driven cavity}
\label{subsec:test-case-b}
The unsteady case of the flow inside a cavity-shaped channel is characterized by a different behavior to the first test case described in \ref{subsec:test-case-a}. Also in this case, we consider a parametrized setup w.r.t. the Reynolds number.

\subsubsection{Offline stage}
\label{subsubsec:fom-b}

\begin{figure}[htpb!]
    \centering
    \subfloat[Domain with notation]{\includegraphics[width=0.5\textwidth]{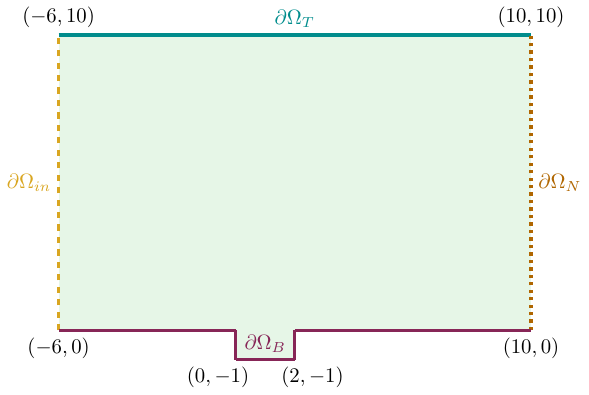}}
    \subfloat[Full order mesh]{\includegraphics[width=0.5\textwidth, trim={4.5cm 0 4.5cm 0cm}, clip]{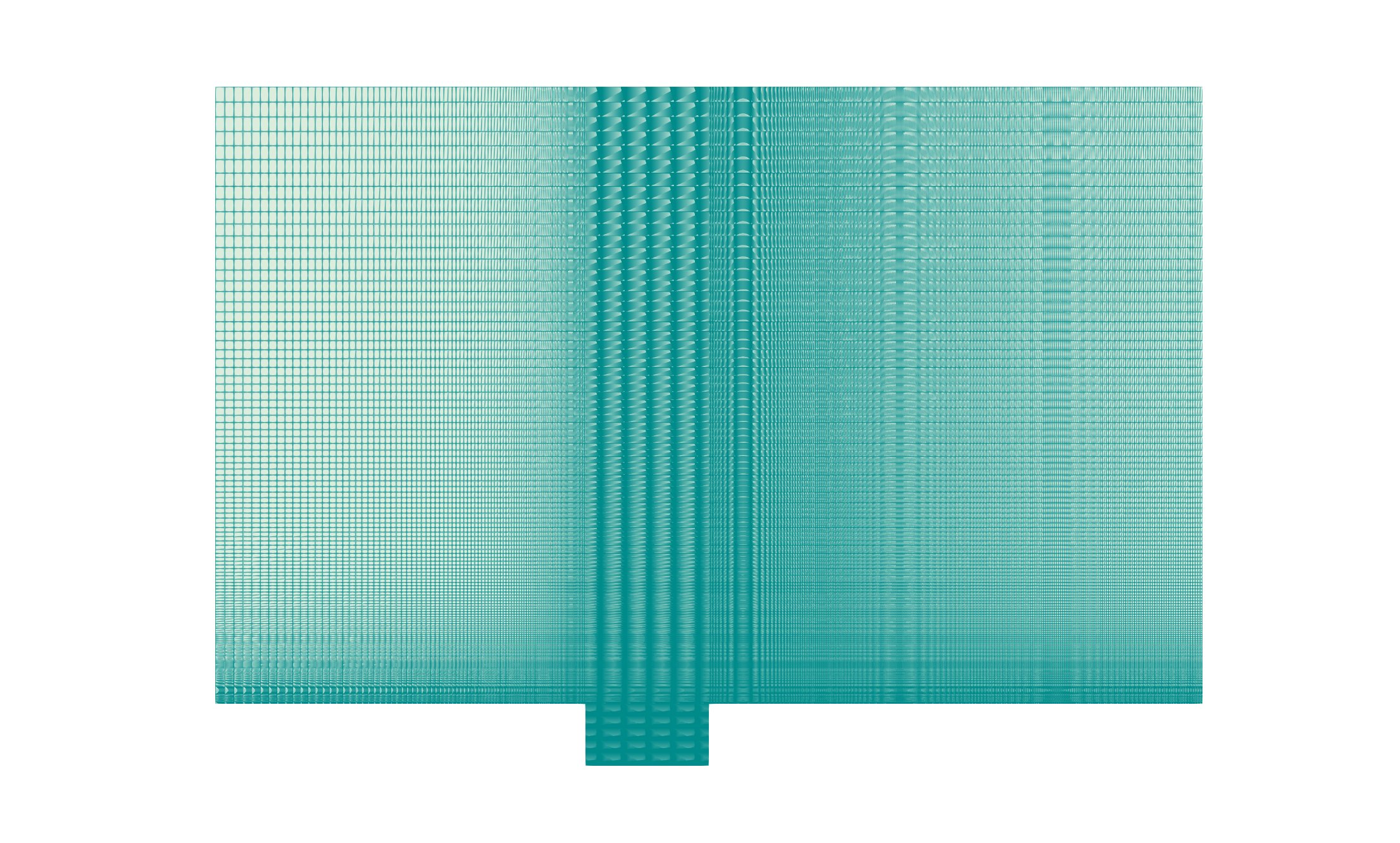}}
    \caption{The domain and full order mesh considered for the channel-driven flow test case.}
    \label{fig:cavity-domain}
\end{figure}

In Figure \ref{fig:cavity-domain} we represent the domain and the mesh considered in the FOM analysis.

Employing the notation in Figure \ref{fig:cavity-domain}, we consider the following boundary conditions: 
\begin{equation*}
    \text{On }\partial \Omega_T:
    \begin{cases}
        \nabla \bm{u} \cdot \bm{n} = 0,\\
        \nabla p \cdot \bm{n} = 0;
    \end{cases}
    \quad
        \text{On }\partial \Omega_B:
    \begin{cases}
       \bm{u}  = \bm{0},\\
       \nabla p \cdot \bm{n} = 0
        
    \end{cases}
\end{equation*}

\begin{equation*}
\text{On }\partial \Omega_{in}:
    \begin{cases}
        \bm{u} = (U^{b}_{in}, 0),\\
        \nabla p \cdot \bm{n} = 0;
    \end{cases}
    \quad
            \text{On }\partial \Omega_N:
    \begin{cases}
       \nabla \bm{u} \cdot \bm{n} = 0,\\
         p = 0.
    \end{cases}
\end{equation*}

For the POD we consider the viscosity range specified in Figure \ref{fig:params-plot}, with a corresponding Reynolds number in range $[\num{1e3}, \num{2e4}]$.
The time snapshots are retained every $\num{0.05}$ seconds in the time interval of 10 seconds, for a total amount of snapshots of $1000$, namely $200$ for each FOM simulation.

\begin{figure}[htpb!]
    \centering
    \subfloat[POD cumulative eigenvalues]{\includegraphics[height=4.7cm]{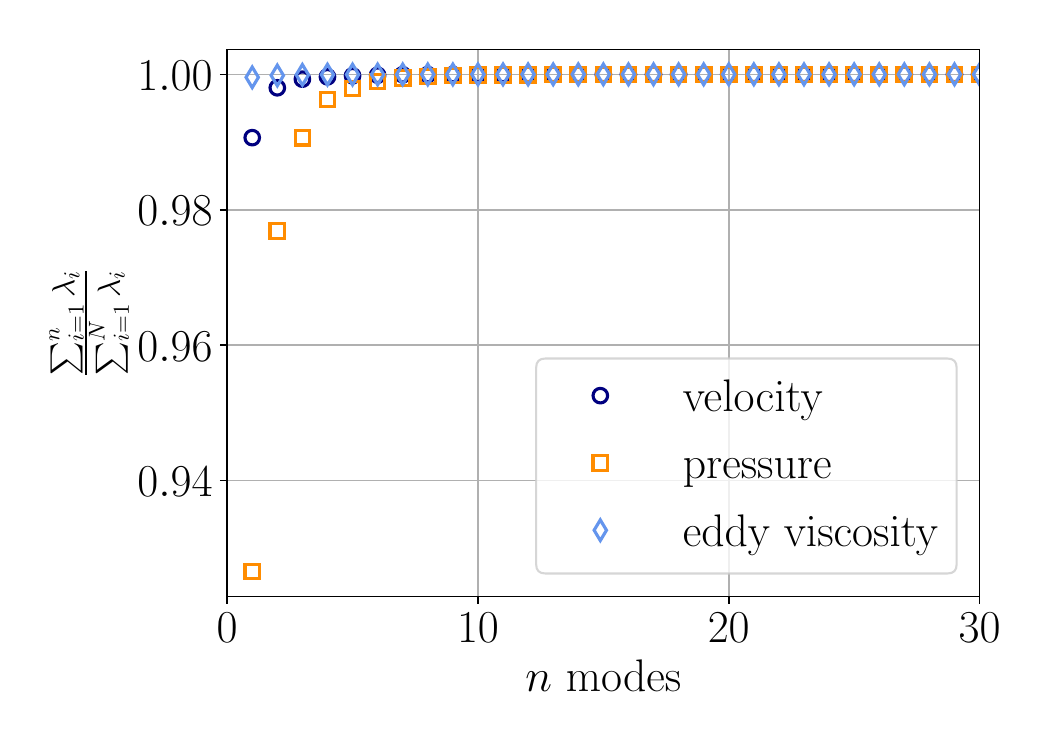}}
    \subfloat[POD eigenvalues decay]{\includegraphics[height=4.7cm]{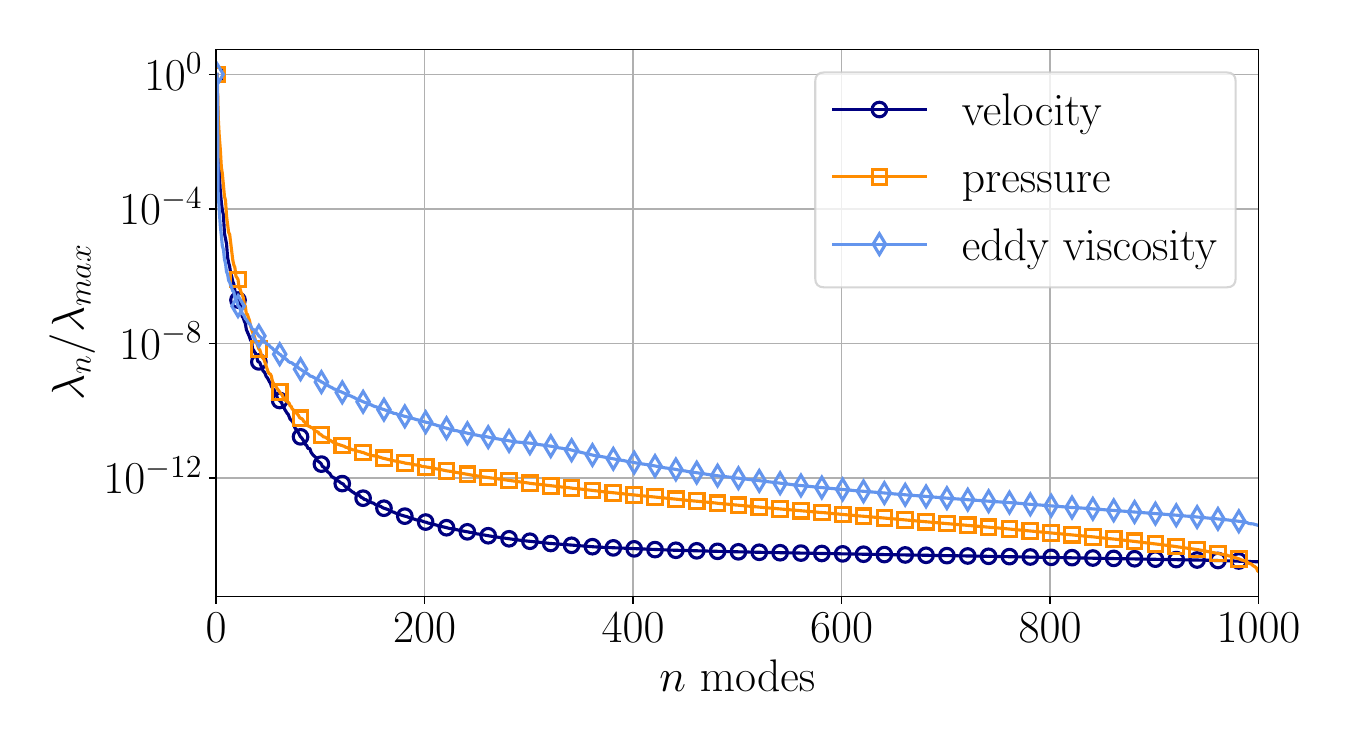}}
    \caption{Cumulative eigenvalues and eigenvalues decay for test case \textbf{b}.}
    \label{fig:eig-b}
\end{figure}

Figure \ref{fig:eig-b} shows the cumulative eigenvalues and the eigenvalues decay of the POD for all the fields considered. The Figure shows a slower decay than in the first test case (Figure \ref{fig:eig-a}). Indeed, this second test case is not periodic leading to a more challenging setting.

We decide to focus our analysis on a number of modes $N_u=N_{\nu_t}=r=2$ and $N_p=q=4$, namely in the \emph{marginally-resolved} regime.

\subsubsection{Neural networks' performance}
\label{subsec:networks-b}
As already done in Subsection \ref{subsec:networks-a}, the goal of this part is to show the performance of the neural networks used to predict both the \emph{eddy viscosity} coefficients and the \emph{correction} coefficients.

\begin{figure}[htpb!]
    \centering
    \includegraphics[width=\textwidth]{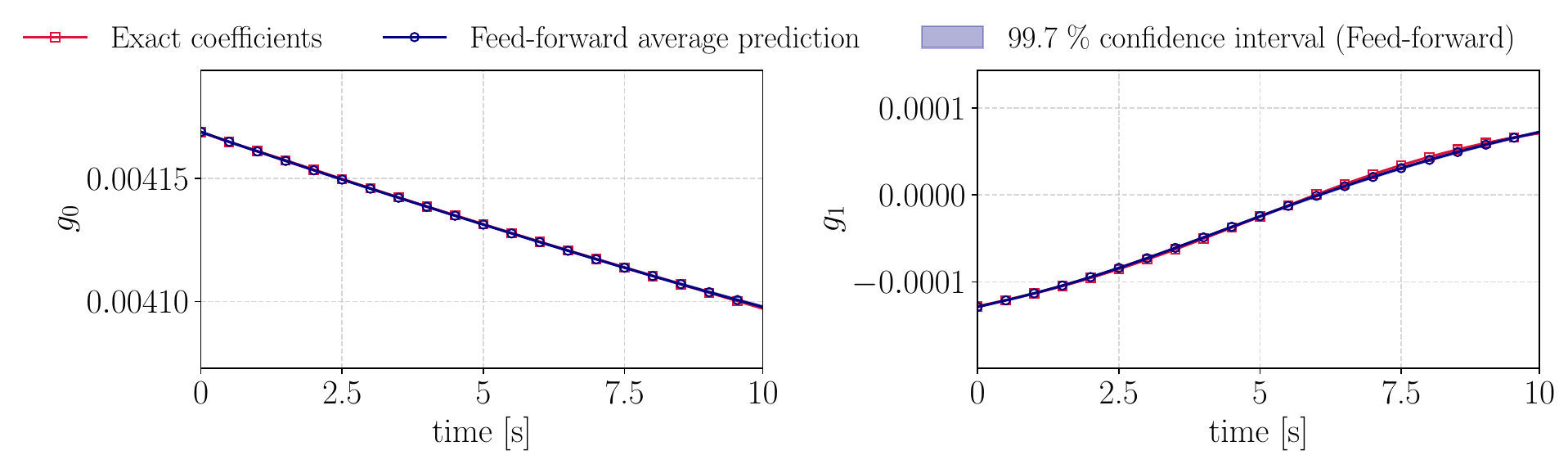}
    \caption{Eddy viscosity coefficients for $\nu=\SI{5e-6}{\metre^2 \per \second}$ (in $\nutrain$): neural networks predictions' average and POD projected exact coefficients.}
    \label{fig:cav-coeffs-eddy-0}
\end{figure}

In particular, the first mapping $\mathcal{G}(\bm{a}, \nu, t)$ is built, as in the first test case, with a multi-layer perceptron architecture. 
Figure \ref{fig:cav-coeffs-eddy-0} represents the performance of the machine learning model considered in a training setup. The resulting average prediction is accurate and the confidence interval is very small in both the coefficients $g_0$ and $g_1$.
Moreover, we notice that the trend of the coefficients is not periodic, but it is monotone, differently from test case \textbf{a}.

For what concerns the mapping $\mathcal{M}$ we only take into account the pressure corrections $\bm{\tau}_p$, since the velocity approximation has not a relevant contribution in this case. For this reason, in the networks we will not take into account the velocity coefficients $\bm{a}$ as input.

\begin{figure}[htpb!]
    \centering
    \includegraphics[width=\textwidth]{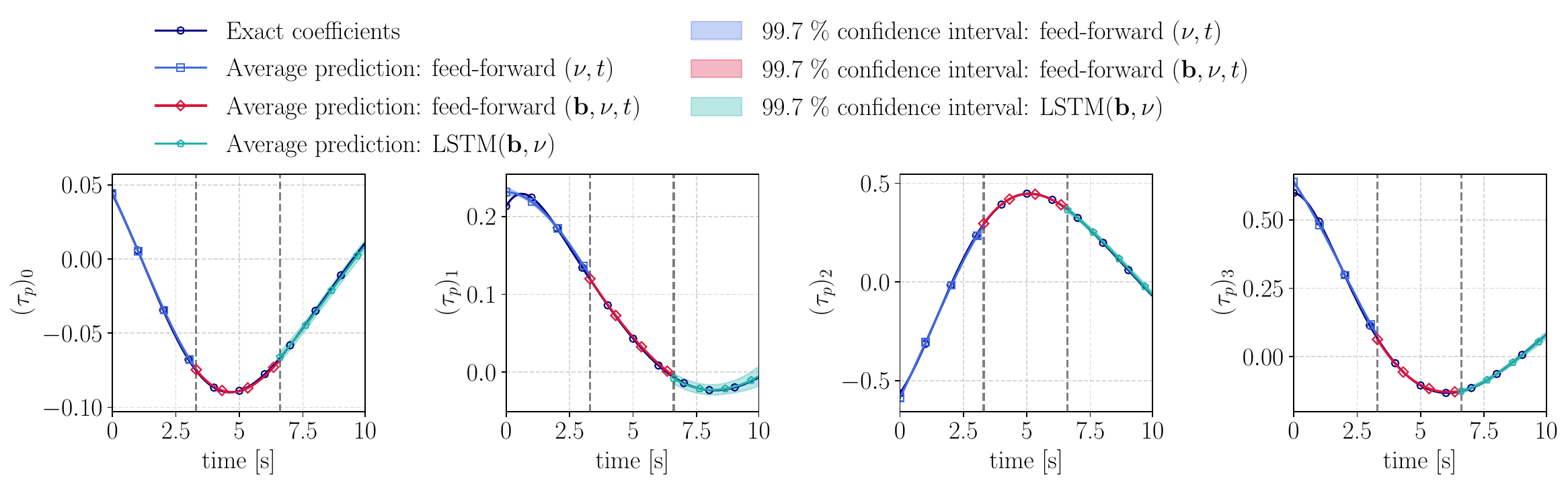}
    \caption{Correction terms' coefficients for $\nu=\SI{5e-6}{\metre^2 \per \second}$ (in $\nutrain$): neural networks predictions' average and POD projected exact coefficients.}
    \label{fig:cav-coeffs-turb-0}
\end{figure}

Figure \ref{fig:cav-coeffs-turb-0} shows
the performance of three different types of neural networks for mapping $\mathcal{M}$. In particular, we named the networks as:
\begin{itemize}
    \item Feed-forward$(\bm{b}, \nu, t)$, where $(\bm{b}, \nu, t)$ are the inputs, with a multi-layer perceptron standard architecture;
    \item LSTM$(\bm{b}, \nu)$, where the time parameter is not considered since it is already embedded in the Long--Short Term Memory framework;
    \item Feed-forward$(\nu, t)$, where $(\nu, t)$ are the only inputs. We chose to test this model for the simple monotone trend in time of the coefficients.
\end{itemize}

The sequence length used for the LSTM, i.e. the entity of time memory, is $20$ time steps, corresponding to $1$ second of simulation.

\subsubsection{DD-ROMs performance}
\label{subsubsec:dd-roms-b}
We decide in this Section to focus only on the hybrid DD-ROM since the purely DD-ROM has not significantly improved the results in test case \textbf{a}.

\begin{figure}[htpb!]
    \centering
    \subfloat[]{\includegraphics[width=\textwidth]{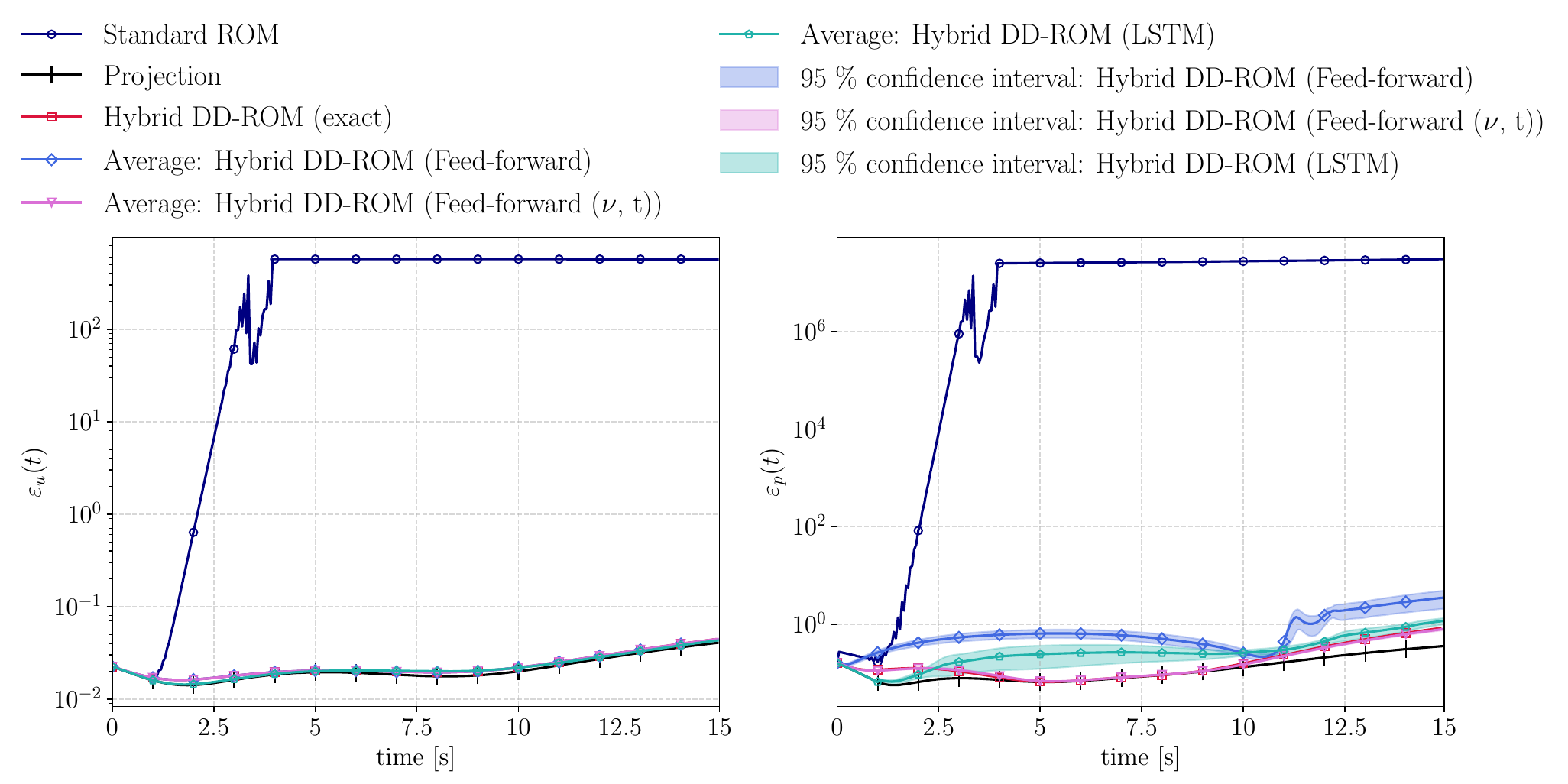}}\\
    \subfloat[]{\includegraphics[width=0.9\textwidth, trim={0 0 0 5cm}, clip]{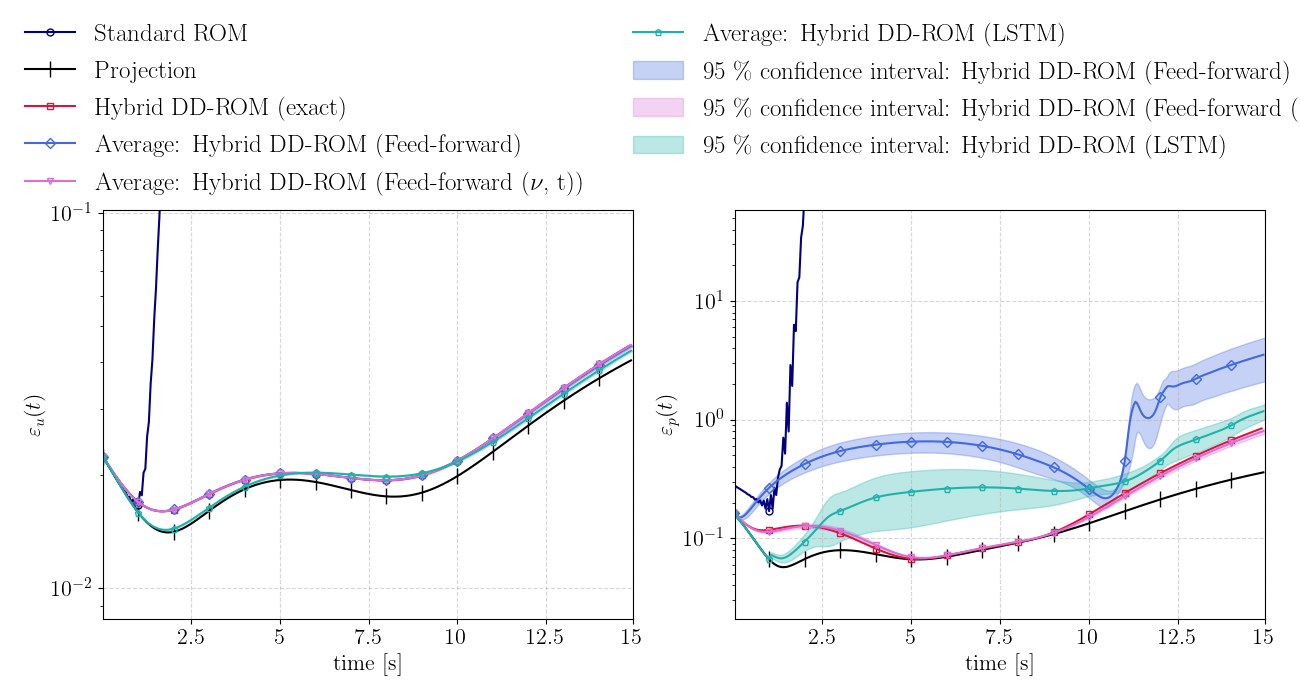}}
    \caption{Relative errors for $\nu \in \nutest$, $\nu=\SI{4e-6}{\metre^2 \per \second}$ for the pressure and the velocity magnitude fields, for the physics-based DD-ROMs, hybrid DD-ROMs, standard POD-Galerkin ROMs (first row). The second row represents a zoomed version of the error.}
    \label{fig:err-test-ext-low-cav}
\end{figure}

\begin{figure}[htpb!]
    \centering
    \subfloat[]{\includegraphics[width=\textwidth]{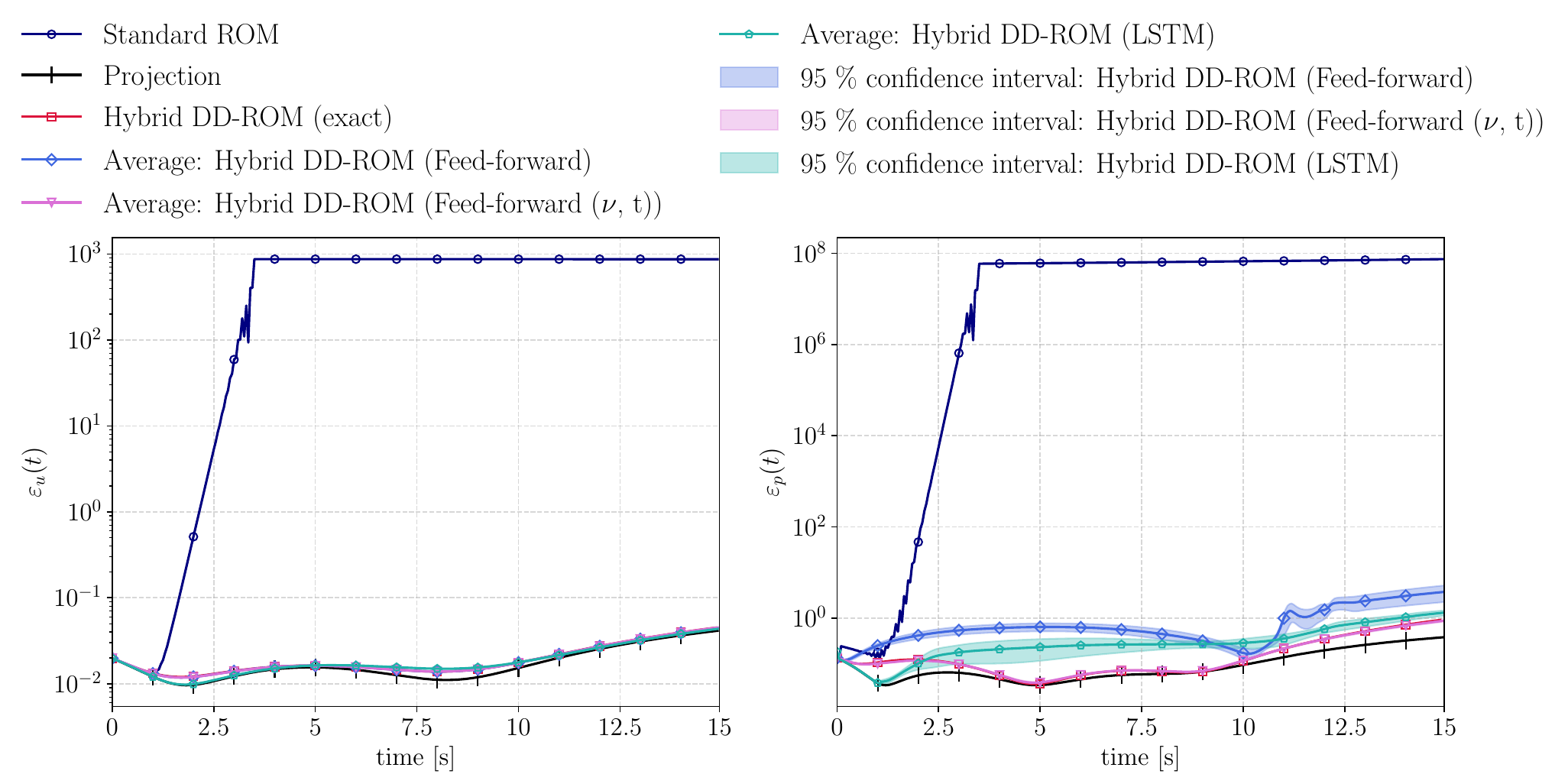}}\\
    \subfloat[]{\includegraphics[width=0.9\textwidth]{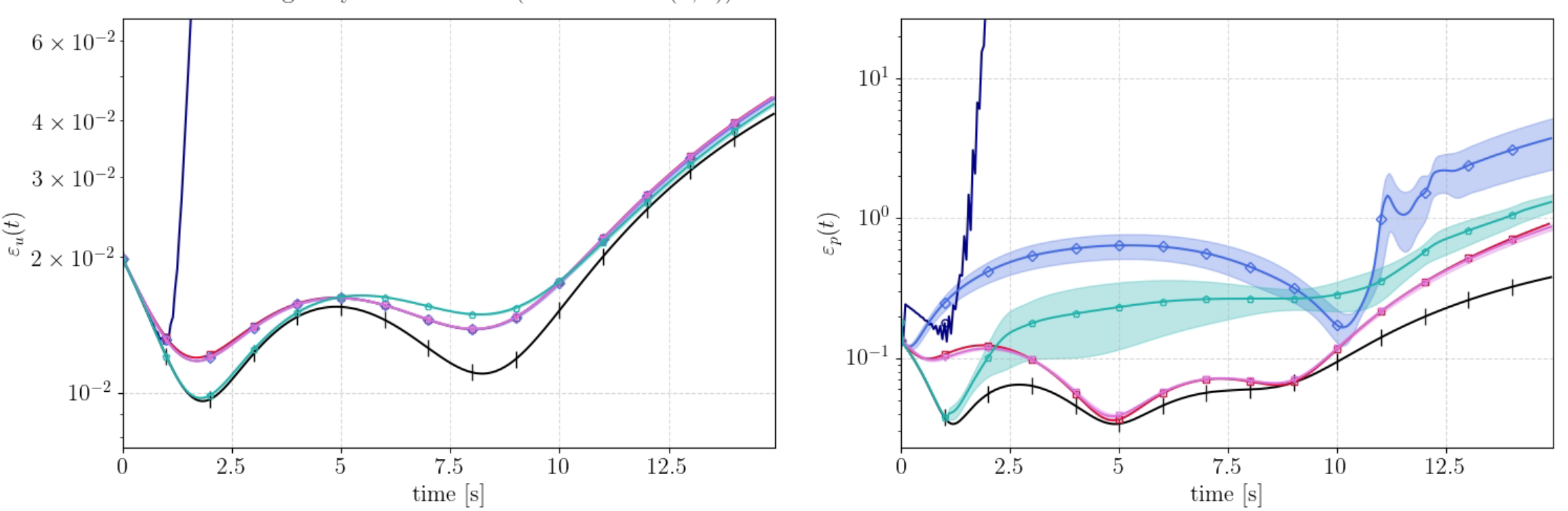}}
    \caption{Relative errors for $\nu \in \nutest$, $\nu=\SI{7e-6}{\metre^2 \per \second}$ for the pressure and the velocity magnitude fields, for the physics-based DD-ROMs, hybrid DD-ROMs, standard POD-Galerkin ROMs (first row). The second row represents a zoomed version of the error.}
    \label{fig:err-test-int-cav}
\end{figure}

Figures \ref{fig:err-test-ext-low-cav} and \ref{fig:err-test-int-cav} show that the relative errors of the velocity and pressure fields w.r.t. the FOM, for two test viscosity values. Moreover, the value for Figure \ref{fig:err-test-ext-low-cav} is outside of the POD range.

The standard ROM approach completely blows up after a few time steps in both cases, because of the ill-conditioning of the system. In this case, the correction terms are necessary not only to improve the accuracy but also to ensure the stability of the reduced formulation.

On the one hand, in Figure \ref{fig:err-test-ext-low-cav} the viscosity value is the lowest taken into account, namely corresponding to the highest Reynolds number. In this case, the performance of all the DD-ROMs leads to a stabilized system with improved accuracy.

In both cases, the best performance is obtained with the multi-layer perceptron only having as input the system parameters, while the worst performance is the feed-forward which also considers as input the pressure coefficients.
The average prediction obtained with the feed-forward$(\nu, t)$ network is similar to the one obtained by considering the exact coefficients, with higher confidence compared to the other neural networks.
On the contrary to what was observed for test case \textbf{a.}, the cavity presents a monotone trend in time, which appears to be well captured through a mapping only depending on parameters $(\nu, t)$. 

In conclusion, the hybrid DD-ROM approach may be used in reduced order systems both to improve the results and to avoid stability issues and ill-conditioning.

However, the method is still limited by the linearity of the POD, which is used to write system \eqref{eq:hyb-dd-rom}. Indeed, the correction terms added into the PPE-ROM system are obtained from the operators coming from a POD-Galerkin projection. Those operators in an advection-dominated test case and a marginally-resolved regime, may not well represent the dynamics of the system.

\begin{figure}[htpb!]
    \centering
    \includegraphics[width=0.6\textwidth]{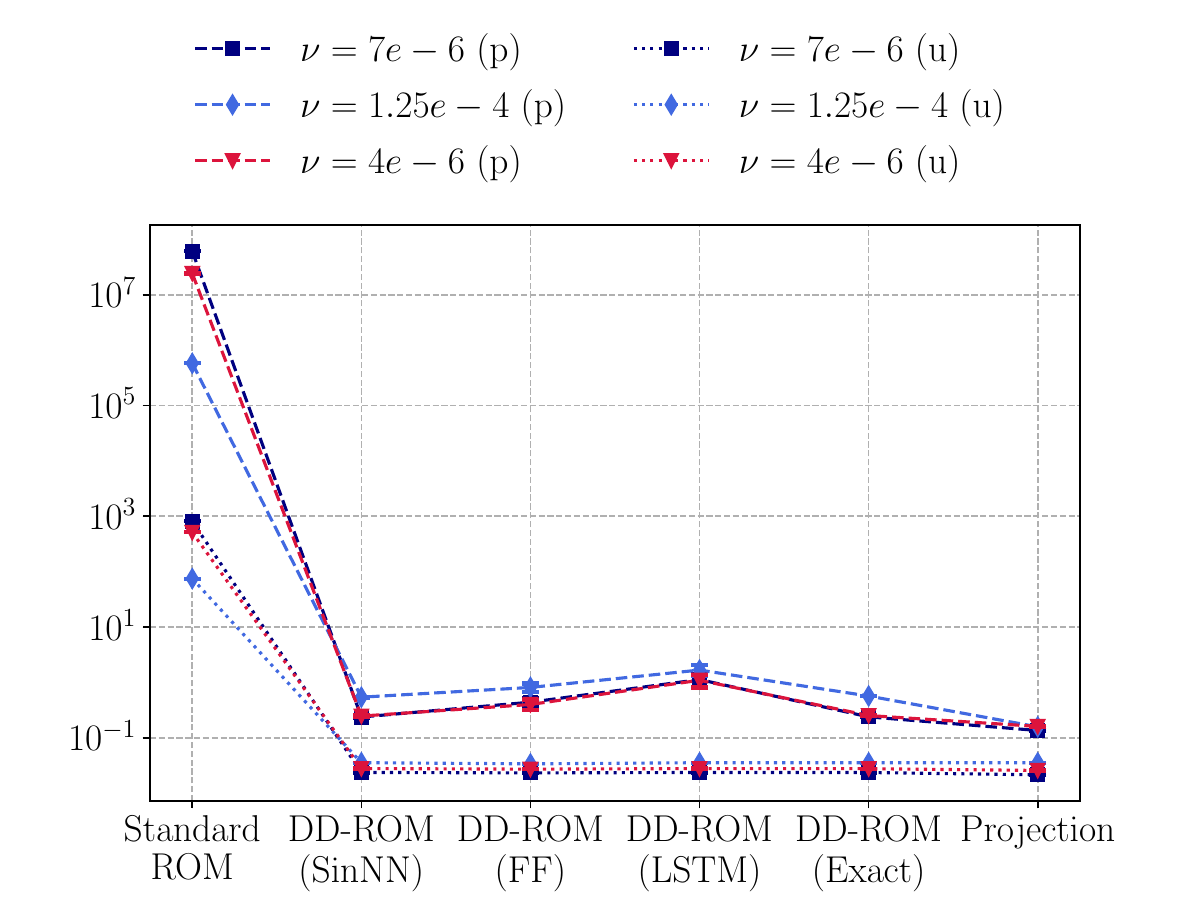}
    \caption{Integrals of relative errors over time for $\nu \in \nutest$ for the pressure and the velocity magnitude fields, in the standard ROM, hybrid DD-ROMs, and for the projection.}
    \label{fig:integrals-cavity}
\end{figure}

\subsubsection{Graphical results}
\label{subsubsec:graph-b}

\begin{figure}[htpb!]
    \centering
    \subfloat[]
    {\includegraphics[width=0.48\textwidth]{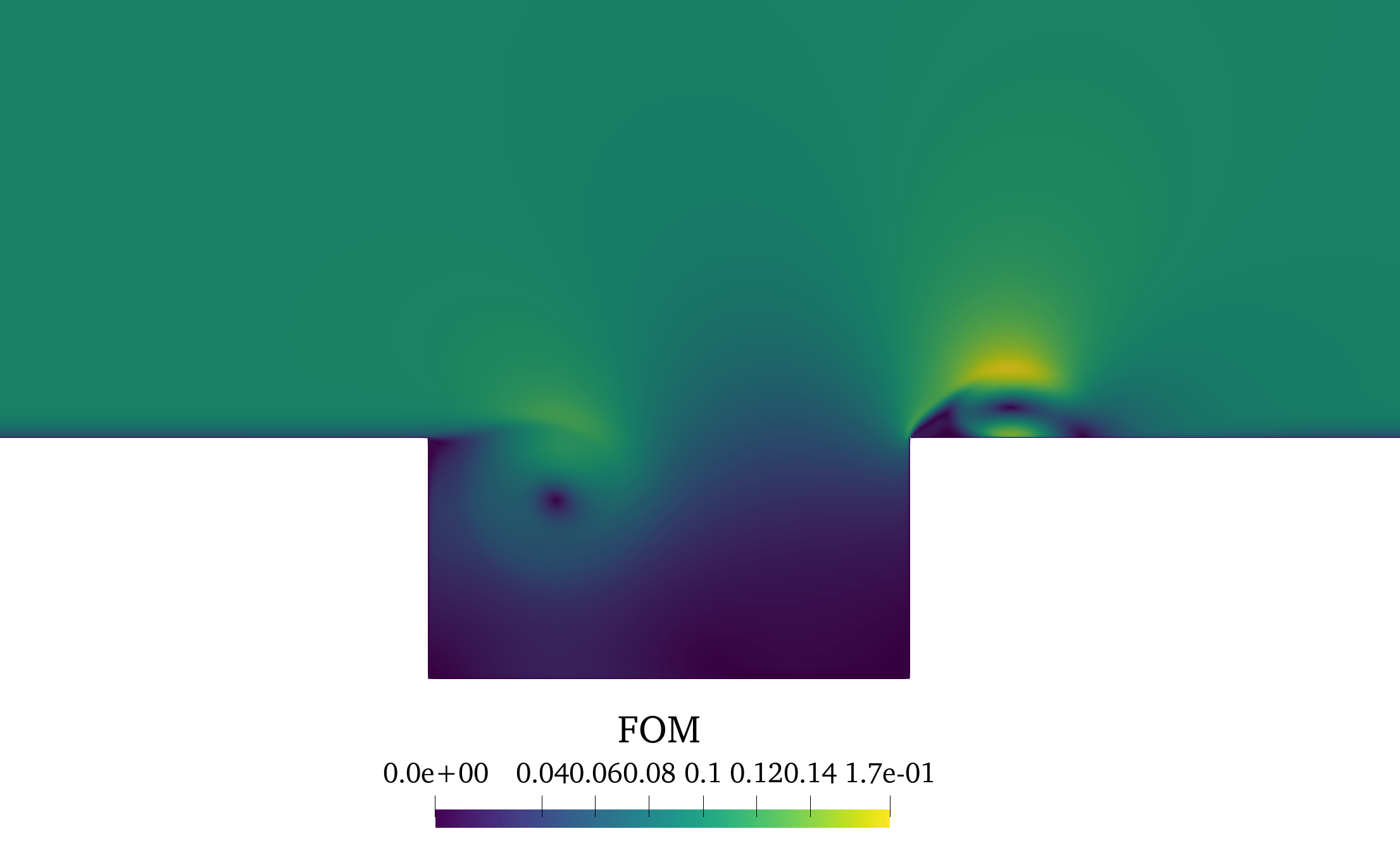}}
    \vspace{0.02\textwidth}
    \subfloat[]
    {\includegraphics[width=0.48\textwidth, trim={1.7cm 0cm 1.7cm 2cm}, clip]{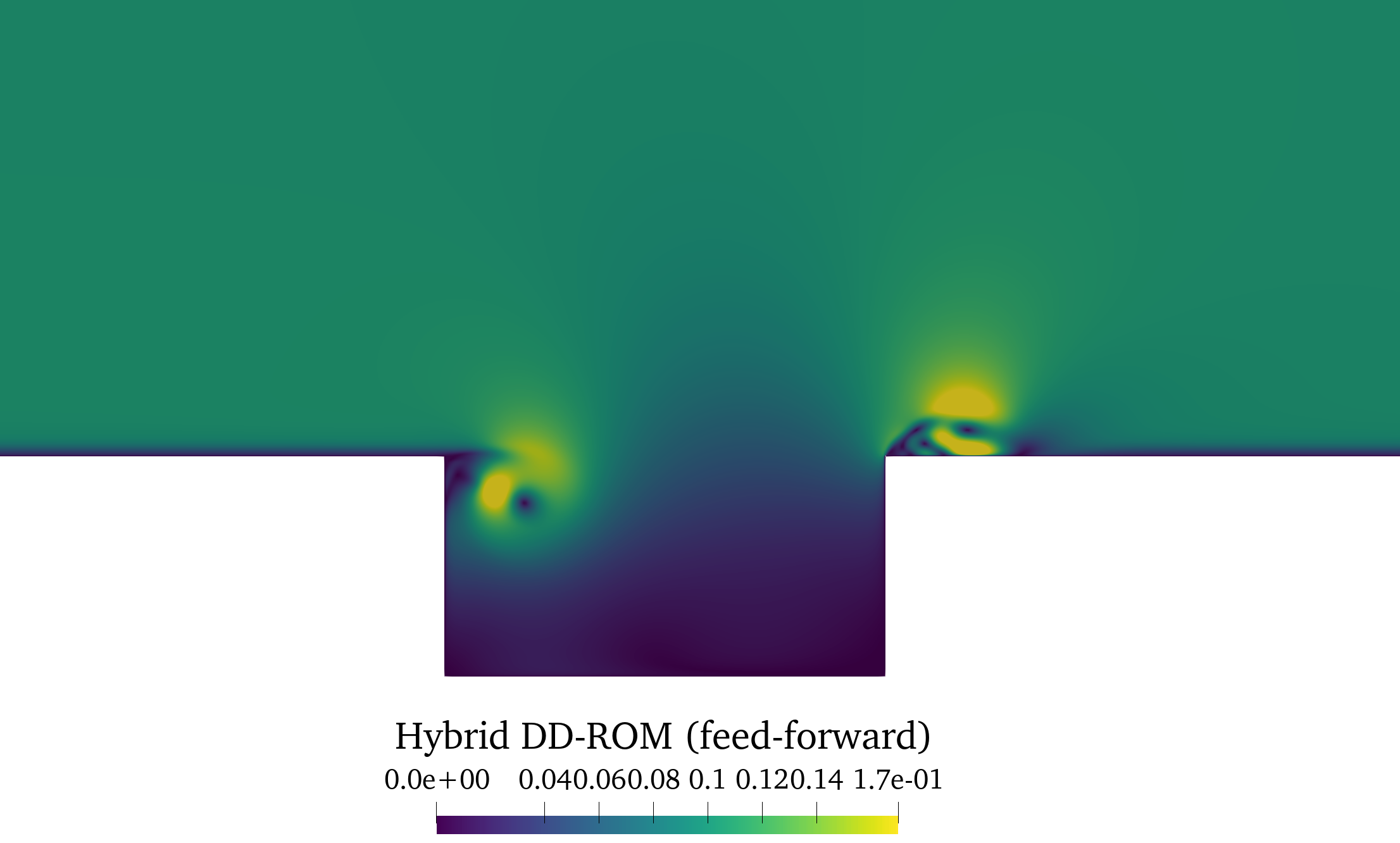}}\\
    \subfloat[]
    {\includegraphics[width=0.48\textwidth, trim={1.7cm 0cm 1.7cm 2cm}, clip]{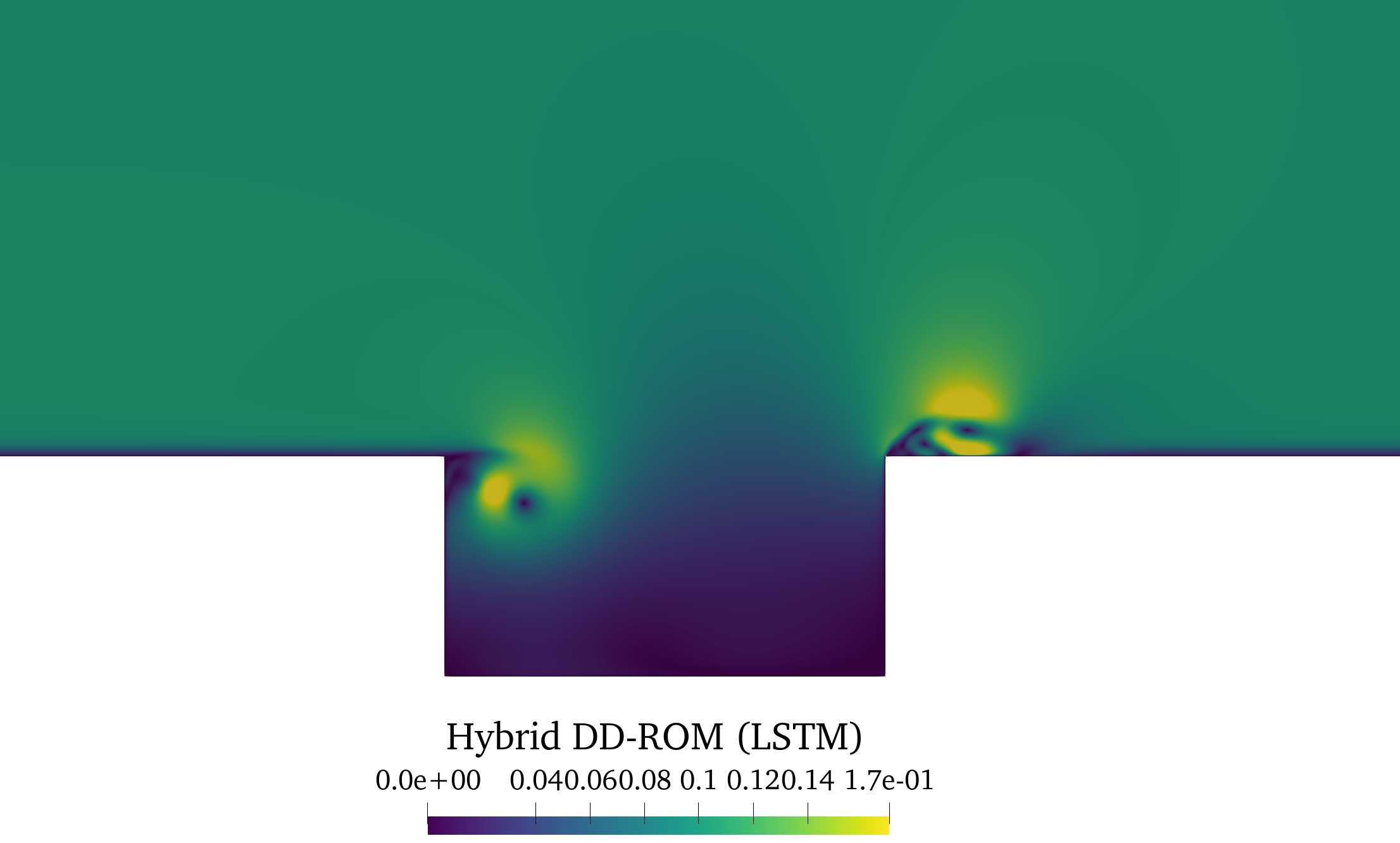}}
    \vspace{0.02\textwidth}
    \subfloat[]
    {\includegraphics[width=0.48\textwidth, trim={1.7cm 0cm 1.7cm 2cm}, clip]{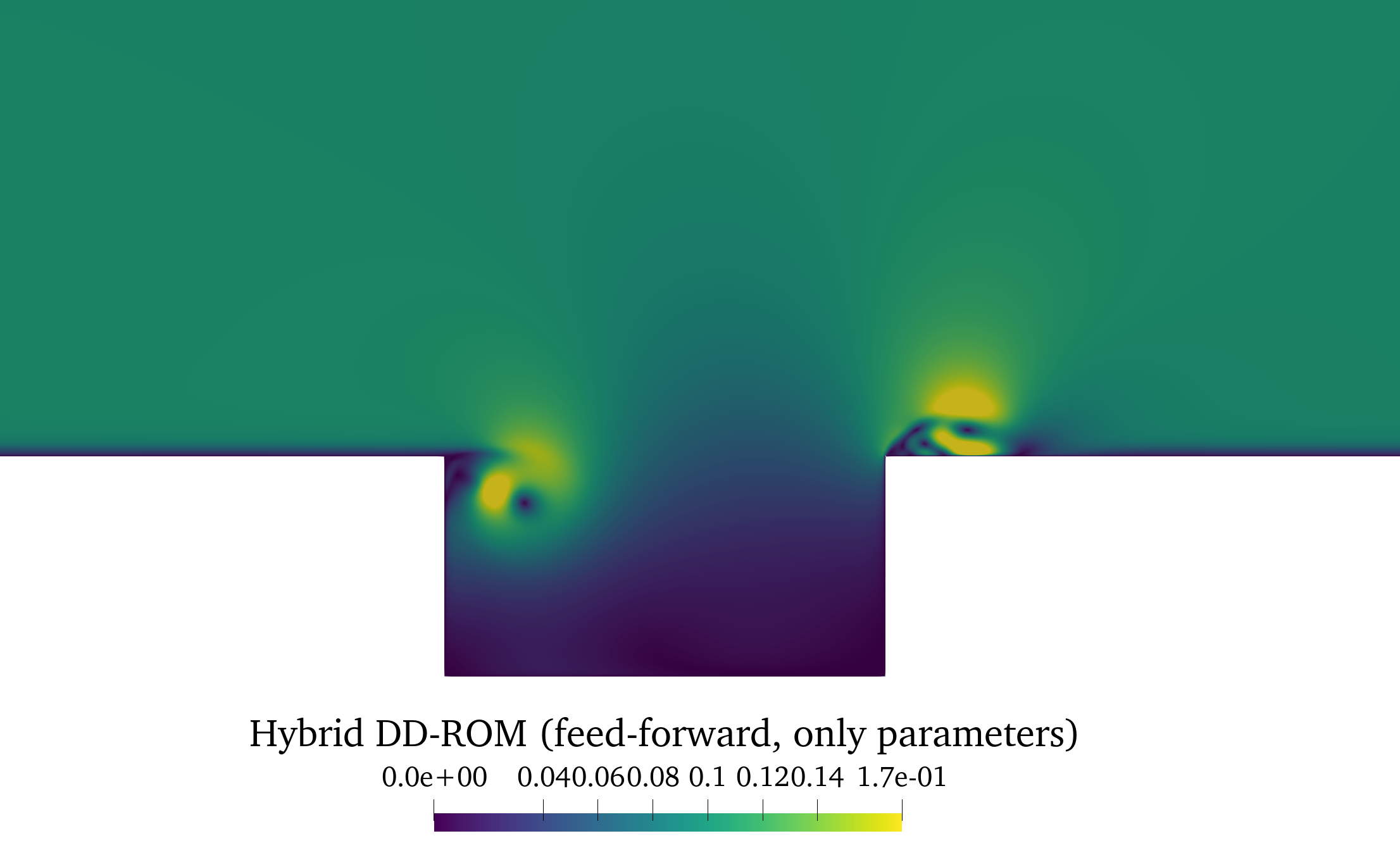}}
    \caption{Velocity magnitude fields at final instant of ROM. The standard ROM results are not represented here because the standard simulation blows up.}
    \label{fig:vel-cav-graph}
\end{figure}

\begin{figure}[htpb!]
    \centering
    \subfloat[]
    {\includegraphics[width=0.48\textwidth, trim={1.7cm 0cm 1.7cm 2cm}, clip]{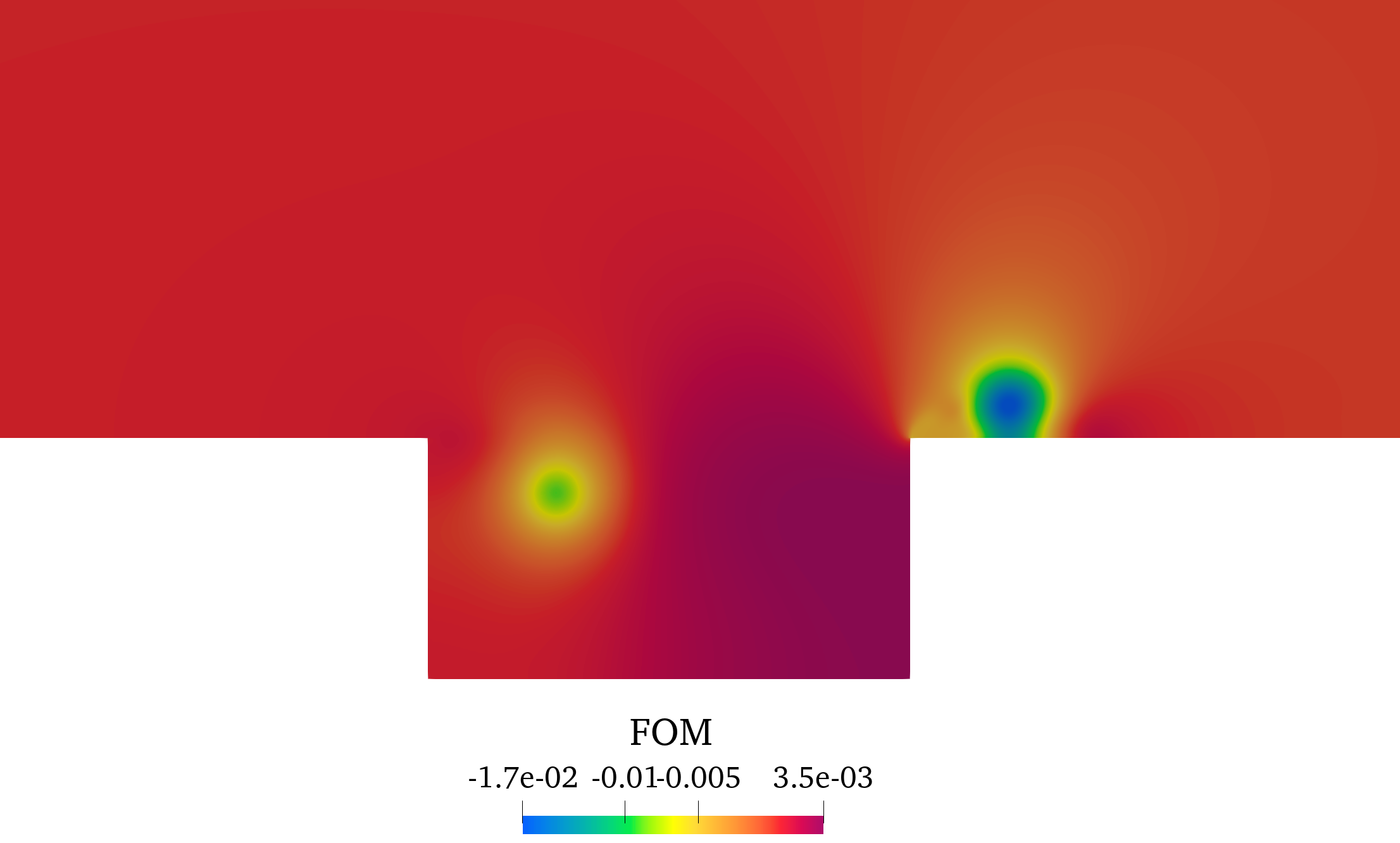}}
    \vspace{0.02\textwidth}
    \subfloat[]
    {\includegraphics[width=0.48\textwidth, trim={1.7cm 0cm 1.7cm 2cm}, clip]{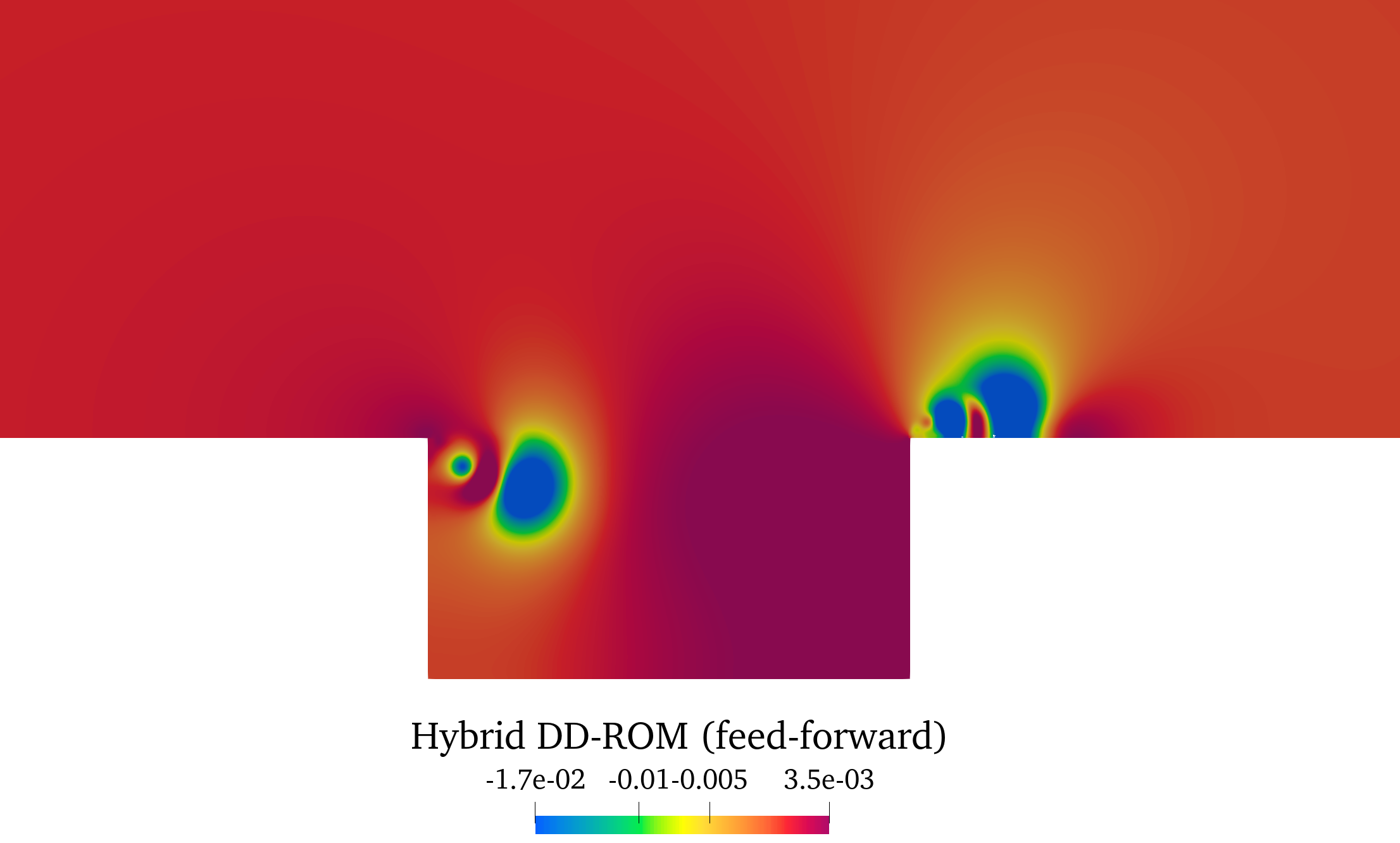}}\\
    \subfloat[]
    {\includegraphics[width=0.48\textwidth, trim={1.7cm 0cm 1.7cm 2cm}, clip]{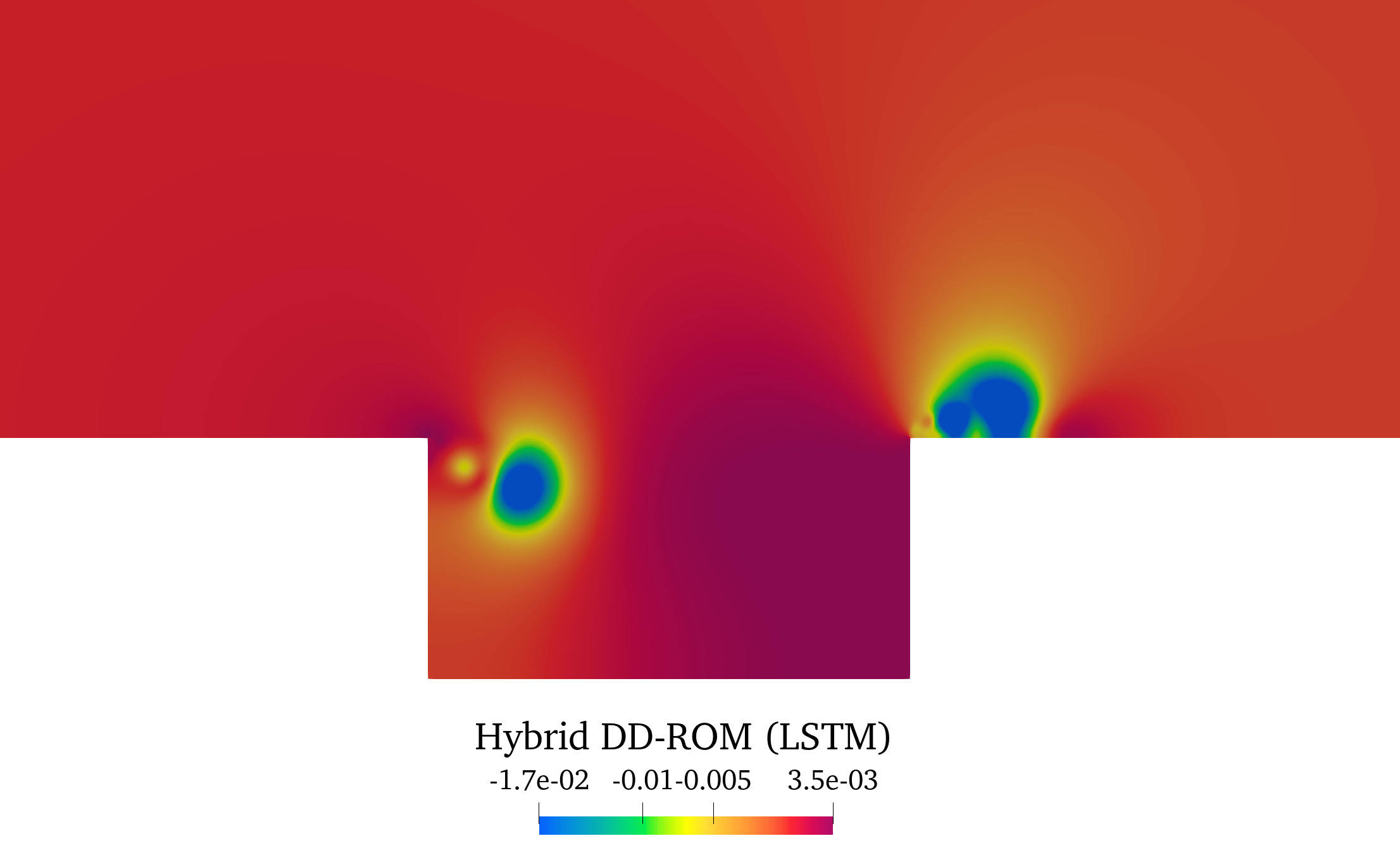}}
    \vspace{0.02\textwidth}
    \subfloat[]
    {\includegraphics[width=0.48\textwidth, trim={1.7cm 0cm 1.7cm 2cm}, clip]{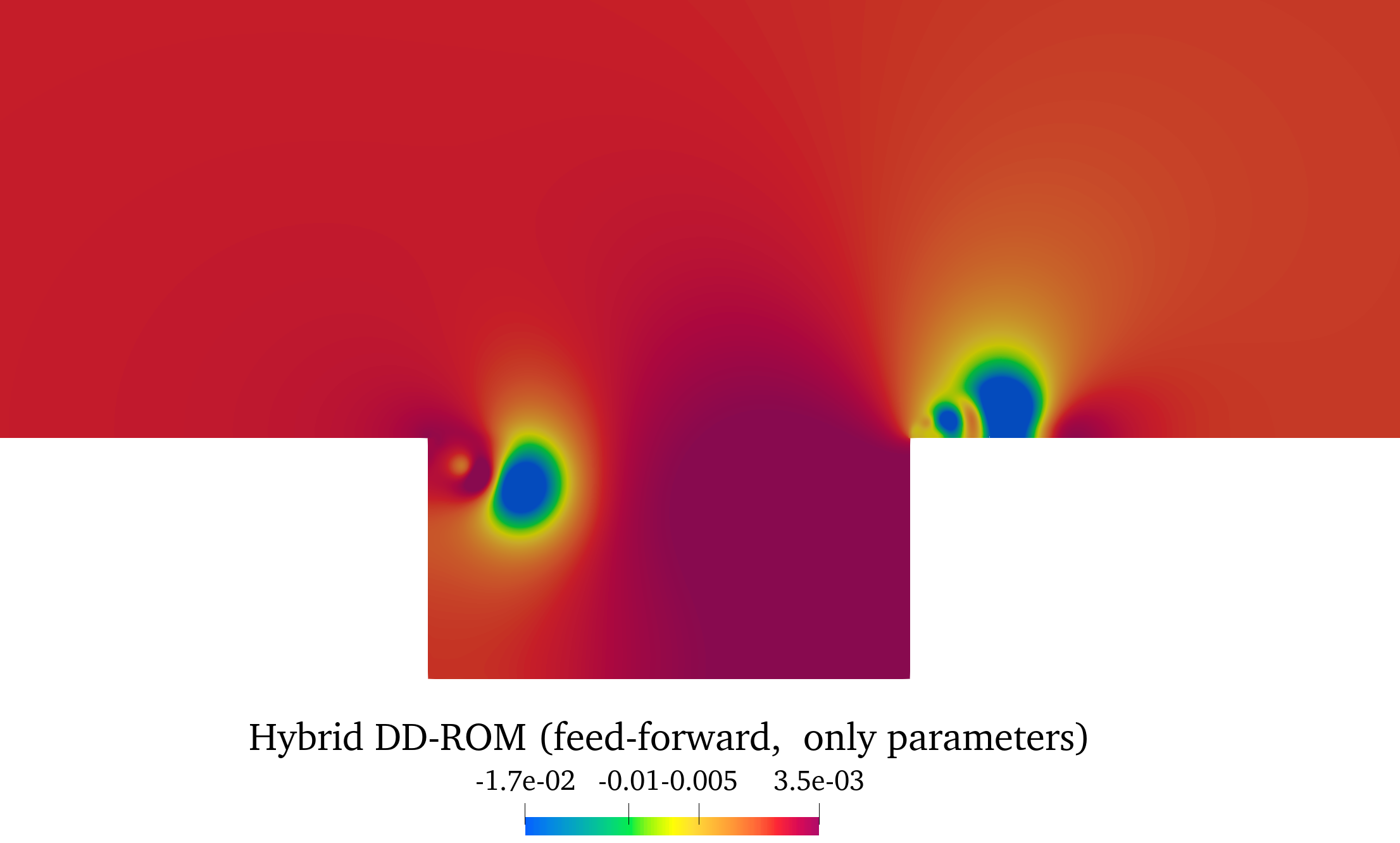}}
    \caption{Pressure fields at final instant of ROM. The standard ROM results are not represented here because the standard simulation blows up.}
    \label{fig:pres-cav-graph}
\end{figure}

The graphical results for the velocity and pressure fields are represented in Figure \ref{fig:vel-cav-graph} and \ref{fig:pres-cav-graph}, respectively. The results are for the final time instance of the online simulation, namely $t=\SI{20}{\second}$.

We exclude the standard ROM from the graphical comparison since it blows up after a few time steps of the online ROM.

As already highlighted in the previous Subsection, Figures \ref{fig:vel-cav-graph} and \ref{fig:pres-cav-graph} show that the DD-ROM approximations are not as accurate as in the cylinder test case.

\subsection{Discussion}
\label{subsec:discussion}

This part of the manuscript is dedicated to a final comparison between the two test cases, for all the different ROMs proposed.

\begin{itemize}
    \item \textbf{POD performance}: In general, as already pointed out, test case \textbf{b} is the most challenging. Indeed, the POD eigenvalues decay is smaller (Figures \ref{fig:eig-a}, \ref{fig:eig-b}) and the reconstruction errors are much higher, as can be seen for instance from the time integrals of the projection errors in Figures \ref{fig:integrals-cylinder} and \ref{fig:integrals-cavity}.
    \item \textbf{Standard ROM}: On the one hand, the standard ROM leads to relative errors increasing in time in the cylinder test case, as shown in Figure \ref{fig:err-turb-test-int}.
    On the other hand, in the cavity test case, the standard POD-Galerkin approach blows up after a few time steps. In this case, the POD-based operators may lead to a ill-conditioned reduced system.
    \item \textbf{Purely DD-ROM}: We only evaluate the performance of this model in test case \textbf{a} and the results slightly improve the standard ROM accuracy. Indeed, if we compare Figures \ref{fig:err-noturb-test-int} (with the purely DD-ROM results) and \ref{fig:err-turb-test-int} (with the hybrid DD-ROM results), we notice that the effect of the correction terms acts in a more massive way on the system if combined with the physics-based DD-ROM. This aspect is deeply analysed in the non-parametric case in \cite{ivagnes2023hybrid}, where the authors showed that the combination of the data-driven methods provides always better results than the individual approaches. This is the reason why we decide to focus only on the hybrid approach in test case \textbf{b}.
    \item \textbf{Hybrid DD-ROM}:
    \begin{itemize}
        \item[$\circ$] As a general consideration, the main difference among the two test cases, is that the DD-ROM acts as a \emph{results-enhancer} for the cylinder, but also as a \emph{system stabilizer} for the cavity case. 
        \item[$\circ$] In the first test case, the LSTM is the best-suited architecture to capture the dynamics of the correction terms. This can be explained with the fact that the LSTM inputs and outputs, namely the velocity/pressure coefficients and the correction coefficients, have the same periodic behavior. Considering a time-memory of about half-period, the network is also well performing in time extrapolation.
        However, when we consider the cavity test case, the best-suitable architecture is the multi-layer perceptron not depending on the pressure coefficients. This is because the pressure coefficients and the correction coefficients have a different behaviour.
        \item[$\circ$] It is important to highlight that all the neural networks taken into account in both the test cases have a good performance in a predictive setting (Figures \ref{fig:cyl-coeffs-turb-0}, \ref{fig:cav-coeffs-turb-0}), but the performance varies when the input coefficients are not anymore the \emph{optimal} ones, namely the projected reduced coefficients. When the inputs are the coefficients $\bm{a}$ and/or $\bm{b}$ coming from the resolution of the dynamical system \eqref{eq:ppe-rom}, the accuracy of the prediction is quite different and depends on the intrinsic relations between inputs and outputs.
    \end{itemize}
    
\end{itemize}

\section{Conclusions and Outlook}
\label{sec:conclusions}
The project presented in this paper aims to enhance the classical ROM approaches using machine learning tools. The general paradigm of DD-ROMs was already presented in \cite{ivagnes2023hybrid}, but it is here extended to a more general parametric setup.

Section \ref{sec:intro} introduces the problem and the issues of standard ROMs in capturing the evolution of the system's dynamics.

Section \ref{sec:methods} is dedicated to the presentation of the methodologies used in this work. In particular, we introduce the FOM used to collect the snapshots in \ref{sec:fom}, the POD-Galerkin ROM approach in \ref{sec:pod-g-roms}, and the machine-learning enhanced ROMs in \ref{sec:ml-roms}. In the last-mentioned part, we consider two types of enhancement: the \emph{physics-based} DD-ROM, aimed to re-introduce the turbulence modeling at the ROM level, is briefly recalled in Subsection \ref{subsec:met-eddy-viscosity}, whereas the \emph{purely} DD-ROM, aimed to close the system through the modeling of the neglected modes, is presented in Subsection \ref{subsec:met-corrections}.

Finally, the numerical results are presented in Section \ref{sec:results}, showing the effects of the machine-learning approaches on the periodic flow past a cylinder (\ref{subsec:test-case-a}), and on the channel-driven cavity flow (\ref{subsec:test-case-b}).

For each test case, we perform the offline stage in a parametric setting, considering as parameters time and the Reynolds number.
After that, a POD eigenvalues analysis is performed (in \ref{subsubsec:fom-a} and \ref{subsubsec:dd-roms-b}), after which we focus our study on the \emph{marginally-resolved} regime. In this framework, standard ROMs produce poor results and do not provide accurate reconstructions of the solutions.

For each test case, we investigate the effects of the machine learning approaches on the reduced dynamical system. In particular, we compared the results of three different types of neural network architectures in both test cases, spacing among LSTM architectures, \emph{ad-hoc} architectures, and more classical multi-layer perceptron structures.

The results are presented in terms of:
(i) neural networks performance directly in a predictive setting (in \ref{subsec:networks-a} and \ref{subsec:networks-b});
(ii) velocity and pressure relative errors with respect to the FOM counterpart for unseen configurations (in \ref{subsubsec:dd-roms-a} and \ref{subsubsec:dd-roms-b});
(iii) graphical views of the fields at the ROM final time instance.

The final Subsection \ref{subsec:discussion} makes a comparison between the two test cases in terms of performances.

In conclusion, we were able to provide in this work an algorithm able to enhance the POD-Galerkin results, acting also as a stabilizer in highly unsteady setups.

As pointed out for the cavity test case, the results' accuracy is influenced by the linearity of the POD approach. Hence, the method may be further improved introducing \emph{nonlinearity} in the projection step, by replacing the POD with more advanced techniques, like autoencoders. This will be the focus of the authors' future research work.

Another extension of interest may be to include an \emph{a-posteriori learning} inside the neural networks used to predict the correction terms. This would increase the training computational time, but improving the final accuracy without affecting the online prediction time.

\section*{Acknowledgements}
\label{sec:acknowledgements}
We acknowledge the support by the European Commission H2020 ARIA (Accurate ROMs for Industrial Applications, GA 872442) project, by MIUR (Italian Ministry for Education University
and Research) and by the European Research Council Consolidator Grant Advanced Reduced Order Methods with Applications in Computational Fluid Dynamics-GA 681447, H2020-ERC COG 2015 AROMA-CFD. 
This work has been conducted within the research activities of the consortium iNEST (Interconnected North-East Innovation Ecosystem), Piano Nazionale di Ripresa e Resilienza (PNRR) – Missione 4 Componente 2, Investimento 1.5 – D.D. 1058 23/06/2022, ECS00000043, supported by the European Union's NextGenerationEU program.
We also acknowledge the support by INdAM-GNCS: Istituto Nazionale di Alta Matematica –– Gruppo Nazionale di Calcolo Scientifico.

Giovanni Stabile acknowledges the financial support under the National Recovery and Resilience Plan (NRRP), Mission 4, Component 2, Investment 1.1, Call for tender No. 1409 published on 14.9.2022 by the Italian Ministry of University and Research (MUR), funded by the European Union – NextGenerationEU– Project Title ROMEU – CUP P2022FEZS3 - Grant Assignment Decree No. 1379 adopted on 01/09/2023 by the Italian Ministry of Ministry of University and Research (MUR).

The main computations in this work were carried out by the usage of ITHACA-FV \cite{ithacasite}, an open-source library and an implementation in OpenFOAM \cite{ofsite} for reduced order modeling techniques. Its developers and contributors are acknowledged.

\newpage

\bibliographystyle{abbrv}
\bibliography{main}

\appendix

\section{Neural networks hyperparameters}
\label{sec:appendix}

In this Section we specify the hyperparameters of the machine learning techniques used in this work.

\subsection{Neural networks' architecture for the hybrid DD-ROM in test case \textbf{a}}
\label{appendix-correction-case-a}
This Section is dedicated to the structure of the neural networks considered in test case \textbf{a} to model the mappings $\mathcal{G}$ (used to compute the eddy-viscosity coefficients), and $\mathcal{M}$ (used to compute the correction coefficients).

Table \ref{tab:networks_first} shows the hyperparameters' settings for the different cases considered. 
In the Table we reported the structure of the hydden layers as a list of numbers, where the length of the list indicates the number of hidden layers and the values in the list indicate the number of neurons for each layer.
In all the models proposed, the final layer is composed of only linear operations.
In particular, for the LSTM network, we considered a sequence length, $N_{\text{seq}}$ in the Table, which is $N_{\text{seq}}\sim T/2$, where $T$ is the flow period.

\begin{figure}[htpb!]
    \centering
    \includegraphics[width=0.9\textwidth, trim={0 16cm 0 0}, clip]{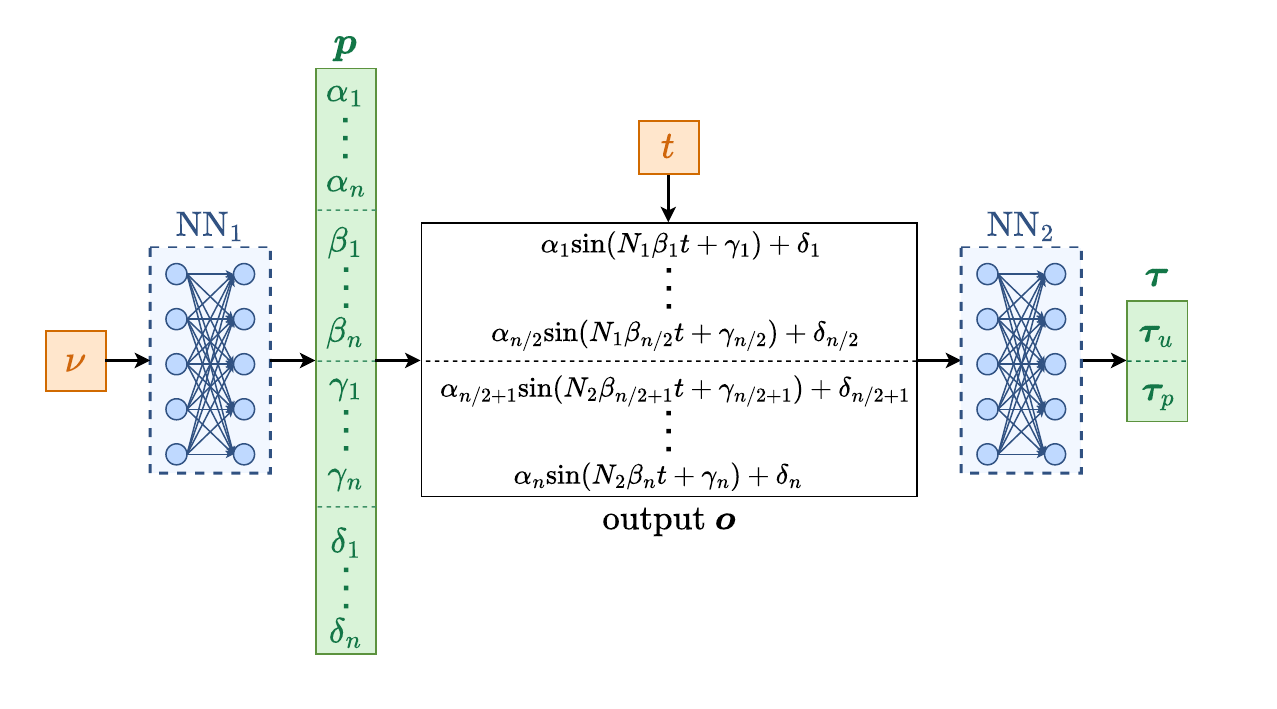}
    \caption{Schematic architecture of the SinNN neural network considered to find the correction coefficients in test case \textbf{a}. The orange boxes are the inputs needed to train the SinNN as a function of the parameters only, namely $\mathcal{M}(\nu, t)$. NN$_1$ and NN$_2$ are two feed-forward fully connected neural networks embedded in the proposed structure and described in Table \ref{tab:networks_first}.}
    \label{fig:sinNN}
\end{figure}

Moreover, for what concerns the SinNN architecture, the values $N_1$ and $N_2$ in \ref{fig:sinNN} have the aim of amplifying the frequency of the periodic expressions proposed, while $n$ is the number of sinusoids considered in the final combination. In our case, $N_1=50$, $N_2=1000$, and $n=20$.

Moreover, in every mapping we consider a regularization in the loss, namely an additional term $\mathcal{L}_{\text{add}}$ depending on the weights' $L^2$ norm. The additional term is expressed as:
\begin{equation}
\mathcal{L}_{\text{add}}=\omega \|\bm{\theta}\|^2_{L^2},
    \label{eq:reg-networks}
\end{equation}
where $\bm{\theta}$ includes all the weights of the neural network, and $\omega=\num{1e-6}$ in our setting.

The optimization employs the \verb|Adam| algorithm in all the networks.

\begin{table}[htpb!]
    \centering
 \caption{Neural networks setting for test case \textbf{a}.}
    \label{tab:networks_first}
    \begin{tabular}{>{\centering\arraybackslash}p{0.18\linewidth}
    >{\centering\arraybackslash}p{0.15\linewidth}
    >{\centering\arraybackslash}p{0.13\linewidth}
    >{\centering\arraybackslash}p{0.15\linewidth}
    >{\centering\arraybackslash}p{0.2\linewidth}
    >{\centering\arraybackslash}p{0.1\linewidth}
    }
    \toprule
    {\textbf{\emph{Mapping}}} &{\textbf{\emph{Network}}}& \textbf{\emph{Hidden layers}} & \textbf{\emph{Non-linearity}} &\textbf{\emph{Learning rate}} & \textbf{\emph{Stopping epoch}}\\
    \midrule
$\bm{g}=\mathcal{G}(\bm{a}, \nu, t)$& Feed-forward &{$[20,20,20]$}&{Softplus}& {$\num{1e-3}$} & {10000}\\ 
   \midrule
   $\bm{\tau}=\mathcal{M}(\bm{a}, \bm{b}, \nu, t)$& Feed-forward &{$[20,20,20]$}&{ReLU}& {$\num{1e-3}$} & {6000}\\
   \midrule
   $\bm{\tau}=\mathcal{M}(\bm{a}, \bm{b}, \nu)$& LSTM &{$[20,20]$, $N_{\text{seq}}$:$100$}&\xmark& {$\num{1e-3}$} & {6000}\\ 
   \midrule
   $\bm{p}=\text{NN}_1(\nu)$& Feed-forward &{$[5, 5]$}&Leaky ReLU& \multirow{2}{*}{$\num{1e-3}$} & \multirow{2}{*}{$6000$}\\ 
   \cline{1-4}
   $\bm{\tau}=\text{NN}_2(\bm{o})$& Feed-forward&{$[20,20]$}&Leaky ReLU& & \\ 
    \bottomrule
    \end{tabular}
\end{table}

\subsection{Neural networks' architecture for the hybrid DD-ROM in test case \textbf{b}}
\label{appendix-correction-case-b}
This Section is dedicated to the specification of the hyperparameters in test case \textbf{b}.
In particular, Table \ref{tab:networks_second} reports the specifics of all the architectures proposed.
Also in this case, we apply a weight regularization, considering the same expression in \eqref{eq:reg-networks}, and with the same weight decay $\omega=\num{1e-6}$. Moreover, as can be seen from Table \ref{tab:networks_second}, the learning rate is decreased with a rate of $1/5$ every $3000$ epochs, in order to avoid picks of the loss during the training process.

\begin{table}[htpb!]
    \centering
 \caption{Neural networks setting for test case \textbf{b}. In the learning rate of the mapping $\mathcal{M}$, we call $t_{\text{epoch}}$ the current epoch, $N_{\text{step}}=3000$ the epochs' interval to decrease the learning rate, $\gamma=0.2$ the corresponding decreasing rate.}
    \label{tab:networks_second}
    \begin{tabular}{>{\centering\arraybackslash}p{0.18\linewidth}
    >{\centering\arraybackslash}p{0.15\linewidth}
    >{\centering\arraybackslash}p{0.13\linewidth}
    >{\centering\arraybackslash}p{0.15\linewidth}
    >{\centering\arraybackslash}p{0.2\linewidth}
    >{\centering\arraybackslash}p{0.1\linewidth}
    }
    \toprule
    {\textbf{\emph{Mapping}}} &{\textbf{\emph{Network}}}& \textbf{\emph{Hidden layers}} & \textbf{\emph{Non-linearity}} &\textbf{\emph{Learning rate}} & \textbf{\emph{Stopping epoch}}\\
    \midrule
$\bm{g}=\mathcal{G}(\bm{a}, \nu, t)$& Standard feed-forward &{$[20,20,20]$}&{Softplus}& {$\num{1e-3}$} & {10000}\\ 
   \midrule
   $\bm{\tau}_p=\mathcal{M}(\bm{b}, \nu, t)$& Feed-forward &{$[20,20,20]$}&{Softplus}& {$(\num{1e-3})\gamma^{\lfloor \frac{t_{\text{epoch}}}{N_{\text{step}}} \rfloor}$} & {10000}\\
   \midrule
   $\bm{\tau}_p=\mathcal{M}(\bm{b}, \nu)$& LSTM &{$[20,20]$, $N_{\text{seq}}$:$20$}&\xmark& {$(\num{1e-3})\gamma^{\lfloor \frac{t_{\text{epoch}}}{N_{\text{step}}} \rfloor}$} & {10000}\\ 
   \midrule
   $\bm{\tau}_p=\mathcal{M}(\nu, t)$&  Feed-forward &{$[5, 5]$}&Softplus& {$(\num{1e-3})\gamma^{\lfloor \frac{t_{\text{epoch}}}{N_{\text{step}}} \rfloor}$} & {$10000$}\\ 
    \bottomrule
    \end{tabular}
\end{table}

We remark that in LSTM networks in both Tables \ref{tab:networks_first} and \ref{tab:networks_second}, no activation function is considered.

The optimization employs the \verb|Adam| algorithm, as in test case \textbf{a}.

\label{appendix-nomenclature}
\printnomenclature

\end{document}